\documentclass[12pt,oneside,reqno]{amsart}
\usepackage{amssymb}
\usepackage{amsmath}
\usepackage{graphicx}
\usepackage{mathrsfs}
\usepackage{stmaryrd}
\usepackage{amsfonts}
\usepackage{amsfonts}
\usepackage{enumerate,amsmath,amssymb,amsthm}
\numberwithin{equation}{section}

\newcommand{\be}{\begin{eqnarray}}
\newcommand{\ee}{\end{eqnarray}}
\newcommand{\ce}{\begin{eqnarray*}}
\newcommand{\de}{\end{eqnarray*}}
\newtheorem{theorem}{Theorem}[section]
\newtheorem{lemma}[theorem]{Lemma}
\newtheorem{remark}[theorem]{Remark}
\newtheorem{definition}[theorem]{Definition}
\newtheorem{proposition}[theorem]{Proposition}
\newtheorem{Examples}[theorem]{Example}
\newtheorem{corollary}[theorem]{Corollary}
\def\eps{\epsilon}

\def\a{\alpha}

\def\b{\beta}
\def\p{\partial}
\def\d{\delta}
\def\g{\gamma}
\def\l{\lambda}

\def\[{{\Big[}}
\def\]{{\Big]}}
\def\<{{\langle}}
\def\>{{\rangle}}
\def\({{\Big(}}
\def\){{\Big)}}

\def\bx{{\mathbf{x}}}

\def\W{{\bf W}}
\def\Ric{{\rm Ricci}}

\def\dif{{\mathord{{\rm d}}}}

\def\div{{\mathord{{\rm div}}}}

\def\u{\mathord{{\bf u}}}
\def\f{\mathord{{\bf f}}}
\def\v{\mathord{{\bf v}}}

\def\b{\mathord{{\bf g}}}

\def\no{\nonumber}
\def\bt{\begin{theorem}}
\def\et{\end{theorem}}
\def\bl{\begin{lemma}}
\def\el{\end{lemma}}
\def\br{\begin{remark}}
\def\er{\end{remark}}
\def\bx{\begin{Examples}}
\def\ex{\end{Examples}}
\def\bd{\begin{definition}}
\def\ed{\end{definition}}
\def\bp{\begin{proposition}}
\def\ep{\end{proposition}}
\def\bc{\begin{corollary}}
\def\ec{\end{corollary}}
\def\cA{{\mathcal A}}
\def\cB{{\mathcal B}}

\def\cE{{\mathcal E}}
\def\cF{{\mathcal F}}

\def\cI{{\mathcal I}}
\def\cJ{{\mathcal J}}

\def\cM{{\mathcal M}}

\def\cO{{\mathcal O}}

\def\mB{{\mathbb B}}
\def\mC{{\mathbb C}}
\def\mD{{\mathbb D}}
\def\mE{{\mathbb E}}

\def\mH{{\mathbb H}}

\def\mN{{\mathbb N}}

\def\mR{{\mathbb R}}
\def\mS{{\mathbb S}}
\def\mT{{\mathbb T}}
\def\mU{{\mathbb U}}

\def\mX{{\mathbb X}}

\def\sA{{\mathscr A}}

\def\sD{{\mathscr D}}

\def\sK{{\mathscr K}}

\def\sP{{\mathscr P}}

\def\geq{\geqslant}
\def\leq{\leqslant}

\def\bH{{\mathbf H}}

\def\fL{{\mathfrak L}}
\def\fS{{\mathfrak S}}
\def\fT{{\mathfrak T}}

\def\fg{{\frak g}}

\def\L{\mathord{{\bf L}}}
\def\C{\mathord{{\bf C}}}

\allowdisplaybreaks

\begin{document}
\title{Stochastic Volterra Equations in Banach Spaces
and Stochastic Partial Differential Equations*}
\date{}
\author{Xicheng ZHANG}




\date{}
\dedicatory{
School of Mathematics and Statistics\\
The University of New South Wales, Sydney, 2052, Australia,\\
Department of Mathematics,
Huazhong University of Science and Technology\\
Wuhan, Hubei 430074, P.R.China\\
email: XichengZhang@gmail.com.}

\begin{abstract}
In this paper, we first study the existence-uniqueness and large deviation estimate
of solutions for stochastic Volterra integral equations
with singular kernels in $2$-smooth Banach spaces. Then, we apply them to a large class of
semilinear stochastic partial differential equations (SPDE) driven by Brownian motions
as well as by fractional Brownian motions, and obtain the existence of
unique maximal strong solutions (in the sense of SDE and PDE) under local Lipschitz conditions.
Lastly, high order SPDEs in a bounded domain of Euclidean space, second order SPDEs on
complete Riemannian manifolds, as well as stochastic Navier-Stokes equations
are investigated.

\end{abstract}

\maketitle

\allowdisplaybreaks
\tableofcontents

\section{Introduction}

The aims of this paper are three folds: First of all, we prove the existence and
uniqueness of solutions with continuous paths for stochastic Volterra integral equations
with singular kernels in a $2$-smooth Banach space. Secondly,
the large deviation principles (abbrev. LDP) of  Freidlin-Wentzell type for
stochastic Volterra equations are established under small perturbations of multiplicative noises.
Thirdly, we apply them to several classes of semilinear stochastic partial
differential equations (abbrev. SPDE). In particular, we give a unified treatment
in certain sense for the LDPs of a large class of SPDEs.

In finite dimensional space, stochastic Volterra integral equations with regular kernels
and driven by Brownian motions were first studied by Berger and Mizel in \cite{bm}.
Later, Protter \cite{Prt} studied the stochastic Volterra  equations driven
by general semimartingales. Using the Skorohod integral, Pardoux and Protter \cite{Pa-Pr}
also investigated the stochastic Volterra  equations with anticipating coefficients.
The study of stochastic Volterra equations with singular kernels can be found in
\cite[etc.]{Co-De,De,wz,La, Nu-Ro}. Recently, the present author  \cite{Zh3} studied the
approximation of Euler type and the LDP of Freidlin-Wentzell type
for stochastic Volterra equations with singular kernels.
In particular, the kernels in \cite{Zh3} may deal with the fractional
Brownian motion kernels as well as the fractional order integral kernels.
The study of LDP for stochastic Volterra equations is also referred to \cite{Nu-Ro, La}.

Since the work of Freidlin and Wentzell \cite{Fr-We}, the theory of
small perturbation large deviations for stochastic differential equations (abbrev. SDE)
has been studied extensively  (cf. \cite[etc.]{Az,Stro}).
In the classical method, to establish such an LDP
for SDE, one usually needs to discretize the time variable and then
prove various necessary exponential continuity and tightness
for approximation equations in different spaces by using comparison principle.
However, such verifications would become rather complicated and even impossible
in some cases, e.g., stochastic evolution equations with multiplicative noises.

Recently, Dupuis and Ellis \cite{de} systematically developed
a weak convergence approach to the theory of large deviation.
The central idea is to prove some variational representation formula
for the Laplace transform of bounded continuous functionals,
which will lead to proving a Laplace principle which is equivalent to
the LDP. In particular, for Brownian functionals, an elegant variational representation
formula has been established by Bou\'e-Dupuis \cite{bd}
and Budhiraja-Dupuis \cite{bd0}. A simplified proof was given by the present author
\cite{Zh2}. This variational representation has already been proved to be very effective
for various finite and infinite dimensional stochastic dynamical systems even
with irregular coefficients  (cf. \cite[etc.]{Rz1,Rz2, Bu-Du-Ma, Zh3,Rz3}).
One of the main advantages of this argument is that one only needs to
make some simple moment estimates (see Section 4 below).

On the other hand, it is well known that in the deterministic case,
many PDE problems of parabolic and hyperbolic types
can be written as Volterra type integral equations in Banach spaces
by using the corresponding semigroup and the variation-of-constants formula  (cf. \cite{Fr, He, Pa}).
An obvious merit of this procedure is that the unbounded operators in PDEs
no longer appear and the analysis is entirely analogous to the ODE case.
Thus, one naturally expects to take the same advantages  for SPDEs in Banach spaces.
However, it is not all Banach spaces in which stochastic integrals are well defined. One can only
work in a class of $2$-smooth Banach spaces. The definition of
stochastic integrals in $2$-smooth Banach spaces
and related properties such as Burkholder-Davis-Gundy's (abbrev. BDG) inequality,
Girsanov's theorem, stochastic Fubini's theorem and the distribution of stochastic integrals
can be found in \cite[etc. ]{Ne,b1,b2,On}. Thus, similar to the deterministic case,
we can develop a parallel theory in $2$-smooth Banach spaces for SPDEs.
It should be emphasized that besides the usual
SPDEs driven by multiplicative Brownian noises,
a class of stochastic evolutionary integral equations appearing in viscoelasticity
and heat conduction with memory  (cf. \cite{Pr}) as well as a class of SPDEs driven by additive
fractional Brownian noise,  can also be written as abstract stochastic Volterra equations
in Banach spaces.

In the past three decades, the theory of general SPDEs has been developed extensively
by numerous authors mainly  based  on two different approaches: semigroup method
based on the variation-of-constants formula (as said above)
(cf. \cite[etc.]{Wa, DaZa, b1,b2,bp,Zh1})
and variation method based on Galerkin's finite dimensional approximation
(cf. \cite[etc.]{Pa2,Kr-Ro,Ro, Kr1, Mi-Ro, Pr-Ro, Zh4, Go-Ro-Zh}).
A new regularization method is given in \cite{Zh7}.
An overview for the classification
and applications of SPDEs are referred to the recent book of Kotelenez \cite{Ko}.
In the author's knowledge, most of the well known results are primarily concentrated on the mild
or weak solutions, even measure-valued solutions.
Such notions of solutions naturally appear in the study of SPDEs driven by the
space-time white noises, and in this case one cannot obtain any differentiability of the solutions
in the spatial variable.

Nevertheless, when one considers an SPDE driven by
the spatial regular and time white  noises, it is reasonable to require the
existence of spatial regular solutions or classical solutions in the sense of PDE.
For  linear SPDEs, such regular solutions are easy and well known
 (cf. \cite[etc.]{Kr-Ro,Ro,Fl}). However, for  nonlinear SPDEs, there seems to be few results
 (cf. \cite{Kr1, Mi-Ro1, Zh2, Zh7}). A major difficulty to prove the spatial regularity
of solutions is that one cannot use the usual bootstrap method
in the theory of PDE since there is no any differentiability of solutions
with respect to the time variable. The present author \cite{Zh2} (see also \cite{Go-Zh,Zh7})
solves this problem by using a non-linear interpolation result due to Tartar \cite{Ta}.
Obviously, for the regularity theory of SPDEs, by using  Sobolev's embedding theorem  (cf. \cite{Ad}),
it is natural to consider the $L^p$-solution of SPDEs. This is also why we need to work in
 $2$-smooth Banach spaces. It should be remarked that
the $L^p$-theory for SPDEs has been established in \cite[etc.]{b1,b2,bp,Kr1,DeSt,DeMaSt,Zh1}.
But, there are few results to deal with the $L^p$-strong solution in the sense of PDE.
In the present paper, we shall prove a general result about the existence of strong solutions
in the sense of both SDE and PDE (see Theorem \ref{Th55}).

We now describe our structure of this paper: In Section 2, we prepare some preliminaries
for later use, and divide it into four subsections. In Subsection 2.1, we prove a
Gronwall's lemma of Volterra type under rather weak assumptions on  kernel functions.
Moreover, two simple examples are provided to show this lemma.
In Subsection 2.2, we recall the It\^o integral in $2$-smooth Banach spaces and
Burkholder-Davies-Gundy's inequality as well as Kolmogorov's continuity criterion
of random fields in random intervals.
In Subsection 2.3, we recall the properties of analytic semigroups
and prove a local non-linear interpolation lemma, see also \cite{Ta} for  other
related non-linear interpolation results. This lemma will play an important role in
proving the existence of strong solutions (in the PDE's sense) in Theorems \ref{62}
and \ref{Th7} below. In Subsection 2.4, we recall the criterion of Laplace principle established by
Budihiraja and Dupuis  \cite{bd,bd0} (see also \cite{Zh5}).

In Section 3, using the Gronwall  inequality of Volterra type in Subsection 2.1,
we first prove the existence and uniqueness of solutions for stochastic Volterra equations
in $2$-smooth Banach spaces under global Lipschitz conditions and singular kernels.
Next, in Subsection 3.2, we study the regularity of solutions under slightly stronger
assumptions on kernels. Moreover, a BDG type of inequality for stochastic Volterra type integral
is also proved. In Subsection 3.3, employing the usual localizing method, we prove
the existence of a unique maximal solution for stochastic Volterra equation under
local Lipschitz conditions. Lastly, in Subsection 3.4, we discuss the
continuous dependence of solutions with respect to the coefficients.

In Section 4, using the weak convergence method, we prove the Freidlin-Wentzell
large deviation principle for the small perturbations of stochastic Volterra equations
under a compactness assumption and some uniform non-explosion
conditions for the controlled equations. We also refer to \cite{Liu, Rz3}
for the application of weak convergence approach in
the LDPs of stochastic evolution equations (the case of evolution triple).
In the proof of Section 4,
we need to use the Yamada-Watanabe Theorem in infinite dimensional space,
which has been established by Ondrej\'at \cite{On} (see also
\cite{Ro-Sc-Zh} for the case of evolution triple). We want to say that although  Ondrej\'at
only considered the case of convolution semigroup, their proofs are also adapted to
more general stochastic Volterra equations. Moreover, since we are considering
the path continuous solution, the proof in \cite{On} can be simplified .

In Section 5, a simple application in a class of semilinear stochastic
evolutionary integral equations is presented, which has been studied in \cite[etc.]{Cl-Da,Bo-St,Ka-Li}
for additive noises. Such type of stochastic evolution equations
appears in viscoelasticity, heat conduction in materials with memory, and electrodynamics with
memory \cite{Pr}.

In Section 6, we apply our general results to a large class of semilinear
stochastic evolution equations driven by multiplicative Brownian noise
and additive fractional Brownian noise. A basic result in semigroup theory states
that if $f$ is a H\"older continuous function in the Banach space $\mX$, then
$$
t\mapsto\int^t_0\fT_{t-s}f(s)\dif s\mbox{ is continuous in $\sD(\fL)$},
$$
where $\fT_t$ is an analytic semigroup and $\fL$ is the generator of $\fT_t$.
We will use this result to prove the existence of strong solutions (in the sense of PDE)
for semilinear SPDEs. Moreover, we also give a simple result about the SPDE driven by
additive fractional Brownian noises. The corresponding LDPs
are also obtained (see also \cite[etc.]{So,Pe,Bu-Du-Ma,Rz3,Liu}
for the study of LDPs of stochastic evolution equations).
We remark that the skeleton equation for the LDP of SPDEs
driven by fractional Brownian motion is a non-convolution
type of Volterra integral equation.

In Section 7, high order SPDEs in a bounded domain of Euclidean space
are studied. Our stochastic version may be regarded as
a parallel result in the deterministic case  (cf. \cite[p.246, Theorem 4.5]{Pa}).
Moreover, the LDP is also obtained.

In Section 8, we in particular study the second order stochastic parabolic equations on
complete Riemannian manifolds. Under one-side Lipschitz and polynomial growth
conditions, we obtain the global existence-uniqueness of strong solutions.
When the manifold is compact, the LDP also holds in this case. In particular,
stochastic reaction diffusion equations with polynomial growth coefficients are included.

In Section 9, we first prove the existence and uniqueness of local
$L^p$-strong solutions for  stochastic Navier-Stokes equations (SNSE) in any dimensional
case. In the two dimensional case, we also obtain the non-explosion of solutions.
Moveover, the LDPs for $2$-dimensional SNSEs are also established in
the case of both Dirichlet boundary and periodic boundary.
We remark that the $L^p$-solutions for SNSEs have been studied by
Brzezniak and Peszat \cite{Br-Pe} (bounded domain) and
Mikulevicius and Rozovskii \cite{Mi-Ro2} (the whole space).
 The large deviation result for two dimensional SNSEs
with additive noise was proved by Chang \cite{Ch} using Girsanov's transformation.
In \cite{Sr-Su}, the authors also used the
weak convergence method to prove the large deviation estimate for two
dimensional SNSEs with multiplicative noises. But, it seems that there is a gap in their proofs
\cite[p.1655 line 6 and p.1658 line 2]{Sr-Su}. Therein, the $\v_n$
only weakly converges to $\v$ in $S_M$. This seems not enough to derive their limits.

\vspace{2mm}

We conclude this introduction by making the following
\textsc{Convention}: Throughout this paper, the letter $C$ with or without subscripts will
denote a positive constant, whose value may change from one place to another. Moreover, we also
use the notation $E_1\preceq E_2$ to denote $E_1\leq C\cdot E_2$, where $C>0$
is an unimportant constant.

\section{Preliminaries}

\subsection{Gronwall's inequality of Volterra type}
Let $\triangle:=\{(t,s)\in\mR^2_+: s\leq t\}$.
We first recall the following result due to Gripenberg \cite[Theorem 1 and p.88]{Gr}.
\bl\label{Le0}
Let $\kappa: \triangle\to\mR_+$ be a measurable function.
Assume that for any $T>0$
$$
t\mapsto\int^t_0\kappa(t,s)\dif s\in L^\infty(0,T)
$$
and
$$
\limsup_{\eps\downarrow 0}\left\|\int^{\cdot+\eps}_{\cdot}
\kappa(\cdot+\eps,s)\dif s\right\|_{L^\infty(0,T)}<1.
$$
Define
\be
r_1(t,s):=\kappa(t,s), \ \ r_{n+1}(t,s):=\int^t_s \kappa(t,u)r_n(u,s)\dif u,\ \ n\in\mN.
\label{QW1}
\ee
Then for any $T>0$, there exist constants $C_T>0$ and $\gamma\in(0,1)$ such that
\be
\left\|\int^\cdot_0r_n(\cdot,s)\dif s\right\|_{L^\infty(0,T)}\leq C_T n\gamma^n,
\ \ \forall n\in\mN.\label{QW2}
\ee
In particular, the series
\be
r(t,s):=\sum_{n=1}^\infty r_n(t,s)\label{R1}
\ee
converges for almost  all $(t,s)\in\triangle$, and
\be
r(t,s)-k(t,s)=\int^t_sk(t,u)r(u,s)\dif u=\int^t_sr(t,u)k(u,s)\dif u\label{R2}
\ee
and for any $T>0$
\be
t\mapsto\int^t_0 r(t,s)\dif s\in L^\infty(0,T).\label{R3}
\ee
\el

The function $r$ defined by (\ref{R1}) is called the resolvent of $\kappa$.
All the functions $\kappa$ in Lemma \ref{Le0} will be denoted by $\sK$.
In what follows, we shall denote by $\sK_0$ the subclass of $\sK$ with the property that
$$
\limsup_{\eps\downarrow 0}\left\|\int^{\cdot+\eps}_{\cdot}
\kappa(\cdot+\eps,s)\dif s\right\|_{L^\infty(0,T)}=0.
$$
We also denote by $\sK_{>1}$ the set of all positive measurable functions $\kappa$ on
$\triangle$ with the property that for any $T>0$ and some $\beta=\beta(T)>1$
\be
t\mapsto\int^t_0\kappa^\beta(t,s)\dif s\in L^\infty(0,T).\label{C2}
\ee
It is clear that $\sK_{>1}\subset\sK_0\subset\sK$  and
for any $\kappa_1,\kappa_2\in\sK_0\ (\mathrm{resp.}\ \sK_{>1})$ and $C_1,C_2\geq0$,
$$
C_1\kappa_1+C_2\kappa_2\in \sK_0\ (\mathrm{resp.}\ \sK_{>1}).
$$

Let $0\leq h\in L^1_{loc}(\mR_+)$. If $\kappa(t,s)=h(s)$, then $\kappa\in\sK_0$ and
$$
r(t,s)=h(s)\exp\left\{\int^t_s h(u)\dif u\right\};
$$
if $\kappa(t,s)=h(t-s)$, then $\kappa\in\sK_0$ and
\be
r(t,s)=a(t-s):=\sum_{n=1}^\infty a_n(t-s),\label{Ah}
\ee
where
$$
a_1(t)=h(t),\ \ a_ {n+1}(t):=\int^t_0 h(t-s)a_n(s)\dif s.
$$
When $0\leq h\in L^1(\mR_+)$, a classical result due to Paley and Wiener
 (cf. \cite[p.207, Theorem 5.2]{Mi}) says that
\be
a\in L^1(\mR_+)\mbox{ if and only if }
\int^\infty_0h(t)\dif t<1.\label{Wi}
\ee
In this case, $\hat a(s)=\hat h(s)/(1-\hat h(s))$, where
the hat denotes the Laplace transform, i.e.:
$$
\hat h(s):=\int^\infty_0 e^{-st}h(t)\dif t, \ \ s\geq 0.
$$
We want to say that (\ref{Wi}) is useful in the study of
large time asymptotic behavior of solutions for
Volterra equations. An important extension to nonintegrable convolution kernel can be found
in \cite{Sh-Wa, Jo-Wh} (see also \cite{Gr}). A simple example is provided in  Example \ref{Exam} below.

We now prove the following Gronwall lemma of Volterra type (see also \cite[Lemma 7.1.1]{He}
for a case of special convolution kernel).
\bl\label{Le1}
Let $\kappa\in\sK$ and $r_n$ and $r$ be defined respectively by (\ref{QW1}) and (\ref{R1}).
Let $f, g:\mR_+\to\mR_+$ be two measurable functions satisfying that
for any $T>0$ and some $n\in\mN$
\be
t\mapsto\int^t_0r_n(t,s) f(s)\dif s\in L^\infty(0,T)\label{QW4}
\ee
and for almost all $t\in(0,\infty)$
\be
\int^t_0r(t,s) g(s)\dif s<+\infty.\label{Es2}
\ee
If for almost all $t\in(0,\infty)$
\be
f(t)\leq g(t)+\int^t_0 \kappa(t,s)f(s)\dif  s,\label{E22}
\ee
then for almost all $t\in(0,\infty)$
\be
f(t)\leq g(t)+\int^t_0r(t,s)g(s)\dif s.\label{Es3}
\ee
\el
\begin{proof}
First of all, if we define
$$
h(t):=g(t)+\int^t_0r(t,s)g(s)\dif s,
$$
then by (\ref{R2}) and (\ref{Es2})
$$
h(t)=g(t)+\int^t_0\kappa(t,s) h(s)\dif s\ \ \mbox{ for a.a. $t\in(0,\infty)$}.
$$
Thus, by (\ref{E22}) we have
\be
f(t)-h(t)\leq\int^t_0\kappa(t,s)(f(s)-h(s))\dif s\ \ \mbox{ for a.a. $t\in(0,\infty)$}.\label{QW3}
\ee

Set $\tilde f(t):=f(t)-h(t)$ and define
$$
\tilde f^*(t):=\mathrm{ess}\!\!\sup_{s\in[0,t]}\tilde f(s),\ \ t>0
$$
and
$$
\tau_0:=\inf\{t>0: \tilde f^*(t)>0\}.
$$
Clearly, $t\mapsto\tilde  f^*(t)$ is increasing and
\be
\tilde f(t)\leq 0\ \ \mbox{ for a.a. $t\in[0,\tau_0)$}.\label{E33}
\ee
We want to prove that
$$
\tau_0=+\infty.
$$

Suppose $\tau_0<+\infty$. Iterating inequality (\ref{QW3}), we obtain
$$
\tilde f(t)\leq\int^t_0r_n(t,s)\tilde f(s)\dif s\leq\int^t_0r_n(t,s) f(s)\dif s,\ \ \forall n\in\mN.
$$
By (\ref{QW4}), one knows that $0<\tilde f^*(t)<+\infty$ for any $t>\tau_0$.
Moreover, we have
$$
\tilde f(t)\stackrel{(\ref{E33})}{\leq}
\int^t_{\tau_0} r_n(t,s)\tilde f(s)\dif  s
\leq \tilde f^*(t)\int^t_{\tau_0} r_n(t,s)\dif  s,\ \ \forall n\in\mN.
$$
So, for any $T>\tau_0$
$$
0<\tilde f^*(T)\leq \tilde f^*(T)\cdot\left\|\int^\cdot_{\tau_0}r_n(\cdot,s)\dif  s
\right\|_{L^\infty(\tau_0,T)}\stackrel{(\ref{QW2})}{\longrightarrow} 0
$$
as $n\to\infty$, which is impossible. So, $\tau_0=+\infty$.
\end{proof}

The following two examples show that (\ref{Es3}) is sensitive to $\kappa\in\sK$.
\bx{\rm For $C_0>0$, set
$$
\kappa_{C_0}(t,s):=\frac{C_0}{\sqrt{t^2-s^2}},\ \  s<t.
$$
It is clear that
$$
\int^t_s\kappa_{C_0}(t,u)\dif u=C_0((\pi/2)-\arcsin(s/t)).
$$
From this, one sees that
\ce
\left\{
\begin{aligned}
&\kappa_{C_0}\notin\sK, \ \ \mbox{ if $C_0\geq 2/\pi$};\\
&\kappa_{C_0}\in\sK\cap\sK^c_0, \ \ \mbox{ if $0<C_0<2/\pi$}.
\end{aligned}
\right.
\de
Consider the following Volterra equation
$$
x(t)=\int^t_0\kappa_{C_0}(t,s)x(s)\dif s,\ \  t\geq 0.
$$
If $C_0=1$, there are at least two solutions $x(t)\equiv0$ and $x(t)=t$;
if $C_0=\frac{2}{\pi}$, there are infinitely many solutions $x(t)\equiv\mathrm{constant}$;
if $0<C_0<2/\pi$, by Lemma \ref{Le1} there is only one solution $x(t)\equiv0$ in $L^\infty_{loc}(\mR_+)$.
}
\ex

\bx{\rm
For $C_0>0$ and $\a,\beta\in[0,1)$, set
$$
\kappa_{C_0}^{\a,\beta}(t,s):=\frac{C_0}{(t-s)^\a s^\beta},\ \  s<t.
$$
It is clear that
\be
\int^t_u\kappa_{C_0}^{\a,\beta}(t,s)\dif s=C_0 t^{1-\a-\beta}\int^1_{u/t}\frac{1}{(1-s)^\a s^\beta}
\dif s.\label{Es6}
\ee
From this, one sees that
\ce
\left\{
\begin{aligned}
&\kappa_{C_0}^{\a,\beta}\notin\sK, \ \ \mbox{ if $\a+\beta>1$ and $C_0>0$};\\
&\kappa_{C_0}^{\a,\beta}\notin\sK, \ \ \mbox{ if $\a+\beta=1$ and
$C_0\geq \int^1_0\frac{1}{(1-s)^\a s^\beta}\dif s$};\\
&\kappa_{C_0}^{\a,\beta}\in\sK\cap\sK^c_0, \ \ \mbox{ if $\a+\beta=1$ and
$C_0<\int^1_0\frac{1}{(1-s)^\a s^\beta}\dif s$};\\
&\kappa_{C_0}^{\a,\beta}\in\sK_{>1}, \ \ \mbox{ if $\a+\beta<1$ and
$C_0>0$}.
\end{aligned}
\right.
\de
Consider the following Volterra equation
\ce
x(t)=\int^t_0k^{\a,\beta}_{C_0}(t,s)x(s)\dif s,\ \  t\geq 0.
\de
If $\a+\beta<1$, by Lemma \ref{Le1} there is only one solution $x(t)\equiv0$ in $L^\infty_{loc}(\mR_+)$;
if $\a=\beta=C_0=1/2$, there are at least two solutions $x(t)\equiv0$ and $x(t)=\sqrt{t}$.
}
\ex
\subsection{Ito's integral in $2$-smooth Banach spaces}
Throughout this paper, we shall fix a stochastic basis
$(\Omega,\cF,P;(\cF_t)_{t\geq 0})$, i.e., a complete probability space with a family of
right-continuous filterations. In what follows, without special declarations, all expectations
$\mE$ are  taken with respect to the probability measure $P$.

Let $\{W^k(t): t\geq 0, k\in\mN\}$ be a sequence of independent one dimensional
standard Brownian motions on $(\Omega,\cF,P;(\cF_t)_{t\geq 0})$.
Let $l^2$ be the usual Hilbert space of all square summable real number sequences,
$\{e_k, k\in\mN\}$ the usual orthonormal basis of $l^2$.
Let $\mX$ be a separable Banach space, and $ L(l^2;\mX)$ the set of all bounded linear
operators from $l^2$ to $\mX$. For an operator $B\in L(l^2;\mX)$, we also write
$$
B=(B_1,B_2,\cdots)\in\mX^\mN,\ \ B_k=Be_k.
$$

\bd\label{Def2}
An operator $B\in L(l^2;\mX)$ is called radonifying if
$$
\mbox{ the series }\sum_{k}B e_k\cdot W^k(1)\ \mbox{ converges in $\ L^2(\Omega;\mX)$}.
$$
We shall denote by $ L_2(l^2;\mX)$ the space of all radonifying operators, and write for
$B\in L_2(l^2;\mX)$
\be
\|B\|_{ L_2(l^2;\mX)}:=\left(\mE\big\| B e_k\cdot W^k(1)\big\|^2_\mX\right)^{1/2}.\label{Nor}
\ee
Here and below, we use the convention that the repeated indices will be summed.
\ed
The following proposition is well known, and a detailed proof was given
in \cite[Proposition 2.5]{On}.
\bp
The space $ L_2(l^2;\mX)$ with norm (\ref{Nor}) is a separable Banach space.
\ep

In order to introduce the stochastic integral of an $\mX$-valued measurable $(\cF_t)$-adapted
process with respect to $W$, in the sequel, we assume that $\mX$ is $2$-smooth  (cf. \cite{Pi}),
i.e., there exists a constant $C_\mX\geq 2$ such that for all $x,y\in\mX$
$$
\|x+y\|^2_\mX+\|x-y\|^2_\mX\leq 2\|x\|^2_\mX+C_\mX\|y\|^2_\mX.
$$

Let now $s\mapsto B(s)$ be an $ L_2(l^2;\mX)$-valued measurable and $(\cF_t)$-adapted
process with
$$
\int^T_0\|B(s)\|^2_{ L_2(l^2;\mX)}\dif s<+\infty,\ \ a.s.,\ \ \forall T>0.
$$
One can define the It\^o stochastic integral  (cf. \cite[Section 3]{On})
$$
t\mapsto \cI_t(B):=\int^t_0 B(s)\dif W(s)=\int^t_0 B_k(s)\cdot\dif W^k(s)\in\mX
$$
such that $t\mapsto\cI_t(B)$ is an $\mX$-valued  continuous local $(\cF_t)$-martingale.
Moreover, let $\tau$ be any ($\cF_t$)-stopping time, then
$$
\int^{t\wedge\tau}_0 B(s)\dif W(s)=\int^{t}_0 1_{\{s<\tau\}}\cdot B(s)\dif W(s).
$$

The following BDG inequality for $\cI_t(B)$ holds  (cf. \cite[Section 5]{On}).
\bt
For any $p>0$, there exists a constant $C_p>0$ depending only on $p$ such that
\be
\mE\left(\sup_{t\in[0,T]}\Big\|\int^t_0 B(s)\dif W(s)\Big\|^p_\mX\right)
\leq C_p\mE\left(\int^T_0\|B(s)\|^2_{ L_2(l^2;\mX)}\dif s\right)^{p/2}.\label{BDG}
\ee
\et

The following two typical examples of $2$-smooth Banach spaces are usually met in applications.
\bx{\rm
Let $\mX$ be a separable Hilbert space. Clearly, $\mX$
is $2$-smooth. In this case, $ L_2(l^2;\mX)$ consists of all Hilbert-Schmidt
operators of mapping $l^2$ into $\mX$, and
\ce
\|B\|_{ L_2(l^2;\mX)}=\left(\sum_{k=1}^\infty\|B e_k\|_\mX^2\right)^{1/2}.
\de
}
\ex
\bx{\rm
Let $(E,\cE,\mu)$ be a measure space, $\mH$ a separable Hilbert space. For $p\geq 2$, let
$L^p(E,\mu;\mH)$ be the usual $\mH$-valued $L^p$-space over $(E,\cE,\mu)$.
Then $\mX=L^p(E,\mu;\mH)$ is $2$-smooth  (cf. \cite{Pi,b1}). In this case,
by BDG's inequality for Hilbert space valued martingale we have
\be
\|B\|^2_{ L_2(l^2;\mX)}&=&\mE\left(\int_E\big\|B_k(x)\cdot
W^k(1)\big\|^p_\mH\mu(\dif x)\right)^{2/p}\no\\
&\leq&\left(\int_E\mE\big\|B_k(x)\cdot W^k(1)\big\|^p_\mH\mu(\dif x)\right)^{2/p}\no\\
&\leq& C_p\left(\int_E\Big(\sum_{k=1}^\infty
\|B_k(x)\|^2_\mH\Big)^{p/2}\mu(\dif x)\right)^{2/p}\no\\
&=&C_p\|B\|^2_{L^p(E,\mu;l^2\otimes\mH)}.\label{Exa}
\ee
Hence
$$
L^p(E,\mu;l^2\otimes\mH)\hookrightarrow L_2(l^2;\mX)=L_2(l^2;L^p(E,\mu;\mH)).
$$
}
\ex

We also recall the following Kolmogorov  continuity criterion, which can be derived
directly by Garsia's inequality  (cf. \cite{Wa}).
\bt\label{KCC}
Let $\{X(t), t\geq 0\}$ be an $\mX$-valued stochastic process, and $\tau$
a bounded random time. Suppose that for some $C_0,p>0$ and $\d>1$
$$
\mE\|(X(t)-X(s))\cdot 1_{\{s,t\in[0,\tau]\}}\|^p_\mX\leq C_0|t-s|^{\d}.
$$
Then there exist constants $C_1>0$ and $a\in(0,(\d-1)/p)$ independent of $C_0$ and a
continuous version $\tilde X$ of $X$ such that
$$
\mE\left(\sup_{s\not=t\in[0,\tau]}\frac{\|\tilde X(t)
-\tilde X(s)\|^p_\mX}{|t-s|^{ap}}\right)\leq C_1\cdot C_0.
$$
\et
\subsection{A local non-linear interpolation lemma}
In what follows, we fix a densely defined closed linear operator $\fL$
on $\mX$ for which
\be
S_\phi:=\{\lambda\in\mC: 0<\phi\leq|\arg\lambda|\leq\pi\}\subset\rho(\fL),\label{Sec}
\ee
and for some $C\geq 1$
$$
\|(\lambda-\fL)^{-1}\|_{L(\mX)}\leq \frac{C}{1+|\l|},\ \ \l\in S_\phi,
$$
where $\rho(\fL)$ denotes the resolvent set of $\fL$. The above operator $\fL$
 is also called sectorial (cf. \cite[p.18]{He}). It is well known that
$\fL$ generates an analytic semigroup
$$
\fT_t=e^{-\fL t},\ \  t\geq 0.
$$
Moreover, we also assume that
$\fL^{-1}$ is a bounded linear operator on $\mX$, i.e.,
$$
0\in\rho(\fL).
$$
Thus,  for any $\a\in\mR$, the fractional power $\fL^{\a}$ is well defined
 (cf. \cite{He,Pa}). For $\a>0$, we define the fractional Sobolev space $\mX_\a$ by
$$
\mX_\a:=\sD(\fL^\a)
$$
with the norm
$$
\|x\|_{\mX_\a}:=\|\fL^\a x\|_\mX.
$$
For $\a<0$,  $\mX_\a$ is defined as the completion of $\mX$ with respect to the above norm.
It is clear that $\mX_\a$ is still $2$-smooth, and
$B\in L_2(l^2;\mX_\a)$ if and only if $\fL^\a B\in L_2(l^2;\mX)$, i.e.,
\be
\|B\|_{L_2(l^2;\mX_\a)}=\|\fL^\a B\|_{L_2(l^2;\mX)}.\label{EP2}
\ee

The following properties are well known  (cf. \cite[p.24-27]{He} or \cite[p.74]{Pa}).
\bp\label{Pr1}
\begin{enumerate}[(i)]
\item $\fT_t: \mX\to\mX_\a$ for each $t>0$ and $\a>0$.

\item For each $t>0$, $\a\in\mR$ and every $x\in\mX_\a$, $\fT_t\fL^\a x=\fL^\a \fT_t x$.

\item For some $\d>0$ and each $t,\a>0$,
the operator $\fL^\a \fT_t$ is bounded in $\mX$ and
$$
\|\fL^\a \fT_tx\|_\mX\leq C_\a t^{-\a} e^{-\d t}\|x\|_\mX, \ \ \forall x\in\mX.
$$

\item Let $\a\in(0,1]$ and $x\in\mX_\a$, then
$$
\|\fT_t x-x\|_\mX\leq C_\a t^\a\|x\|_{\mX_\a}.
$$
\item For any $0\leq \beta<\a$
$$
\|x\|_{\mX_\beta}\leq C_{\a,\beta}\|x\|^{1-\frac{\beta}{\a}}\|x\|^{\frac{\beta}{\a}}_{\mX_\a},
\ \ \forall x\in\mX_\a.
$$
\end{enumerate}
\ep

We need the following embedding result.
\bp
For any $0<\theta<1$ and $\a>0$
\be
(L_2(l^2;\mX),L_2(l^2;\mX_{\a}))_{\theta,1}
\subset L_2(l^2;(\mX,\mX_{\a})_{\theta,1})
\subset L_2(l^2;\mX_{\theta\a}),\label{Inter}
\ee
where $(\cdot,\cdot)_{\theta,1}$ stands for the real interpolation space
between two Banach spaces.
\ep
\begin{proof}
We only prove the first embedding. The second embedding follows
from \cite[p.101, (d) and (f)]{Tr}, i.e.,
$$
(\mX,\mX_{\a})_{\theta,1}\subset\mX_{\theta\a}.
$$

Let
$$
B\in(L_2(l^2;\mX),L_2(l^2;\mX_{\a}))_{\theta,1}=:\mB_{\theta,1}.
$$
By the $K$-method of real interpolation space, we have (cf. \cite[p.24]{Tr})
$$
\|B\|_{\mB_{\theta,1}}=\int^\infty_0\frac{t^{-\theta}K(t,B)}{t}\dif t,
$$
where the $K$-function of $B$ is defined by
$$
K(t,B):=\inf_{B=B_1+B_2}\Big\{\|B_1\|_{L_2(l^2;\mX)}+t\|B_2\|_{L_2(l^2;\mX_{\a})}\Big\},
\ \ t\geq 0.
$$

By Definition \ref{Def2} we have
\ce
K(t,B)&=&\inf_{B=B_1+B_2}\Big\{\Big(\mE\|B_1e_k\cdot W^k(1)\|^2_{\mX}\Big)^{\frac{1}{2}}
+t\Big(\mE\|B_2e_k\cdot W^k(1)\|^2_{\mX_\a}\Big)^{\frac{1}{2}}\Big\}\\
&\geq&\inf_{B=B_1+B_2}\Big\{\Big(\mE\big[\|B_1e_k\cdot W^k(1)\|_{\mX}+
t\|B_2e_k\cdot W^k(1)\|_{\mX_\a}\big]^2\Big)^{\frac{1}{2}}\Big\}\\
&\geq&\Big(\mE\Big[\inf_{B=B_1+B_2}\Big\{\|B_1e_k\cdot W^k(1)\|_{\mX}+
t\|B_2e_k\cdot W^k(1)\|_{\mX_\a}\Big\}\Big]^2\Big)^{\frac{1}{2}}\\
&\geq&\Big(\mE\Big[K\big(t,Be_k\cdot W^k(1)\big)\Big]^2\Big)^{\frac{1}{2}},
\de
where
$$
K\big(t,Be_k\cdot W^k(1)\big):=\inf_{Be_k\cdot W^k(1)=x_1+x_2}\Big\{\|x_1\|_{\mX}+t\|x_2\|_{\mX_\a}\Big\}.
$$

Therefore, by Minkowski's inequality we obtain
\ce
\|B\|_{\mB_{\theta,1}}&\geq&
\left(\mE\left[\int^\infty_0\frac{t^{-\theta}
K\big(t,Be_k\cdot W^k(1)\big)}{t}\dif t\right]^{2}\right)^{\frac{1}{2}}\\
&=&\left(\mE\|Be_k\cdot W^k(1)\|^2_{(\mX,\mX_\a)_{\theta,1}}\right)^{\frac{1}{2}}
=\|B\|_{L_2(l^2;(\mX,\mX_{\a})_{\theta,1})}.
\de
The result follows.
\end{proof}

The following local non-linear interpolation lemma will play a crucial role in
the proofs of Theorems \ref{62} and \ref{Th7} below. We refer to \cite{Ta} for some
other nonlinear interpolation results.
\bl\label{Le7}
Let $0\leq \a_0<\a_1\leq 1$ and $0\leq\a_2<\a_3\leq 1$.
Let $\Psi: \mX_{\a_0}\to L_2(l^2;\mX_{\a_2})$ be a locally Lipschitz
continuous map, and satisfy that for all $R>0$
and $x\in\mX_{\a_1}$ with $\|x\|_{\mX_{\a_0}}\leq R$
$$
\|\Psi(x)\|_{L_2(l^2;\mX_{\a_3})}\leq C_R(1+\|x\|_{\mX_{\a_1}}).
$$
Then for any $0<\theta'<\theta<1$ and $R>0$
$$
\sup_{\|x\|_{\mX_{\a_0+\theta(\a_1-\a_0)}}\leq R}
\|\Psi(x)\|_{L_2(l^2;\mX_{\a_2+\theta'(\a_3-\a_2)})}\leq C_R.
$$
\el
\begin{proof}
By (\ref{EP2}), we may assume that $\a_2=0$.
Fix $R>0$ and $x\in\mX_{\a_0+\theta(\a_1-\a_0)}$ with
$$
\|x\|_{\mX_{\a_0+\theta(\a_1-\a_0)}}\leq R.
$$
Set for $t\geq0$
$$
K(t,\Psi(x)):=\inf_{\Psi(x)=\Psi_1+\Psi_2}
\Big\{\|\Psi_1\|_{L_2(l^2;\mX)}+t\|\Psi_2\|_{L_2(l^2;\mX_{\a_3})}\Big\}.
$$
For $\d>0$ and $t\in[0,1]$, noting that
$$
\|\fT_{t^\d}x\|_{\mX_{\a_0}}\preceq \|x\|_{\mX_{\a_0}}\preceq
\|x\|_{\mX_{\a_0+\theta(\a_1-\a_0)}}\preceq R,
$$
by the assumptions and (iii) and (iv) of Proposition \ref{Pr1}
we have
\ce
K(t,\Psi(x))&\leq& \|\Psi(x)-\Psi(\fT_{t^\d} x)\|_{L_2(l^2;\mX)}
+t\|\Psi(\fT_{t^\d} x)\|_{L_2(l^2;\mX_{\a_3})}\\
&\leq&C_R \|\fT_{t^\d} x-x\|_{\mX_{\a_0}}+C_R t\cdot(1+\|\fT_{t^\d}x\|_{\mX_{\a_1}})\\
&\leq&C_R t^{\d\theta(\a_1-\a_0)}\cdot\|x\|_{\mX_{\a_0+\theta(\a_1-\a_0)}}\\
&&+C_Rt\cdot\big(1+t^{-\d(1-\theta)(\a_1-\a_0)}\|x\|_{\mX_{\a_0+\theta(\a_1-\a_0)}}\big)\\
&\leq& C_R\big(t^{\d\theta(\a_1-\a_0)}+t+t^{1-\d(1-\theta)(\a_1-\a_0)}\big).
\de
Letting $\d=\frac{1}{\a_1-\a_0}$, we obtain that for $t\in[0,1]$
$$
K(t,\Psi(x))\leq C_R(t^\theta+t)\leq C_R t^\theta.
$$
Moveover, it is clear that for $t\geq 1$
$$
K(t,\Psi(x))\leq \|\Psi(x)\|_{L_2(l^2;\mX)}\leq C_R\|x\|_{\a_0}+\|\Psi(0)\|_{L_2(l^2;\mX)}\leq C_R.
$$

Hence, for any $0<\theta'<\theta<1$
\ce
\|\Psi(x)\|_{(L_2(l^2;\mX),L_2(l^2;\mX_{\a_3}))_{\theta',1}}
&=&\int^\infty_0\frac{t^{-\theta'} K(t;\Psi(x)) }{t}\dif t\\
&\leq&C_R\left[\int^1_0\frac{t^{\theta-\theta'}}{t}\dif t
+\int^\infty_1\frac{t^{-\theta'}}{t}\dif t\right]\leq C_R.
\de
The result follows by (\ref{Inter}).
\end{proof}
\subsection{A criterion for Laplace principles}\label{Sec2}
It is well known that there exists a Hilbert space so that $l^2\subset\mU$
is Hilbert-Schmidt with embedding operator $J$ and
$\{W^k(t),k\in\mN\}$ is a  Brownian motion with values in $\mU$, whose covariance operator is
given by $Q=J\circ J^*$.
For example, one can take $\mU$ as the completion of $l^2$ with respect to
the norm generated by scalar product
$$
\<h,h'\>_\mU:=\left(\sum_{k=1}^\infty \frac{h_k h'_k}{k^2}\right)^{\frac{1}{2}}, \ \ h,h'\in l^2.
$$

For $T>0$ and a Banach space $\mB$, we denote by $\cB(\mB)$ the Borel $\sigma$-field, and
by $\mC_T(\mB)$ the continuous function space from $[0,T]$
to $\mB$, which is endowed with the uniform norm.
Define
\be
\ell^2_T:=\left\{h=\int^\cdot_0\dot h(s)\dif s: ~~\dot h\in L^2(0,T;l^2)\right\}\label{H2}
\ee
with the norm
$$
\|h\|_{\ell^2_T}:=\left(\int^T_0\|\dot h(s)\|_{l^2}^2\dif s\right)^{1/2},
$$
where the dot denotes the generalized derivative.
Let $\mu$ be the law of the Brownian motion $W$ in $\mC_T(\mU)$. Then
$$
(\mC_T(\mU),\ell^2_T,\mu)
$$
forms an abstract Wiener space.

For $T,N>0$,  set
$$
\mD_N:=\{h\in \ell^2_T: \|h\|_{\ell^2_T}\leq N\}
$$
and
\be
\cA^T_N:=\left\{
\begin{aligned}
&\mbox{ $h: [0,T]\to l^2$ is a continuous and
$(\cF_t)$-adapted }\\
&\mbox{ process, and for almost all $\omega$},\ \ h(\cdot,\omega)\in\mD_N
\end{aligned}
\right\}.\label{Op2}
\ee
It is well known that with respect to the weak convergence topology in $\ell^2_T$  (cf. \cite{km}),
\be
\mbox{$\mD_N$ is metrizable as a compact Polish space}.\label{Metr}
\ee

Let $\mS$ be a Polish space. A function $I: \mS\to[0,\infty]$ is given.
\bd
The function $I$  is called a rate function if for every $a<\infty$, the set
$\{f\in\mS: I(f)\leq a\}$ is compact in $\mS$.
\ed

Let $\{Z_\eps: \mC_T(\mU)\to\mS,\eps\in(0,1)\}$
be a family of measurable mappings. Assume that
there is a measurable map $Z_0: \ell^2_T\mapsto \mS$ such that
\begin{enumerate}[{\bf (LD)$_\mathbf{1}$}]
\item
For any $N>0$, if a family $\{h^\eps, \eps\in(0,1)\}\subset\cA^T_N$ (as random variables in $\mD_N$)
converges in distribution to  $h\in \cA^T_N$, then
for some subsequence $\eps_k$, $Z_{\eps_k}\Big(\cdot+\frac{h^{\eps_k}(\cdot)}{\sqrt{\eps_k}}\Big)$
converges in distribution to $Z_0(h)$ in $\mS$.
\end{enumerate}
\begin{enumerate}[{\bf (LD)$_\mathbf{2}$}]
\item
 For any $N>0$, if $\{h_n,n\in\mN\}\subset \mD_N$ weakly converges to $h\in\ell^2_T$,
then for some subsequence $h_{n_k}$, $Z_0(h_{n_k})$ converges to $Z_0(h)$ in $\mS$.
\end{enumerate}

For each $f\in\mS$, define
\be
I(f):=\frac{1}{2}\inf_{\{h\in\ell^2_T:~f=Z_0(h)\}}\|h\|^2_{\ell^2_T},\label{ra}
\ee
where $\inf\emptyset=\infty$ by convention. Then under {\bf (LD)$_\mathbf{2}$},
$I(f)$ is a  rate function.
In fact, assume that $I(f_n)\leq a$. By the definition of $I(f_n)$, there exists a
sequence $h_n\in\ell_2$ such that $Z_0(h_n)=f_n$ and
$$
\frac{1}{2}\|h_n\|_{\ell^2_T}^2\leq a+\frac{1}{n}.
$$
By the weak compactness of $\mD_{2a+2}$, there exist a subsequence $n_k$ (still denoted by $n$)
and $h\in\ell^2_T$ such that $h_n$ weakly converges to $h$ and
$$
\|h\|_{\ell^2_T}^2\leq\varliminf_{n\rightarrow\infty}\|h_n\|_{\ell^2_T}^2\leq 2a.
$$
Hence, by {\bf (LD)$_\mathbf{2}$} we have
$$
\lim_{k\to\infty}\|Z_0(h_{n_k})-Z_0(h)\|_\mS=0
$$
and
$$
I(Z_0(h))\leq a.
$$

We recall the following result due to \cite{bd, bd0} (see also \cite[Theorem 4.4]{Zh2}).
\bt\label{Th2}
Under {\bf (LD)$_\mathbf{1}$} and {\bf (LD)$_\mathbf{2}$},
$\{Z_\eps,\eps\in(0,1)\}$ satisfies the Laplace principle with
the rate function $I(f)$ given by (\ref{ra}). More precisely, for each real bounded continuous
function  $g$ on $\mS$:
\be
\lim_{\eps\rightarrow 0}\eps\log\mE^{\mu}\left(\exp\left[-\frac{g(Z_\eps)}{\eps}\right]\right)
=-\inf_{f\in\mS}\{g(f)+I(f)\}.\label{La}
\ee
In particular, the family of $\{Z_\eps,\eps\in(0,1)\}$
satisfies  the large deviation principle in $(\mS,\cB(\mS))$ with the rate function $I(f)$.
More precisely, let $\nu_\eps$ be the law of $Z_\eps$ in $(\mS,\cB(\mS))$,
then for any $A\in\cB(\mS)$
$$
-\inf_{f\in A^o}I(f)\leq\liminf_{\eps\rightarrow 0}\eps\log\nu_\eps(A)
\leq\limsup_{\eps\rightarrow 0}\eps\log\nu_\eps(A)\leq -\inf_{f\in \bar A}I(f),
$$
where the closure and the interior are taken in $\mS$,
and $I(f)$ is defined by (\ref{ra}).
\et

\section{Abstract stochastic Volterra integral equations}

In this section, we consider the following stochastic Volterra integral equation
in a $2$-smooth Banach space $\mX$:
\be
X(t)=g(t)+\int^t_0A(t,s,X(s))\dif s+\int^t_0B(t,s,X(s))\dif W(s),\label{Eq}
\ee
where $g(t)$ is an $\mX$-valued  measurable $(\cF_t)$-adapted process, and
$$
A:\triangle\times\Omega\times\mX\to\mX\in\cM_\triangle\times\cB(\mX)/\cB(\mX)
$$
and
$$
B:\triangle\times\Omega\times\mX\to L_2(l^2;\mX)\in\cM_\triangle\times\cB(\mX)/\cB( L_2(l^2;\mX)).
$$
Here and below, $\triangle:=\{(t,s)\in\mR^2_+: s\leq t\}$, and
$\cM_\triangle$ denotes the progressively measurable $\sigma$-field on $\triangle\times\Omega$
generated by the sets $E\in \cB(\triangle)\times\cF$ with properties:
$1_E(t,s,\cdot)\in\cF_s$ for all $(t,s)\in \triangle$, and
$s\mapsto 1_E(t,s,\omega)$ is right continuous for any $t\in\mR_+$ and $\omega\in\Omega$.

We start with the global existence and uniqueness of solutions of
\textsc{Eq.}(\ref{Eq}) under global Lipschitz conditions and singular kernels.

\subsection{Global existence and uniqueness}
In this subsection, we make the following global Lipschitz and linear growth conditions
on the coefficients:
\begin{enumerate}[{\bf (H1)}]
\item For some $p\geq 2$ and any $T>0$
$$
\mathrm{ess}\!\!\sup_{t\in[0,T]}\int^t_0[\kappa_1(t,s)+\kappa_2(t,s)]
\cdot\mE\|g(s)\|_\mX^p\dif s<+\infty,
$$
where $\kappa_1$ and $\kappa_2$ are from {\bf (H2)} and {\bf (H3)} below.

\item There exists  $\kappa_1\in\sK_0$ such that for all $(t,s)\in\triangle$,
$\omega\in\Omega$ and $x\in\mX$
$$
\|A(t,s,\omega,x)\|_\mX\leq \kappa_1(t,s)\cdot(\|x\|_\mX+1)
$$
and
$$
\|B(t,s,\omega,x)\|^2_{L_2(l^2;\mX)}\leq \kappa_1(t,s)\cdot(\|x\|_\mX^2+1).
$$
\item There exists  $\kappa_2\in\sK_0$ such that for all $(t,s)\in\triangle$,
$\omega\in\Omega$ and $x,y\in\mX$
$$
\|A(t,s,\omega,x)-A(t,s,\omega,y)\|_\mX\leq \kappa_2(t,s)\cdot\|x-y\|_\mX
$$
and
$$
\|B(t,s,\omega,x)-B(t,s,\omega,y)\|^2_{L_2(l^2;\mX)}\leq \kappa_2(t,s)\cdot\|x-y\|_\mX^2.
$$
\end{enumerate}

We now prove the following basic existence and uniqueness result.
\bt\label{Th0}
Assume that {\bf (H1)}-{\bf (H3)} hold. Then there exists a unique measurable
$(\cF_t)$-adapted process $X(t)$ such that for almost all $t\geq0$,
\be
X(t)=g(t)+\int^t_0A(t,s,X(s))\dif s+\int^t_0B(t,s,X(s))\dif W(s),\ \ \mbox{ $P$-a.s.}
\label{PP1}
\ee
and for any $T>0$ and some $C_{T,p,\kappa_1}>0$,
\be
\mE\|X(t)\|^p_\mX\leq C_{T,p,\kappa_1}\left[\mE\|g(t)\|_\mX^p
+\mathrm{ess\!\!}\sup_{t\in[0,T]}\int^t_0\kappa_1(t,s)\cdot\mE\|g(s)\|_\mX^p\dif s\right]\label{PP2}
\ee
for almost all $t\in[0,T]$, where $p$ is from {\bf(H1)}.
Moreover, if
\be
t\mapsto\int^t_0\kappa_1(t,s)\dif s\in L^\infty(\mR_+),\label{Con3}
\ee
then for almost all $t\geq 0$
\be
\mE\|X(t)\|^p_\mX&\leq& C_{p,\kappa_1}\bigg(\mE\|g(t)\|_\mX^p
+\int^t_0\tilde\kappa_1(t,s)\cdot\mE\|g(s)\|_\mX^p\dif s\no\\
&&+\int^t_0r_{\tilde\kappa_1}(t,u)\cdot\left[\int^u_0\tilde\kappa_1(u,s)
\cdot\mE\|g(s)\|_\mX^p\dif s\right]\dif u\bigg),\label{PP02}
\ee
where $\tilde\kappa_1=\tilde C_{p,\kappa_1}\cdot\kappa_1$,
$r_{\tilde\kappa_1}$ is defined by (\ref{R1}) in terms of $\tilde\kappa_1$,
 and $C_{p,\kappa_1},\tilde C_{p,\kappa_1}$ are constants
only depending on $p,\kappa_1$.
\et
\begin{proof}
We use Picard's iteration to prove the existence.
Let $X_1(t):=g(t)$, and define recursively for $n\in\mN$
\be
X_{n+1}(t)=g(t)+\int^t_0A(t,s,X_n(s))\dif s+\int^t_0B(t,s,X_n(s))\dif W(s).\label{Es4}
\ee
Fix $T>0$ below. By {\bf (H2)}, the BDG inequality (\ref{BDG}) and H\"older's inequality we have
\be
\mE\|X_{n+1}(t)\|^p_\mX&\preceq&\mE\|g(t)\|_\mX^p
+\mE\left(\int^t_0\|A(t,s,X_n(s))\|_\mX\dif s\right)^p\no\\
&&+\mE\left\|\int^t_0B(t,s,X_n(s))\dif W(s)\right\|^p_\mX\no\\
&\preceq&\mE\|g(t)\|_\mX^p+\mE\left(\int^t_0\kappa_{1}(t,s)\cdot(\|X_n(s)\|_\mX+1)\dif s\right)^p\no\\
&&+\mE\left(\int^t_0\|B(t,s,X_n(s))\|_{L_2(l^2;\mX)}^2\dif s\right)^{\frac{p}{2}}\no\\
&\preceq&\mE\|g(t)\|_\mX^p+\int^t_0\kappa_{1}(t,s)\cdot\mE(\|X_n(s)\|^p_\mX+1)\dif s
\cdot\left(\int^t_0\kappa_{1}(t,s)\dif s\right)^{p-1}\no\\
&&+\int^t_0\kappa_{1}(t,s)\cdot\mE(\|X_n(s)\|^p_{\mX}+1)\dif s
\cdot\left(\int^t_0\kappa_{1}(t,s)\dif s\right)^{\frac{p}{2}-1}\no\\
&\preceq&\mE\|g(t)\|_\mX^p+C_{T,p}\cdot C_T
+C_{T,p}\int^t_0\kappa_1(t,s)\cdot\mE\|X_n(s)\|^p_\mX\dif s,\label{Con4}
\ee
where $C_T:=\mathrm{ess}\sup_{t\in[0,T]}|\int^t_0\kappa_{1}(t,s)\dif s|$
and $C_{T,p}:=C_T^{p-1}+C^{(p-2)/2}_T$.

Set
$$
f_m(t):=\sup_{n=1,\cdots, m}\mE\|X_n(t)\|^p_\mX.
$$
Then
$$
f_m(t)\leq C_{T,p,\kappa_1}\Big(\mE\|g(t)\|_\mX^p+1\Big)+\int^t_0\tilde\kappa_1(t,s)\cdot f_m(s)\dif s,
$$
where $\tilde\kappa_1=C_{T,p,\kappa_1}\cdot\kappa_1$
and the constant $C_{T,p,\kappa_1}$ is independent of $m$.

Let $r_{\tilde\kappa_1}$ be defined by (\ref{R1}) in terms of $\tilde\kappa_1$.
Note that by (\ref{R2})
\ce
&&\int^t_0r_{\tilde\kappa_1}(t,s)\cdot\mE\|g(s)\|_\mX^p\dif s
-\int^t_0\tilde\kappa_1(t,s)\cdot\mE\|g(s)\|_\mX^p\dif s\no\\
&&\quad=\int^t_0\left(\int^t_sr_{\tilde\kappa_1}(t,u)
\tilde\kappa_1(u,s)\dif u\right)\cdot\mE\|g(s)\|_\mX^p\dif s\no\\
&&\quad=\int^t_0r_{\tilde\kappa_1}(t,u)\left(\int^u_0
\tilde\kappa_1(u,s)\cdot\mE\|g(s)\|_\mX^p\dif s\right)\dif u.
\de
Hence, by Lemma \ref{Le1} and {\bf (H1)},  we obtain that for almost all $t\in[0,T]$
\be
\sup_{n\in\mN}\mE\|X_{n}(t)\|^p_\mX
&=&\lim_{m\to\infty}f_m(t)\leq C_{T,p,\kappa_1}\left(\mE\|g(t)\|_\mX^p
+\int^t_0r_{\tilde\kappa_1}(t,s)\cdot\mE\|g(s)\|_\mX^p\dif s\right)\no\\
&\leq&C_{T,p,\kappa_1}\bigg(\mE\|g(t)\|_\mX^p
+\int^t_0\tilde\kappa_1(t,s)\cdot\mE\|g(s)\|_\mX^p\dif s\no\\
&&\qquad+\int^t_0r_{\tilde\kappa_1}(t,u)\left(\int^u_0
\tilde\kappa_1(u,s)\cdot\mE\|g(s)\|_\mX^p\dif s\right)\dif u\bigg)\label{PP5}\\
&\stackrel{(\ref{R3})}{\leq}&C_{T,p,\kappa_1}\left[\mE\|g(t)\|_\mX^p
+\mathrm{ess\!\!}\sup_{t\in[0,T]}\int^t_0\kappa_1(t,s)\cdot\mE\|g(s)\|_\mX^p\dif s\right].\label{Prio}
\ee

On the other hand, set
$$
Z_{n,m}(t):=X_n(t)-X_m(t).
$$
As the above calculations, by {\bf (H3)} we have
\ce
\mE\|Z_{n+1,m+1}(t)\|^2_\mX&\preceq&
\mE\left\|\int^t_0(A(t,s,X_n(s))-A(t,s,X_m(s)))\dif s\right\|^2_\mX\\
&&+\mE\left\|\int^t_0(B(t,s,X_n(s))-B(t,s,X_m(s)))\dif W(s)\right\|^2_\mX\\
&\preceq&\int^t_0\kappa_2(t,s)\cdot\mE\|Z_{n,m}(s)\|^2_\mX\dif s.
\de
Set
$$
f(t):=\limsup_{n,m\rightarrow\infty}\mE\|Z_{n,m}(t)\|^2_\mX.
$$
By (\ref{Prio}), {\bf (H1)} and using Fatou's lemma, we get
$$
f(t)\preceq\int^t_0\kappa_2(t,s)\cdot f(s)\dif s.
$$
By Lemma \ref{Le1} again,  we have for almost all $t\in[0,T]$
$$
f(t)=\limsup_{n,m\rightarrow\infty}\mE\|Z_{n,m}(t)\|^2_\mX=0.
$$
Hence, there exists an $\mX$-valued $(\cF_t)$-adapted process $X(t)$ such that
for almost all $t\in[0,T]$
$$
\lim_{n\rightarrow\infty}\mE\|X_{n}(t)-X(t)\|^2_\mX=0.
$$
Taking limits for (\ref{Es4}), one finds that (\ref{PP1}) holds.

Moreover, the estimate (\ref{PP2})  follows from (\ref{Prio}).
Note that when (\ref{Con3}) is satisfied, the constant $C_{T,p}$ in (\ref{Con4})
is independent of $T$. Hence, the estimate (\ref{PP02}) is direct from (\ref{PP5}).
The uniqueness follows by similar calculations as above.
\end{proof}
\bx{\rm\label{Exam}
Let for $\d>0$
$$
h(s):=\frac{e^{-\d s}}{s\log^2 s},\ \ t>s\geq0.
$$
It is easy to see that $h\in L^1(\mR_+)$.
Consider the following stochastic Volterra equation:
$$
X(t)=x_0\sqrt{|\log(t\wedge 1)|}+\int^t_0h(t-s)A(X(s))\dif s+\int^t_0\sqrt{h(t-s)}B(X(s))\dif W(s),
$$
where $A:\mX\to\mX$ and $B:\mX\to L_2(l^2;\mX)$ are global Lipschitz continuous functions.
By  elementary calculations, one finds that
$$
\sup_{t\geq 0}\int^t_0\frac{e^{-\d(t-s)}|\log (s\wedge 1)|}{(t-s)\log^2(t-s)}\dif s<+\infty.
$$
So, {\bf (H1)}-{\bf (H3)} are satisfied with $p=2$. Moreover, by (\ref{Wi}) and
(\ref{PP02}), one finds that if $\d$ is large enough, then for any $T>0$
$$
\sup_{t\geq T}\mE\|X(t)\|^2_\mX<+\infty.
$$
We remark that in this example, $X(0)=\infty$.
}
\ex

\subsection{Path continuity of solutions}
In this subsection,
in addition to {\bf (H2)} and {\bf (H3)}, we also assume that
\begin{enumerate}[{\bf (H1)$'$}]
\item
The process $t\mapsto g(t)$ is continuous and $(\cF_t)$-adapted,
and for any $p\geq 2$ and $T>0$
$$
\mE\left(\sup_{t\in[0,T]}\|g(t)\|^p_\mX\right)<+\infty.
$$
\end{enumerate}

\begin{enumerate}[{\bf (H4)}]
\item For all $s<t<t'$, $\omega\in\Omega$ and $x\in\mX$
$$
\|A(t',s,\omega,x)-A(t,s,\omega,x)\|_\mX\leq \lambda(t',t,s)\cdot(\|x\|_\mX+1)
$$
and
$$
\|B(t',s,\omega,x)-B(t,s,\omega,x)\|^2_{L_2(l^2;\mX)}\leq \lambda(t',t,s)\cdot(\|x\|_\mX^2+1),
$$
where $\l$ is a positive measurable function satisfying that for any $T>0$
and some $\g=\g(T),C=C(T)>0$
\be
\int^t_0\lambda(t',t,s)\dif s\leq C|t'-t|^\g,\ \ 0\leq t<t'\leq T.\label{C1}
\ee
\end{enumerate}

\bt\label{Th1}
Assume that {\bf (H1)$'$} and {\bf (H2)}-{\bf (H4)} hold, and
the kernel function $\kappa_1$ in {\bf (H2)} belongs to $\sK_{>1}$.
Then there exists a unique
$\mX$-valued continuous $(\cF_t)$-adapted process $X(t)$ such that $P$-a.s.,
for all $t\geq0$
\be
X(t)=g(t)+\int^t_0A(t,s,X(s))\dif s+\int^t_0B(t,s,X(s))\dif W(s)
\label{PPP1}
\ee
and for any $p\geq 2$ and $T>0$,
\be
\mE\left(\sup_{t\in[0,T]}\|X(t)\|^p_\mX\right)<+\infty.\label{PPP2}
\ee
Moreover, if for some $\delta>0$ and any $p\geq 2, T>0$, it holds that
$$
\mE\|g(t')-g(t)\|^p_\mX\leq C_{T,p}|t'-t|^{\delta p},
$$
then, $t\mapsto X(t)$ admits a H\"older continuous modification
and for any $p\geq 2, T>0$ and some $a>0$
$$
\mE\left(\sup_{t\not= t'\in[0,T]}
\frac{\|X(t')-X(t)\|^p_{\mX}}
{|t'-t|^{ap}}\right)\leq C_{T,p,a}.
$$
\et

\begin{proof}
First of all, for any $p\geq 2$ and $T>0$, by {\bf (H1)$'$} and (\ref{PP2})  we have
\be
\mathrm{ess\!\!}\sup_{t\in[0,T]}\mE\|X(t)\|^p_\mX<+\infty.\label{PP0}
\ee
Set
$$
J(t):=\int^t_0B(t,s,X(s))\dif W(s)
$$
and write for $0\leq t<t'\leq T$
\ce
J(t')-J(t)&=&\int^t_0\big[B(t',s,X(s))-B(t,s,X(s))\big]\dif W(s)\\
&&+\int^{t'}_tB(t',s,X(s))\dif W(s)=:J_1(t',t)+J_2(t',t).
\de
In view of $\kappa_1\in\sK_{>1}$,  (\ref{C2}) holds for some $\beta>1$.
Fix $p\geq 2(\beta^*:=\beta/(\beta-1))$.
By the BDG inequality (\ref{BDG}),
{\bf (H2)} and H\"older's inequality  we have
\ce
\mE\|J_2(t',t)\|^p_\mX&\preceq&\mE\left(\int^{t'}_t\kappa_1(t',s)
\cdot(\|X(s)\|^2_\mX+1)\dif s\right)^{\frac{p}{2}}\\
&\preceq&\left(\int^{t'}_tk^\beta_1(t',s)\dif s\right)^{\frac{p}{2\beta}}
\mE\left(\int^{t'}_t(\|X(s)\|^{2\beta^*}_\mX+1)\dif s\right)^{\frac{p}{2\beta^*}}\\
&\stackrel{(\ref{C2})}{\preceq}&|t'-t|^{\frac{p}{2\beta^*}-1}
\int^{t'}_t(\mE\|X(s)\|^{p}_\mX+1)\dif s\\
&\stackrel{(\ref{PP0})}{\preceq}&|t'-t|^{\frac{p}{2\beta^*}},
\de
and by {\bf (H4)} and Minkowski's inequality
\ce
\mE\|J_1(t',t)\|^p_\mX&\preceq&\mE\left(\int^t_0\lambda(t',t,s)
\cdot(\|X(s)\|^2_\mX+1)\dif s\right)^{\frac{p}{2}}\\
&\preceq&\left(\int^t_0\lambda(t',t,s)\cdot((\mE\|X(s)\|^p_\mX)^{\frac{2}{p}}+1)
\dif s\right)^{\frac{p}{2}}\\
&\stackrel{(\ref{PP0})}{\preceq}&\left(\int^t_0\lambda(t',t,s)\dif s\right)^{\frac{p}{2}}\\
&\stackrel{(\ref{C1})}{\preceq}&|t'-t|^{\frac{\g p}{2}}.
\de
Hence, for all $0\leq t<t'\leq T$
$$
\mE\|J(t')-J(t)\|^p_\mX\preceq |t-t'|^{\frac{\g p}{2}}+|t-t'|^{\frac{p}{2\beta^*}}.
$$

Similarly, we may prove that for all $0\leq t<t'\leq T$ and $p\geq \beta^*$
$$
\mE\left\|\int^{t'}_0A(t',s,X(s))\dif s-\int^t_0A(t,s,X(s))\dif s\right\|^p_\mX
\preceq |t-t'|^{\g p}+|t-t'|^{\frac{p}{\beta^*}}.
$$
The desired conclusions follow from Theorem \ref{KCC}.
\end{proof}

We conclude this subsection by proving a lemma, which will be used frequently later.
We put it here since the proof is similar to Theorem \ref{Th1}.
\bl\label{Le00}
Let $\tau$ be an ($\cF_t$)-stopping time and
$$
G:\triangle\times\Omega\to L_2(l^2;\mX)\in\cM_\triangle/\cB(L_2(l^2;\mX)).
$$
Assume that for all $0\leq s<t<t'$ and $\omega\in\Omega$
\be
\|G(t,s,\omega)\|_{L_2(l^2;\mX)}^2&\leq& \kappa(t,s)\cdot f^2(s,\omega),\label{PL3}\\
\|G(t',s,\omega)-G(t,s,\omega)\|^2_{L_2(l^2;\mX)}&\leq& \l(t',t,s)\cdot f^2(s,\omega),\label{PL4}
\ee
where $\kappa\in\sK_{>1}$ and for any $T>0$ and some $\a>1$ and $\g>0$
$$
\int^t_0\l^{\a}(t',t,s)\dif s\leq C_T|t'-t|^\g,\ \ \forall 0\leq t<t'\leq T,
$$
and $(s,\omega)\mapsto f(s,\omega)$ is a positive measurable process with
$$
\mE\left(\int^{T\wedge\tau}_0f^p(s)\dif s\right)<+\infty,\ \ \forall p\geq 2.
$$
Then $t\mapsto J(t):=\int^t_0G(t,s)\dif W(s)\in\mX$ admits a continuous modification on
$[0,\tau)$, and for any $T>0$ and $p$ large enough
\ce
\mE\left(\sup_{t\in[0,T\wedge\tau]}\left\|\int^t_0G(t,s)\dif W(s)\right\|^p_\mX\right)\leq
C_T\mE\left(\int^{T\wedge\tau}_0f^p(s)\dif s\right),
\de
where the constant $C_T$ is independent of $f$ and $\tau$.
\el
\begin{proof}
Fix $T>0$ and write for $0\leq t<t'\leq T$
\ce
J(t')-J(t)&=&\int^{t'}_{t}G(t',s)\dif W(s)+\int^{t}_0[G(t',s)-G(t,s)]\dif W(s)\\
&=:&J_1(t',t)+J_2(t',t).
\de
In view of $\kappa\in\sK_{>1}$ and (\ref{C2}),
by the BDG inequality (\ref{BDG}) and H\"older's inequality we have,
for some $\beta>1$ and $p\geq 2(\beta^*=\beta/(\beta-1))$,
\ce
\mE\|J_1(t',t)\cdot 1_{\{t',t\in[0,\tau)\}}\|^p_\mX
&\leq& \mE\left\|\int^{t'\wedge\tau}_{t\wedge\tau}G(t',s)\dif W(s)\right\|_\mX^{p}\\
&\preceq& \mE\left(\int^{t'\wedge\tau}_{t\wedge\tau}
\|G(t',s)\|^2_{L_2(l^2;\mX)}\dif s\right)^{p/2}\\
&\stackrel{(\ref{PL3})}{\preceq}& \mE\left(\int^{t'\wedge\tau}_{t\wedge\tau}
\kappa(t',s)\cdot f^2(s)\dif s\right)^{p/2}\\
&\preceq& \left(\int^{t'}_{t}\kappa^{\beta}(t',s)\dif s\right)^{\frac{p}{2\beta}}\cdot
\mE\left(\int^{t'\wedge\tau}_{t\wedge\tau} f^{2\beta^*}(s)\dif s\right)^{\frac{p}{2\beta^*}}\\
&\preceq&|t'-t|^{\frac{ p}{2\beta^*}-1}\cdot
\mE\left(\int^{T\wedge\tau}_0f^p(s)\dif s\right)
\de
and for $p\geq 2(\a^*=\a/(\a-1))$,
\ce
\mE\|J_1(t',t)\cdot 1_{\{t',t\in[0,\tau)\}}\|^p_\mX&\preceq& \mE\left(\int^{t\wedge\tau}_0
\|G(t',s)-G(t,s)\|^2_{L_2(l^2;\mX)}\dif s\right)^{p/2}\\
&\stackrel{(\ref{PL4})}{\preceq}& \mE\left(\int^{t\wedge\tau}_0
\l(t',t,s)\cdot f^2(s)\dif s\right)^{p/2}\\
&\preceq& \left(\int^{t}_0\l^{\a}(t',t,s)\dif s\right)^{\frac{p}{2\a}}\cdot
\mE\left(\int^{t\wedge\tau}_0f^{2\a^*}(s)\dif s\right)^{\frac{p}{2\a^*}}\\
&\preceq& |t'-t|^{\frac{\gamma p}{2\a}}\cdot
\mE\left(\int^{T\wedge\tau}_0f^p(s)\dif s\right).
\de
Hence, for any $p\geq 2(\a^*\vee\beta^*)$ and $0\leq t<t'\leq T$,
$$
\mE\|(J(t')-J(t))\cdot 1_{\{t',t\in[0,\tau)\}}\|^p_\mX\preceq
|t'-t|^{\left(\frac{ p}{2\beta^*}-1\right)\wedge\frac{\gamma p}{2\a}}\cdot
\mE\left(\int^{T\wedge\tau}_0f^p(s)\dif s\right).
$$
The desired result now follows by Theorem \ref{KCC}.
\end{proof}

\subsection{Local existence and uniqueness}

In this subsection, we assume that
\begin{enumerate}[{\bf (H2)$'$}]
\item For any $R>0$, there exists  $\kappa_{1,R}\in\sK_{>1}$
such that for all $(t,s)\in\triangle$, $\omega\in\Omega$ and $x\in\mX$
with $\|x\|_\mX\leq R$
$$
\|A(t,s,\omega,x)\|_\mX+\|B(t,s,\omega,x)\|^2_{L_2(l^2;\mX)}\leq \kappa_{1,R}(t,s).
$$

\end{enumerate}
\begin{enumerate}[{\bf (H3)$'$}]

\item For any $R>0$, there exists  $\kappa_{2,R}\in\sK_0$
such that for all $(t,s)\in\triangle$, $\omega\in\Omega$ and $x,y\in\mX$
with $\|x\|_\mX,\|y\|_\mX\leq R$
$$
\|A(t,s,\omega,x)-A(t,s,\omega,y)\|_\mX\leq \kappa_{2,R}(t,s)\cdot\|x-y\|_\mX
$$
and
$$
\|B(t,s,\omega,x)-B(t,s,\omega,y)\|^2_{L_2(l^2;\mX)}\leq \kappa_{2,R}(t,s)\cdot\|x-y\|_\mX^2.
$$
\end{enumerate}
\begin{enumerate}[{\bf (H4)$'$}]
\item For any $R>0$, there exists a measurable function $\l_R$ satisfying that
for  any $T>0$ and some $\g, C>0$
$$
\int^t_0\lambda_R(t',t,s)\dif s\leq C|t'-t|^\g,\ \ 0\leq t<t'\leq T,
$$
such that for all $s<t<t'$, $\omega\in\Omega$ and $x\in\mX$ with $\|x\|_\mX\leq R$,
\ce
&&\|A(t',s,\omega,x)-A(t,s,\omega,x)\|_\mX
+\|B(t',s,\omega,x)-B(t,s,\omega,x)\|^2_{L_2(l^2;\mX)}\\
&&\qquad\qquad\leq \lambda_R(t',t,s).
\de
\end{enumerate}

We first introduce the following notion of local solutions.
\bd\label{Def1}
Let $\tau$ be an $(\cF_t)$-stopping time, and $\{X(t); t\in[0,\tau)\}$
an $\mX$-valued continuous $(\cF_t)$-adapted process.
The pair of $(X,\tau)$ is called a local solution of \textsc{Eq.}(\ref{Eq}) if
$P$-a.s., for all $t\in[0,\tau)$
\ce
X(t)=g(t)+\int^{t}_0A(t,s,X(s))\dif s+\int^{t}_0B(t,s,X(s))\dif W(s);
\de
$(X,\tau)$ is called a maximal solution of \textsc{Eq.}(\ref{Eq}) if
$$
\lim_{t\uparrow\tau(\omega)}\|X(t,\omega)\|_\mX=+\infty \
\mbox{ on $~~\{\omega: \tau(\omega)<+\infty\}$},\ P-a.s..
$$
We call $(X,\tau)$ a non-explosion solution of \textsc{Eq.}(\ref{Eq}) if
$$
P\{\omega: \tau(\omega)<+\infty\}=0.
$$
\ed
\br
The stochastic integral in the above definition is defined on $[0,\tau)$ by
$$
\int^{t}_0B(t,s,X(s))\dif W(s)=\lim_{n\rightarrow\infty}
\int^{t\wedge\tau_n}_0B(t,s,X(s))\dif W(s),\ \ t<\tau,
$$
where $\tau_n:=\inf\{t>0: \|X(t)\|_\mX>n\}\nearrow\tau$.
\er

We now prove the following main result in this section.
\bt\label{Main}
Under {\bf (H1)$'$}-{\bf (H4)$'$}, there exists a unique
maximal solution $(X,\tau)$ for \textsc{Eq.}(\ref{Eq}) in the sense of Definition \ref{Def1}.
\et
\begin{proof}
For $n\in\mN$, let $\chi_n$ be a positive smooth function on $\mR_+$
with $\chi_n(s)=1, s\leq n$ and $\chi_n(s)=0, s\geq  n+1$.
Define
\ce
A_n(t,s,\omega,x)&:=&A(t,s,\omega,x)\cdot \chi_n(\|x\|_\mX)\\
B_n(t,s,\omega,x)&:=&B(t,s,\omega,x)\cdot \chi_n(\|x\|_\mX).
\de
It is easy to see that for $A_n$ and $B_n$,
{\bf (H2)} holds with $\kappa_{1,n+1}$,
{\bf (H4)} holds with $\lambda_{n+1}$, and  {\bf (H3)} holds with some
$\kappa_{3,n}\in\sK_0$.
Thus, by Theorem \ref{Th1} there exists a unique continuous ($\cF_t$)-adapted process
$X_n(t)$ such that for any $p\geq 2$ and $T>0$
$$
\mE\left(\sup_{t\in[0,T]}\|X_n(t)\|^p_\mX\right)\leq C_{T,p,n}
$$
and
\be
X_n(t)=g(t)+\int^t_0A_n(t,s,X_n(s))\dif s+\int^t_0B_n(t,s,X_n(s))\dif W(s).\label{App}
\ee

We have the following claim:
$$
\mbox{\it Let $\tau$ be any stopping time.
The uniqueness holds for (\ref{App}) on $[0,\tau)$.}
$$
We remark that when $\tau=T$ is non-random, it follows from Theorem \ref{Th0}.
Let $X_i(t), i=1,2$ be two $\mX$-valued continuous $(\cF_t)$-adapted processes, and
satisfy on $[0,\tau)$
$$
X_i(t)=g(t)+\int^{t}_0
A_n(t,s,X_i(s))\dif s+\int^{t}_0B_n(t,s,X_i(s))\dif W(s),\ i=1,2.
$$
Set
$$
Z(t):=X_1(t)-X_2(t).
$$
Since $\kappa_{3,n}\in\sK_0$, as the calculations in (\ref{Con4}), by
the BDG  inequality (\ref{BDG}) and {\bf (H3)} for $A_n$ and $B_n$, we have
\be
\mE\|Z(t)\cdot 1_{\{t<\tau\}}\|^p_\mX&\preceq&\mE\left(\int^{t\wedge\tau}_0\kappa_{3,n}(t,s)
\cdot\|Z(s)\|_\mX\dif s\right)^p\no\\
&&+\mE\left(\int^{t\wedge\tau}_0\kappa_{3,n}(t,s)
\cdot\|Z(s)\|^2_\mX\dif s\right)^{\frac{p}{2}}\no\\
&=&\mE\left(\int^{t}_0\kappa_{3,n}(t,s)\cdot 1_{\{s<\tau\}}
\cdot\|Z(s)\|_\mX\dif s\right)^p\no\\
&&+\mE\left(\int^{t}_0\kappa_{3,n}(t,s)\cdot 1_{\{s<\tau\}}
\cdot\|Z(s)\|^2_\mX\dif s\right)^{\frac{p}{2}}\no\\
&\preceq&\int^{t}_0\kappa_{3,n}(t,s)\cdot \mE\|Z(s)\cdot 1_{\{s<\tau\}}\|^p_\mX\dif s.\label{Es5}
\ee
By Lemma \ref{Le1}, we get
$$
\mE\|Z(t)\cdot 1_{\{t<\tau\}}\|^p_\mX=0\ \ \mbox{ for almost all $t\in[0,T]$},
$$
which implies by the arbitrariness of $T$ and the continuities of $X_i(t), i=1,2$,
$$
X_1(\cdot)|_{[0,\tau)}=X_2(\cdot)|_{[0,\tau)}.
$$
The claim is proved.

\vspace{2mm}

Now, for $n\in\mN$, define the stopping times
$$
\tau_n:=\inf\{t>0: \|X_n(t)\|_\mX>n\}
$$
and
$$
\sigma_n:=\inf\{t>0: \|X_{n+1}(t)\|_\mX>n\}.
$$
By the above claim, we have
$$
X_n(\cdot)|_{[0,\tau_n\wedge\sigma_n)}=X_{n+1}(\cdot)|_{[0,\tau_n\wedge\sigma_n)},
$$
which implies
$$
\tau_n\leq\sigma_n\leq\tau_{n+1},\ \ a.e..
$$

Hence, we may define
$$
\tau(\omega):=\lim_{n\to\infty}\tau_n(\omega)
$$
and for all $t<\tau(\omega)$
$$
X(t,\omega):=X_n(t,\omega), \ \ \mbox{ if $t<\tau_n(\omega)$}.
$$
Clearly, $(X,\tau)$ is a maximal solution of \textsc{Eq.}(\ref{Eq}) in the sense of Definition \ref{Def1}.

We next prove the uniqueness. Let $(\tilde X,\tilde\tau)$ be another maximal
solution of \textsc{Eq.}(\ref{Eq}) in the sense of Definition \ref{Def1}.
Define the stopping times
$$
\tilde\tau_n:=\inf\{t>0: \|\tilde X(t)\|_\mX>n\}
$$
and
$$
\hat\tau_n:=\tau_n\wedge\tilde\tau_n,\ \ \hat\tau:=\tau\wedge\tilde\tau.
$$
It is clear that
$$
\hat\tau_n\nearrow\hat\tau\ \  a.s.\ \mbox{ as $n\to\infty$}
$$
and
\ce
1_{[0,\hat\tau_n)}(t)\cdot\tilde X(t)&=&1_{[0,\hat\tau_n)}(t)\cdot g(t)
+1_{[0,\hat\tau_n)}(t)\cdot\int^{t}_0A(t,s,\tilde X(s))\dif s\\
&&+1_{[0,\hat\tau_n)}(t)\cdot\int^{t}_0B(t,s,\tilde X(s))\dif W(s)\\
&=&1_{[0,\hat\tau_n)}(t)\cdot g(t)
+1_{[0,\hat\tau_n)}(t)\cdot\int^{t}_0A_n(t,s,\tilde X(s))\dif s\\
&&+1_{[0,\hat\tau_n)}(t)\cdot\int^{t}_0B_n(t,s,\tilde X(s))\dif W(s).
\de
By the above claim again, we have
$$
X(\cdot)|_{[0,\hat\tau_n)}=\tilde X(\cdot)|_{[0,\hat\tau_n)}.
$$
So
$$
X(\cdot)|_{[0,\hat\tau)}=\tilde X(\cdot)|_{[0,\hat\tau)}.
$$
By the definition of maximal solution we must have $\hat\tau=\tau=\tilde\tau$.
\end{proof}

We have the following simple criterion of non explosion.
\bt\label{Th4}
Assume that {\bf (H1)$'$}, {\bf (H2)} and {\bf (H4)} hold, and
$\kappa_1$ in  {\bf (H2)} belongs to $\sK_{>1}$.
Then there is no explosion for \textsc{Eq.}(\ref{Eq}).
\et
\begin{proof}Let $(X,\tau)$ be a maximal solution of \textsc{Eq.}(\ref{Eq}). Define
$$
\tau_n:=\inf\{t>0: \|X(t)\|_\mX\geq n\}.
$$
By the BDG inequality (\ref{BDG}) and H\"older's inequality, and using the same method
as estimating (\ref{Es5}), we have, for any $T>0$, some $\beta>1$ and
$p\geq 2(\beta^*=\beta/(\beta-1))$
\ce
\mE\|X(t)\cdot 1_{\{t\leq\tau_n\}}\|^p_\mX&\preceq&\mE\|g(t)\|_\mX^p
+\mE\left(\int^{t\wedge\tau_n}_0\|A(t,s,X(s))\|_\mX\dif s\right)^p\\
&&+\mE\left\|\int^{t\wedge\tau_n}_0
B(t,s,X(s))\dif W(s)\right\|^p_\mX\\
&\preceq&\mE\|g(t)\|_\mX^p+\mE\left(\int^{t\wedge\tau_n}_0
\kappa_{1}(t,s)\cdot(\|X(s)\|_\mX+1)\dif s\right)^p\\
&&+\mE\left(\int^{t\wedge\tau_n}_0\|B(t,s,X(s))\|_{L_2(l^2;\mX)}^2\dif s\right)^{\frac{p}{2}}\\
&\preceq&\mE\|g(t)\|_\mX^p
+\mE\left(\int^{t\wedge\tau_n}_0(\|X(s)\|^{\beta^*}_\mX+1)\dif s\right)^{\frac{p}{\beta^*}}\\
&&+\mE\left(\int^{t\wedge\tau_n}_0(\|X(s)\|^{2\beta^*}_\mX+1)\dif s\right)^{\frac{p}{2\beta^*}}\\
&\leq&C_{T,p}\left[\mE\|g(s)\|_\mX^p+1
+\int^t_0\mE\|X(s)\cdot 1_{\{s\leq\tau_n\}}\|^p_\mX\dif s\right],
\de
where the constant $C_{T,p}$ is independent of $n$.

By Gronwall's inequality, we obtain
$$
\sup_{t\in[0,T]}\mE\|X(t)\cdot 1_{\{t\leq\tau_n\}}\|^p_\mX\leq C_{T,p}.
$$
Using this estimate, as in the proofs of Theorem \ref{Th1} and Lemma \ref{Le00}, we can prove
that for any $T>0$ and $p\geq 2$
$$
\sup_{n\in\mN}\mE\left(\sup_{t\in[0,T\wedge\tau_n]}\|X(t)\|^p_\mX\right)\leq C_{T,p}.
$$
Hence,
\ce
\lim_{n\to\infty}P\{\tau_n\leq T\}&=&
\lim_{n\to\infty}P\left\{\sup_{t\in[0,T\wedge\tau_n]}\|X(t)\|_\mX\geq n\right\}\\
&\leq&\lim_{n\to\infty}\mE\left(\sup_{t\in[0,T\wedge\tau_n]}\|X(t)\|^p_\mX\right)/n^p\\
&\leq&\lim_{n\to\infty}C_{T,p}/n^p=0,
\de
which produces the non-explosion, i.e., $P\{\tau<\infty\}=0$.
\end{proof}

\br
One cannot directly prove
$$
\sup_{n\in\mN}\mE\|X(t\wedge\tau_n)\|^p_\mX<+\infty,\ \ \forall t\geq0
$$
to obtain the non-explosion, because it does not in general make sense to write
$$
\int^{t\wedge\tau_n}_0B(t\wedge\tau_n,s,X(s))\dif W(s).
$$
\er
\subsection{Continuous dependence of solutions with respect to data}
In this subsection, we study the continuous dependence of solutions for \textsc{Eq.}(\ref{Eq})
with respect to the  coefficients.

Let $\{(g_m, A_m, B_m), m\in\mN\}$ be a sequence of coefficients associated to
\textsc{Eq.}(\ref{Eq}). Assume that for each $m\in\mN$,
$(g_m, A_m, B_m)$ satisfies {\bf (H1)$'$}-{\bf (H4)$'$}
with the same $\kappa_{1,R},\kappa_{2,R}$ and $\l_R$ as $(g, A, B)$, and
for each $p\geq 2$
\be
\lim_{m\to\infty}\sup_{t\in[0,T]}\mE\|g_m(t)-g(t)\|^p_\mX=0\label{V1}
\ee
and for each $T,R>0$,
\be
&&\lim_{m\to\infty}\sup_{t\in[0,T],\|x\|_{\mX}\leq R}\int^t_0\|A_m(t,s,x)-A(t,s,x)\|_\mX\dif s=0,\label{V2}\\
&&\lim_{m\to\infty}\sup_{t\in[0,T],\|x\|_{\mX}\leq R}\int^t_0
\|B_m(t,s,x)-B(t,s,x)\|^2_{L_2(l^2;\mX)}\dif s=0.\label{V3}
\ee
Let $(X_m,\tau_m)$ (resp. $(X,\tau)$)
be the unique maximal solution associated with $(g_m, A_m, B_m)$
(resp. $(g, A, B)$). For each $R>0$ and $m\in\mN$, define
$$
\tau^R_m:=\inf\{t>0: \|X(t)\|_\mX,\|X_m(t)\|_\mX>R\}.
$$
Suppose that for each $t>0$
\be
\lim_{R\to\infty}\sup_mP\{\tau^R_m<t\}=0.\label{V4}
\ee
Then we have:
\bt
For each $t>0$ and $\eps>0$
$$
\lim_{m\to\infty}P\big\{\|X_m(t)-X(t)\|_\mX
\geq \eps\big\}=0.
$$
\et
\begin{proof}
For $R>0$ and $m\in\mN$, set
$$
Z^R_m(t):=(X_m(t)-X(t))\cdot 1_{\{t\leq\tau^R_m\}}.
$$
Then
$$
Z^R_m(t)=J^R_{1,m}(t)+J^R_{2,m}(t)+J^R_{3,m}(t)+J^R_{4,m}(t)+J^R_{5,m}(t),
$$
where
\ce
J^R_{1,m}(t)&:=&1_{\{t\leq\tau^R_m]}\cdot[g_m(t)-g(t)],\\
J^R_{2,m}(t)&:=&1_{\{t\leq\tau^R_m]}\cdot\int^{t\wedge\tau^R_m}_0\big[A_m(t,s,X_n(s))
-A_m(t,s,X(s))\big]\dif s,\\
J^R_{3,m}(t)&:=&1_{\{t\leq\tau^R_m]}\cdot
\int^{t\wedge\tau^R_m}_0\big[A_m(t,s,X(s))-A(t,X(s))\big]\dif s,\\
J^R_{4,m}(t)&:=&1_{\{t\leq\tau^R_m]}\cdot\int^{t\wedge\tau^R_m}_0\big[B_m(t,s,X_m(s))
-B_m(t,s,X(s))\big]\dif W(s),\\
J^R_{5,m}(t)&:=&1_{\{t\leq\tau^R_m]}\cdot\int^{t\wedge\tau^R_m}_0\big[B_m(t,s,X(s))
-B(t,s,X(s))\big]\dif W(s).
\de
Fix $T>0$. Clearly, for any $p\geq 2$ and $t\in[0,T]$
$$
\mE\|J^R_{1,m}(t)\|^p_\mX\leq\sup_{t\in[0,T]}\mE\|g_m(t)-g(t)\|^p_\mX=:\cJ_{1,m}.
$$
For $J^R_{2,m}(t)$, by {\bf (H3)$'$} and H\"older's inequality we have, for $p$ large enough
($\kappa_{2,R}\in\sK_{>1}$)
\ce
\mE\|J^R_{2,m}(t)\|^p_\mX&\leq& \mE\left(\int^{t\wedge\tau^R_m}_0\kappa_{2,R}(t,s)
\cdot\|X_m(s)-X(s)\|_\mX\dif s\right)^p\\
&\leq&\left[\int^{t}_0\kappa^\beta_{2,R}(t,s)\dif s\right]^{\frac{p}{\beta}}
\cdot\mE\left[\int^{t }_0\|Z^R_m(s)\|^{\beta^*}_\mX\dif s\right]^{\frac{p}{\beta^*}}\\
&\leq&  C \int^t_0\mE\|Z^R_m(s)\|^p_\mX\dif s.
\de
For $J^R_{3,m}(t)$,  we have
\ce
\mE\|J^R_{3,m}(t)\|^p_\mX&\leq&\mE\left(\sup_{\|x\|_\mX\leq R}
\int^{t\wedge\tau^R_m}_0\|A_m(t,s,x)-A(t,s,x)\|_\mX\dif s\right)^p\\
&\leq&\left(\sup_{t\in[0,T]}\sup_{\|x\|_\mX\leq R}
\int^t_0\|A_m(t,s,x)-A(t,s,x)\|_\mX\dif s\right)^p=:\cJ^R_{3,m}.
\de
Similarly, by the BDG inequality (\ref{BDG}) we have, for $p$ large enough
$$
\mE\|J^R_{4,m}(t)\|^p_\mX\leq C \int^t_0\mE\|Z^R_m(s)\|^p_\mX\dif s
$$
and
$$
\mE\|J^R_{5,m}(t)\|^p_\mX\leq C_p
\left(\sup_{t\in[0,T]}\sup_{\|x\|_\mX\leq R}
\int^t_0\|B_m(t,s,x)-B(t,s,x)\|^2_{L_2(l^2;\mX)}\dif s\right)^{\frac{p}{2}}=:\cJ^R_{5,m}.
$$
Combining the above calculations, we get
$$
\mE\|Z^R_m(t)\|^p_\mX\leq \cJ_{1,m}+\cJ^R_{3,m}+\cJ^R_{5,m}
+C \int^t_0\mE\|Z^R_m(s)\|^p_\mX\dif s.
$$
By Gronwall's inequality and (\ref{V1})-(\ref{V3}) we get, for any $R>0$ and $p$ large enough
\ce
\lim_{m\to\infty}\mE\|Z^R_m(t)\|^p_\mX=0.
\de
Hence
\ce
P\big\{\|X_m(t)-X(t)\|_\mX\geq\eps\big\}
&\leq&P\big\{\|X_m(t)
-X(t)\|_\mX\cdot 1_{\{t\leq\tau^R_m\}}\geq\eps\big\}
+P\big\{\tau^R_m<t\big\}\\
&\leq&\mE\|Z^R_m(t)\|^p_\mX/\eps^p+P\big\{\tau^R_m<t\big\}.
\de
First letting $m\to\infty$, then $R\to\infty$, we then get the desired limit by (\ref{V4}).
\end{proof}

\section{Large deviation for stochastic Volterra equations}

In this section, we study the large deviation of small perturbations for
stochastic Volterra equations.
In addition to {\bf (H2)$'$}, {\bf (H3)$'$} and {\bf (H4)$'$},
we assume that $g$ and $A,B$ are non-random, and

\begin{enumerate}[{\bf (H1)$''$}]
\item For any $T>0$ and some $\d>0$,
$$
\|g(t)-g(t')\|_\mX\leq C|t-t'|^\d,\ \ t,t'\in[0,T]
$$
and for some $\a>0$,
$$
\sup_{t\in[0,T]}\|g(t)\|_{\mX_\a}<+\infty.
$$
\end{enumerate}

\begin{enumerate}[{\bf (H2)$''$}]
\item For the same $\a$ as in {\bf (H1)$''$}  and any $R>0$,
there exists a kernel function $\kappa_{\a,R}\in\sK_0$
such that for all $(t,s)\in\triangle$ and $x\in\mX$ with $\|x\|_\mX\leq R$
$$
\|A(t,s,x)\|_{\mX_\a}+\|B(t,s,x)\|^2_{L_2(l^2;\mX_{\frac{\a}{2}})}\leq \kappa_{\a,R}(t,s).
$$
\end{enumerate}
\br
If the $\kappa_{\a,R}$ in {\bf (H2)$''$} belongs to $\sK_{>1}$,
then {\bf (H2)$''$} implies {\bf (H2)$'$} in view of $\mX_\a\hookrightarrow\mX$.
\er

Consider the following small perturbation of stochastic Volterra equation (\ref{Eq})
\be
X_\eps(t)=g(t)+\int^t_0A(t,s,X_\eps(s))\dif s
+\sqrt{\eps}\int^t_0B(t,s,X_\eps(s))\dif W(s),\label{Eqq1}
\ee
where $\eps\in(0,1)$. By Theorem \ref{Main}, there exists
a unique maximal solution $(X_\eps,\tau_\eps)$
for \textsc{Eq.}(\ref{Eqq1}).
Below, we fix $T>0$ and work in the finite time interval $[0,T]$, and
assume that for each $\eps\in(0,1)$
$$
\tau_\eps>T,\ \ a.s..
$$
By Yamada-Watanabe's theorem  (cf. \cite{On,Ro-Sc-Zh}), there exists a measurable mapping
$$
\Phi_\eps:\mC_T(\mU)\to\mC_T(\mX)
$$
such that
$$
X_\eps(t,\omega)=\Phi_\eps(W(\cdot,\omega))(t).
$$
It should be noticed that although the equation considered in \cite{On} is
a little different from \textsc{Eq.}(\ref{Eq}),
 the proof is obviously adapted to our more general equation.

We now fix a family of processes $\{h^\eps,\eps\in(0,1)\}$ in $\cA^T_N$
(see (\ref{Op2}) for the definition of $\cA^T_N$), and put
$$
X^\eps(t,\omega):=\Phi_\eps\Big(W(\cdot,\omega)+\frac{h^{\eps}(\cdot,\omega)}{\sqrt{\eps}}\Big)(t).
$$
Here, we have used a little confused notations $X_\eps$ and $X^\eps$, but they are clearly different.
By Girsanov's theorem  (cf. \cite[Section 7]{On}),
$X^\eps(t)$ solves the following  stochastic Volterra equation (also called control equation):
\be
X^\eps(t)&=&g(t)+\int^t_0A(t,s,X^\eps(s))\dif s
+\int^t_0B(t,s,X^\eps(s))\dot h^\eps(s)\dif s\no\\
&&+\sqrt{\eps}\int^t_0B(t,s,X^\eps(s))\dif W(s).\label{Eq2}
\ee

Although  $h$ is defined only on $[0,T]$,
we can extend it to $\mR_+$ by setting $\dot h(t)=0$ for $t>T$
so that \textsc{Eq.}(\ref{Eq2}) can be considered on $\mR_+$. We shall always use this extension below.
Let $\tau^\eps$ be the explosion time of \textsc{Eq.}(\ref{Eq2}).
For $n\in\mN$, define
\be
\tau^\eps_n:=\inf\{t\geq 0: \|X^\eps(t)\|_\mX>n\}.\label{Stop}
\ee
Then $\tau^\eps_n\nearrow\tau^\eps$, and we have:
\bl\label{Le2}
For any $\a_0\in(0,\a)$, there is an $a>0$ such that for $p$ sufficiently large
$$
\sup_{\eps\in(0,1)}\mE\left(\sup_{t\not= t'\in[0,T\wedge\tau^\eps_n]}
\frac{\|X^\eps(t')-X^\eps(t)\|^p_{\mX_{\a_0}}}
{|t'-t|^{a p}}\right)\leq C_{N,n,T,p,\kappa_{\a,n},\a_0}.
$$
\el
\begin{proof}
Note that
\ce
\|X^\eps(t)\cdot 1_{\{t\leq\tau^\eps_n\}}\|_{\mX_\a}&\leq&\|g(t)\|_{\mX_\a}
+\int^{t\wedge\tau^\eps_n}_0\|A(t,s,X^\eps(s))\|_{\mX_\a}\dif s\\
&&+\int^{t\wedge\tau^\eps_n}_0\|B(t,s,X^\eps(s))\dot h^\eps(s)\|_{\mX_\a}\dif s\\
&&+\sqrt{\eps}\left\|\int^{t\wedge\tau^\eps_n}_0B(t,s,X^\eps(s))\dif W(s)\right\|_{\mX_\a}\\
&=:&J_1(t)+J_2(t)+J_3(t)+J_4(t).
\de
By {\bf (H2)$''$} and (\ref{Stop}) we have
$$
\mE|J_2(t)|^p\leq C_n\mE\left(\int^{t\wedge\tau^\eps_n}_0
\kappa_{\a,n}(t,s)\dif s\right)^p\leq C_{n,T,p,\kappa_{\a,n}}
$$
and by H\"older's inequality
\ce
\mE|J_3(t)|^p&\leq&\mE\left(\int^{t\wedge\tau^\eps_n}_0
\|B(t,s,X^\eps(s))\dot h^\eps(s)\|_{\mX_\a}\dif s\right)^{p}\\
&\leq&\mE\left(\int^{t\wedge\tau^\eps_n}_0
\|B(t,s,X^\eps(s))\|_{L_2(l^2;\mX_\a)}\cdot
\|\dot h^\eps(s)\|_{l^2}\dif s\right)^{p}\\
&\leq&N^{\frac{p}{2}}\mE\left(\int^{t\wedge\tau^\eps_n}_0
\|B(t,s,X^\eps(s))\|^2_{L_2(l^2;\mX_\a)}\dif s\right)^{\frac{p}{2}}\\
&\leq& C_{N,n,T,p,\kappa_{\a,n}},
\de
where we have used that $h^\eps\in\cA^T_N$.

Similarly, by the BDG inequality (\ref{BDG})
and {\bf (H2)$''$} we have
$$
\mE|J_4(t)|^p\leq C_p\mE\left(\int^{t\wedge\tau^\eps_n}_0\|B(t,s,X^\eps(s))\|^2_{L_2(l^2;\mX_\a)}
\dif s\right)^{\frac{p}{2}}\leq C_{n,T,p,\kappa_{\a,n}}.
$$
Combining the above calculations, we get
\be
\sup_{\eps\in(0,1)}\sup_{t\in[0,T]}\mE\|X^\eps(t)\cdot 1_{\{t\leq\tau^\eps_n\}}\|_{\mX_\a}^p
\leq C_{N,n,T,p,\kappa_{\a,n}},\ \ p\geq 2.\label{Q1}
\ee

Moreover, as in the proofs of Theorem \ref{Th1} and Lemma \ref{Le00},
by {\bf (H1)$''$}, {\bf (H2)$'$} and {\bf (H4)$'$},
for some $\beta_3>1$ and $p\geq 2(\beta_3^*:=\beta_3/(\beta_3-1)$),
we  have that for any $0\leq t<t'\leq T$
$$
\sup_{\eps\in(0,1)}\mE\|(X^\eps(t')
-X^\eps(t))\cdot 1_{\{t',t\leq\tau^\eps_n\}}\|^p_{\mX}\leq C_{T,p,n}\Big(|t-t'|^{\d p}
+|t-t'|^{\frac{\g p}{2}}+|t-t'|^{\frac{p}{2\beta^*_3}}\Big).
$$
Thus, by (v) of Proposition \ref{Pr1} and (\ref{Q1}),
for any $\a_0\in(0,\a)$ and $p$ large enough we have
\ce
 &&\sup_{\eps\in(0,1)}\mE\|(X^\eps(t')
-X^\eps(t))\cdot 1_{\{t',t\leq T\wedge\tau^\eps_n\}}\|^p_{\mX_{\a_0}}\\
&&\qquad\leq C_{N,n,T,p,\kappa_{\a,n},\a_0}\Big(|t-t'|^{\d p}
+|t-t'|^{\frac{\g p}{2}}+|t-t'|^{\frac{p}{2\beta^*}}\Big)^{1-\frac{\a_0}{\a}}.
\de
The desired estimate now follows by Theorem \ref{KCC}.
\end{proof}

In order to obtain the tightness of the laws of
$\{X^\eps,\eps\in(0,1)\}$
in $\mC_T(\mX)$,  we assume that
\begin{enumerate}[{\bf (C1)}]
\item $\fL^{-1}$ is a  compact operator on $\mX$.
\item
$\lim_{n\rightarrow\infty}\sup_{\eps\in(0,1)}P\{\omega:\tau^\eps_n(\omega)<T\}=0.$
\end{enumerate}
Note that {\bf (C2)} implies
$$
P\{\omega:\tau^\eps(\omega)>T\}=1.
$$

We now prove the following key lemma for the large deviation principle of \textsc{Eq}.(\ref{Eqq1}).
\bl\label{Le3}
Under {\bf (C1)} and {\bf (C2)},
there exist   subsequence $\eps_k\downarrow 0$, a probability space
$(\tilde\Omega,\tilde\cF,\tilde P)$ and a sequence
$\{(\tilde h^{k}, \tilde X^{k},\tilde W^{k})\}_{k\in\mN}$
as well as $(h, X^h, \tilde W)$ defined on this probability
space and taking values in $\mD_N\times \mC_T(\mX)\times \mC_T(\mU)$ such that
\begin{enumerate}[(i)]
\item $(\tilde h^{k}, \tilde X^{k},\tilde W^{k})$ has the same law as
$(h^{\eps_k}, X^{\eps_k},W)$ for each $k\in\mN$;

\item $(\tilde h^{k}, \tilde X^{k},\tilde W^{k})
\rightarrow (h, X^h,\tilde W)$ in $\mD_N \times \mC_T(\mX)\times \mC_T(\mU)$,
$\tilde P$-a.s. as $k\to\infty$;

\item $(h,X^h)$ uniquely solves the following Volterra equation:
\be
X^h(t)=g(t)+\int^t_0A(t,s,X^h(s))\dif s+\int^t_0B(t,s,X^h(s))\dot h(s)\dif s.\label{Eq3}
\ee
\end{enumerate}
In particular, {\bf (LD)$_\mathbf{1}$} in Subsection \ref{Sec2} holds.
\el
\begin{proof}
Let $\a_0\in(0,\a)$ and $a>0$ be as in Lemma \ref{Le2}. For $R>0$, set
$$
K_{R}:=\left\{x\in\mC_T(\mX): \sup_{t\in[0,T]}\|x(t)\|_\mX
+\sup_{s\not=t\in[0,T]}\frac{\|x(t)-x(s)\|_{\mX_{\a_0}}}{|t-s|^a}\leq R\right\}.
$$
By {\bf (C1)}, $\mX_{\a_0}\hookrightarrow\mX$ is compact (cf. \cite[p.29, Theorem 1.4.8]{He}).
Thus, by Ascoli-Arzel\`a's theorem  (cf. \cite{Ka}),
the set $K_{R}$ is compact in $\mC_T(\mX)$.
For any $\d>0$, by {\bf (C2)} we can choose $n$ sufficiently large such that
\ce
\sup_{\eps\in(0,1)}P\Big\{\omega:\tau^{\eps}_n(\omega)<T\Big\}\leq \d.
\de
By Lemma \ref{Le2} and Chebyschev's inequality, for any $R>n$ we have
\ce
P\{X^{\eps}(\cdot)\notin K_{R}\}&=&P\{X^{\eps}(\cdot)\notin K_{R},\tau^{\eps}_n\geq T\}
+P\{X^{\eps}(\cdot)\notin K_{R}, \tau^{\eps}_n<T\}\\
&\leq&P\left\{\sup_{s\not=t\in[0,T\wedge\tau^{\eps}_n]}\frac{\|X^{\eps}(t)
-X^{\eps}(s)\|_{\mX_{\a_0}}}{|t-s|^a}\geq R-n\right\}
+P\{\tau^{\eps}_n<T\}\\
&\leq&\mE\left[\sup_{s\not=t\in[0,T\wedge\tau^{\eps}_n]}\frac{\|X^{\eps}(t)
-X^{\eps}(s)\|^p_{\mX_{\a_0}}}{|t-s|^{ap}}\right]/(R-n)^p+\d\\
&\leq&C_{N,n,T,p,\kappa_{\a,n},\a_0}/(R-n)^p+\eps'.
\de
Therefore, for $R$ large enough we have
$$
\sup_{\eps\in(0,1)}P\{X^{\eps}(\cdot)\notin K_{R}\}\leq 2\d.
$$
Thus, by the compactness of $\mD_N$ (see (\ref{Metr})),
the laws of  $(h^{\eps}, X^{\eps}, W)$
in $\mD_N\times \mC_T(\mX)\times \mC_T(\mU)$ is tight.
By Skorohod's embedding theorem (cf. \cite{Ka}),
the conclusions (i) and (ii) hold.

We now prove (iii). Note that by (i)  (cf. \cite[Section 8]{On})
\ce
\tilde X^{k}(t)&=&g(t)+\int^t_0A(t,s,\tilde X^{k}(s))\dif s+
\int^t_0B(t,s,\tilde X^{k}(s))\dot{ \tilde{h}}^{k}(s)\dif s\\
&&+\sqrt{\eps_k}\int^t_0B(t,s,\tilde X^{k}(s))\dif \tilde W^{k}(s)\\
&=:&g(t)+J^{k}_1(t)+J^{k}_2(t)+J^{k}_3(t),\ \ \tilde P-a.s..
\de
Set
$$
\tilde\tau^{k}_n:=\inf\{t\geq 0: \|\tilde X^{k}(t)\|_\mX>n\}.
$$
Then for any $\d>0$, by (i) and {\bf (C2)}
there exists an $n$ large enough such that
\ce
\sup_{k\in\mN}\tilde P\{\tilde\tau^{k}_n<T\}
&=&\sup_{k\in\mN}\tilde P\left\{\sup_{s\in[0,T)}\|\tilde X^{k}(s)\|_\mX>n\right\}\\
&=&\sup_{k\in\mN}P\left\{\sup_{s\in[0,T)}\| X^{\eps_k}(s)\|_\mX>n\right\}\\
&=&\sup_{k\in\mN}P\{\tau^{\eps_k}_n<T\}\leq\d.
\de
Hence, for any $\d'>0$, by the BDG inequality (\ref{BDG})
and {\bf(H2)$'$} we have
\ce
\tilde P\Big\{\|J^{k}_3(t)\|_\mX\geq \d'\Big\}&\leq&
\tilde P\Big\{J^{k}_3(t)\geq \d'; \tilde\tau^{k}_n\geq T\Big\}
+\tilde P\Big\{\tilde\tau^{k}_n<T\Big\}\\
&\leq& \frac{\mE^{\tilde P}\|J^{k}_3(t)\cdot 1_{\{t\leq\tilde\tau^{k}_n\}}\|^2_\mX}{ \d'^2}
+\d\\
&\leq& \frac{\eps_k \cdot C_n\mE^{\tilde P}\left(\int^{t\wedge\tilde\tau^{k}_n}_0
\kappa_{1,n}(t,s)\dif s\right)}{ \d'^2}+\d\\
&\leq&\frac{\eps_k\cdot C_{n,t}}{ \d'^2}+\d.
\de
Thus, we get
$$
\lim_{k\to\infty}\tilde P\Big\{\|J^{k}_3(t)\|_\mX\geq \d'\Big\}=0.
$$

Let $J_i(t), i=1,2$ be the corresponding terms in  \textsc{Eq.}(\ref{Eq3}).
In order to prove that $X^h$ solves \textsc{Eq.}(\ref{Eq3}), it is now enough to show
that for any $t\in[0,T]$ and $y\in\mX^*$
$$
\lim_{k\to\infty}{}_\mX\<J^{k}_i(t)-J_i(t),y\>_{\mX^*}=0,\ \ i=1,2,\ \ \tilde P-a.s..
$$

Observe that
\ce
|{}_\mX\<J^{k}_2(t)-J_2(t),y\>_{\mX^*}|
&\leq&\|y\|_{\mX^*}\cdot\int^t_0\|[B(t,s,\tilde X^{k}(s))-B(t,s,X^h(s))]
\dot{ \tilde{h}}^{k}(s)\|_\mX\dif s\\
&&+\left|\int^t_0{}_\mX\<B(t,s,X^h(s))[\dot{ \tilde{h}}^{k}(s)
-\dot h(s)],y\>_{\mX^*}\dif s\right|\\
&=:&\|y\|_{\mX^*}\cdot J^{k}_{21}(t)+J^{k}_{22}(t).
\de
By the weak convergence of $\tilde{h}^{k}$ to $h$ in $\mD_N$, we have
$$
\lim_{k\to\infty}J^{k}_{22}(t)=0.
$$
Noting that by (ii), for almost all $\tilde\omega\in\tilde\Omega$ and some $K(\tilde\omega)\in\mN$
$$
n(\tilde\omega):=\sup_{s\in[0,T]}\|X^h(s,\tilde\omega)\|_\mX
\vee\sup_{k\geq K(\tilde\omega)}\sup_{s\in[0,T]}\|\tilde X^{k}(s,\tilde\omega)\|_\mX<+\infty,
$$
we have, by H\"older's inequality and {\bf (H3)$'$}
\ce
J^{k}_{21}(t,\tilde\omega)&\leq&\|\tilde{h}^{k}(\tilde\omega)\|_{\ell^2_T}
\cdot\left(\int^t_0\big\|B(t,s,\tilde X^{k}(s,\tilde\omega))
-B(t,s,X^h(s,\tilde\omega))\big\|^2_{L_2(l^2;\mX)}\dif s\right)^{1/2}\\
&\leq&N\cdot\left(\int^t_0\kappa_{2,n(\tilde\omega)}(t,s)\cdot\|\tilde X^{k}(s,\tilde\omega)
-X^h(s,\tilde\omega)\|^2_\mX\dif s\right)^{1/2}\\
&\stackrel{(ii)}{\to}&0,\ \ \mbox{ as $k\to\infty$},
\de
where we have used $\tilde h^{k}(\tilde\omega)\in \mD_N$.

Similarly, we have
$$
\lim_{k\to\infty}\|J^{k}_1(t)-J_1(t)\|_{\mX}=0, \ \tilde P-a.s..
$$
Combining the above estimates, we find that $X^h$ solves \textsc{Eq.}(\ref{Eq3}).
\end{proof}

Let $I(f)$ be defined by
\be
I(f):=\frac{1}{2}\inf_{\{h\in\ell^2_T:~f=X^h\}}\|h\|^2_{\ell^2_T},\ \ f\in\mC_T(\mX),\label{Rate}
\ee
where $X^h$ is defined by \textsc{Eq.}(\ref{Eq3}).
In order to identify $I(f)$, we assume that
\begin{enumerate}[{\bf(C3)}]
\item For any $N\in\mN$
$$
\sup_{h\in \mD_N}\sup_{t\in[0,T]}\|X^h(t)\|_\mX<+\infty.
$$
\end{enumerate}
Similar to the proof of Lemma \ref{Le3}, we can prove that:
\bl\label{Le4} Under {\bf (C3)},
{\bf (LD)$_\mathbf{2}$} in Subsection \ref{Sec2} holds.
\el

Thus, by Theorem \ref{Th2} we have proven:
\bt\label{Large}
Assume that {\bf (H1)$''$}-{\bf (H2)$''$}, {\bf (H2)$'$}-{\bf (H4)$'$}
and {\bf (C1)}-{\bf (C3)} hold.
Then, $\{X_\eps,\eps\in(0,1)\}$ satisfies the large deviation principle in $\mC_T(\mX)$
with the rate function $I(f)$ given by (\ref{Rate}).
\et

\br\label{R4}
The conditions {\bf(C2)} and {\bf(C3)} are satisfied
if {\bf (H1)$''$}, {\bf (H2)} and {\bf (H4)} hold,
and $\kappa_1$ in  {\bf (H2)} belongs to $\sK_{>1}$.
In fact, we can prove as the proof of Theorem \ref{Th4}
$$
\sup_{n\in\mN}\sup_{\eps\in(0,1)}\mE\left(\sup_{t \in[0,T\wedge\tau^\eps_n]}
\|X^{\eps}(t)\|^p_{\mX}\right)\leq C_{T,p,\kappa_{1}},
$$
which then implies  {\bf(C2)}.
The condition {\bf(C3)} is more direct in this case.
\er

\section{Semilinear stochastic evolutionary integral equations}

In this section, we consider the following semilinear stochastic evolutionary integral equation:
\be
X(t)=x_0-\int^t_0a(t-s)\fL X(s)\dif s+\int^t_0\Phi(s,X(s))\dif s
+\int^t_0\Psi(s,X(s))\dif W(s),\label{Eqe7}
\ee
where $a:\mR_+\to\mR_+$ is a measurable function, and
$$
\Phi: \mR_+\times\Omega\times\mX\to\mX\in\cM\times\cB(\mX)/\cB(\mX)
$$
and
$$
\Psi: \mR_+\times\Omega\times\mX\to L_2(l^2;\mX)
\in\cM\times\cB(\mX)/\cB(L_2(l^2;\mX)).
$$
Here and below, $\cM$ stands for the progressively
measurable $\sigma$-algebra over $\mR_+\times\Omega$.

Consider first the following deterministic integral equation:
\be
x(t)=x_0-\int^t_0 a(t-s)\fL x(s)\dif s.\label{Res}
\ee
The solution of this equation is called the resolvent of $(a,\fL)$,
and denoted by $\fS_t x_0=x(t)$.
Note that in general
$$
\fS_{t+s}\not=\fS_t\circ\fS_s.
$$

We make the following assumptions:
\begin{enumerate}[{\bf (S1)}]
\item
The resolvent $\{\fS_t:t\geq 0\}$ is of analyticity type $(\omega_0,\theta_0)$
in the sense of \cite[Definition 2.1]{Pr}, where $\omega_0\in\mR$ and $\theta_0\in(0,\pi/2]$.
\item For any $R>0$, there exist $C_R>0$ and $\beta\in[0,1)$
such that for all $s>0$, $\omega\in\Omega$ and $x,y\in\mX$ with $\|x\|_\mX,\|y\|_\mX\leq R$
$$
\|\Phi(s,\omega,x)\|_{\mX}+\|\Psi(s,\omega,x)\|^2_{L_2(l^2;\mX)}
\leq \frac{C_R}{(s\wedge 1)^\beta},
$$
and
\ce
\|\Phi(s,\omega,x)-\Phi(s,\omega,y)\|_{\mX}&\leq& \frac{C_R}{(s\wedge 1)^\beta}\|x-y\|_{\mX},\\
\|\Psi(s,\omega,x)-\Psi(s,\omega,y)\|^2_{L_2(l^2;\mX)}&\leq&
\frac{C_R}{(s\wedge 1)^\beta}\|x-y\|_{\mX}^2.
\de
\item For all $s>0$, $\omega\in\Omega$ and $x\in\mX$, it holds that
\ce
\|\Phi(s,\omega,x)\|_{\mX}&\leq& \frac{C}{(s\wedge 1)^\beta}(1+\|x\|_\mX),\\
\|\Psi(s,\omega,x)\|^2_{L_2(l^2;\mX)}&\leq& \frac{C}{(s\wedge 1)^\beta}(1+\|x\|^2_\mX).
\de
\end{enumerate}

The following property of analytic resolvent $\{\fS_t:t>0\}$ is crucial for the proof of
Theorem \ref{Th3} below  (cf. \cite[Corollary 2.1]{Pr}).
\bp
Let $\fS_t$ be an analytic resolvent of type $(\omega_0,\theta_0)$.
Then for any $T>0$
\be
\sup_{t\in[0,T]}\|\fS_t\|_{L(\mX;\mX)}\leq C_T\label{PO1}
\ee
and for any $t\in(0,T]$
\be
\|\dot\fS_t\|_{L(\mX;\mX)}\leq C_T t^{-1}, \label{Po}
\ee
where the dot denotes the operator derivative and $\|\cdot\|_{L(\mX;\mX)}$ denotes
the norm of bounded linear operators.
\ep

By a solution of \textsc{Eq.}(\ref{Eqe7}) we mean that $X(t)$ satisfies the
following stochastic Volterra equation:
\be
X(t)=\fS_t x_0+\int^t_0 \fS_{t-s}\Phi(s,X(s))\dif s
+\int^t_0\fS_{t-s}\Psi(s,X(s))\dif W(s).\label{Eqe4}
\ee
Let us define
$$
A(t,s,\omega,x):=\fS_{t-s}\Phi(s,\omega,x),\ \
B(t,s,\omega,x):=\fS_{t-s}\Psi(s,\omega,x).
$$
We have:
\bt\label{Th3}
Under {\bf (S1)} and {\bf (S2)}, there exists a unique maximal solution
$(X,\tau)$ for \textsc{Eq.} (\ref{Eqe4})
in the sense of Definition \ref{Def1}.  Moreover, if {\bf (S3)} holds,
then $\tau=+\infty$, a.s..
\et
\begin{proof}
First of all, it is easy to see by (\ref{PO1}) that {\bf (H2)$'$} and {\bf (H3)$'$} hold with
$$
\kappa_{1,R}(t,s)=\kappa_{2,R}(t,s)=\frac{C_R}{(s\wedge 1)^\beta}\in\sK_{>1}.
$$
For $0\leq s<t<t'$, $\omega\in\Omega$ and $x\in\mX$ with $\|x\|_\mX\leq R$,
we have
\ce
\|A(t',s,\omega,x)-A(t,s,\omega,x)\|_\mX&=&
\|(\fS_{t'-s}-\fS_{t-s})\Phi(s,\omega,x)\|_\mX\\
&\leq&\frac{C_R}{(s\wedge 1)^\beta}\|\fS_{t'-s}-\fS_{t-s}\|_{L(\mX;\mX)}\\
&\leq&\frac{C_R}{(s\wedge 1)^\beta}\int^{t'-s}_{t-s}\|\dot\fS_{r}\|_{L(\mX;\mX)}\dif r\\
&\stackrel{(\ref{Po})}{\leq}&\frac{C_R}{(s\wedge 1)^\beta}\int^{t'-s}_{t-s}\frac{1}{r}\dif r\\
&=&\frac{C_R}{(s\wedge 1)^\beta}\log\left(\frac{t'-s}{t-s}\right)
\de
and
$$
\|B(t',s,\omega,x)-B(t,s,\omega,x)\|^2_{L_2(l^2;\mX)}\leq
\frac{C_R}{(s\wedge 1)^\beta}\log^2\left(\frac{t'-s}{t-s}\right).
$$
Note that the following elementary inequality holds for any $\gamma\in(0,1)$
$$
\log(1+s)\leq C s^\gamma, \ \ \forall s>0.
$$
Therefore, for $0\leq s<t<t'$, $\omega\in\Omega$ and $x\in\mX$ with $\|x\|_\mX\leq R$
\ce
&&\|A(t',s,\omega,x)-A(t,s,\omega,x)\|_\mX
+\|B(t',s,\omega,x)-B(t,s,\omega,x)\|^2_{L_2(l^2;\mX)}\\
&&\qquad\leq \frac{C_R(t'-t)^\gamma}{(s\wedge 1)^\beta(t-s)^{\gamma}}\left[1+
\frac{(t'-t)^\gamma}{(t-s)^{\gamma}}\right]
=:\l_R(t',t,s).
\de
Thus, we find that {\bf (H4)$'$} holds if $\gamma\in(0,(1-\beta)/2)$.

Lastly, if {\bf (S3)} is satisfied, it is clear that {\bf (H2)} holds
with $\kappa_1(t,s)=\frac{C}{(s\wedge 1)^\beta}\in\sK_{>1}$, and {\bf (H4)} also holds from the
above calculations.
The non-explosion now follows from Theorem \ref{Th4}.
\end{proof}

We now turn to the small perturbation of \textsc{Eq.}(\ref{Eqe4}) and assume that $\Phi$ and
$\Psi$ are non-random. Consider
$$
X_\eps(t)=\fS_t x_0+\int^t_0 \fS_{t-s}\Phi(s,X_\eps(s))\dif s
+\sqrt{\eps}\int^t_0\fS_{t-s}\Psi(s,X_\eps(s))\dif W(s).
$$
In order to use Theorem \ref{Large} to get the LDP for $\{X_\eps,\eps\in(0,1)\}$, we also assume
\begin{enumerate}[{\bf (S4)}]
\item Let $\{\fS_t:t\geq 0\}$ be an analytic resolvent of type $(\omega_0,\theta_0)$. Assume that
for some $\omega_1>\omega_0$, $0<\theta_1<\theta_0$, $C>0$ and $\a_1>0$
\be
|\hat a(\lambda)|\geq C(|\lambda-\omega_1|^{\a_1}+1)^{-1},\ \ \forall \lambda\in\mC
\mbox{ with }|\mathrm{arg}(\lambda-\omega)|<\theta_1,\label{Li1}
\ee
where $\hat a$ denotes the Laplace transform of $a$. Moreover, we also assume that
\be
\int^r_0 a(s)\dif s+\int^t_0|a(r+s)-a(s)|\dif s\leq C_T|r|^\d,\label{Li2}
\ee
where $r,t\in[0,T]$ and $T,\d>0$.
\end{enumerate}
We have
\bt
Under {\bf (S1)}-{\bf (S4)} and {\bf (C1)}, for any $x_0\in\sD(\fL)$,
$\{X_\eps,\eps\in(0,1)\}$ satisfies the large deviation principle
in $\mC_T(\mX)$ with the rate function $I(f)$ given by (\ref{Rate}).
\et
\begin{proof}
From the proof of Theorem \ref{Th3}, it is enough to check {\bf (H1)$''$}
and {\bf (H2)$''$}. By (\ref{Li1}) and \cite[p.57, Theorem 2.2 (ii)]{Pr}, we have
$$
\|\fL\fS_t\|_{L(\mX;\mX)}\leq Ce^{\omega_1 t}(1+t^{-\a_1}),\ \ \forall t>0,
$$
which together with (v) of Proposition \ref{Pr1} yields that
for any $\a\in(0,1)$ and $T>0$
$$
\|\fL^\a\fS_t\|_{L(\mX;\mX)}\leq C_T(1+t^{-\a_1\cdot\a}),\ \ \forall t\in(0,T].
$$
Thus, {\bf (H2)$''$} holds by choosing $\a<\frac{1-\beta}{\a_1}$, where $\beta$
is from {\bf (S3)}.

For {\bf (H1)$''$}, since $x_0\in\sD(\fL)=\mX_1$, by (\ref{PO1}) we have
$$
\|\fL\fS_t x_0\|_\mX=\|\fS_t \fL x_0\|_\mX\leq C\|\fL x_0\|_\mX.
$$
On the other hand, by the resolvent equation (\ref{Res})
and (\ref{Li2}) we have, for any $0\leq t<t'\leq T$
\ce
\|\fS_{t'}x_0-\fS_{t}x_0\|_\mX&\leq&\int^t_0|a(t'-s)-a(t-s)|\cdot\|\fL\fS_{s}x_0\|_\mX\dif s\\
&&+\int^{t'}_t |a(t'-s)|\cdot\|\fL\fS_{s}x_0\|_\mX\dif s\\
&\leq& C_T\|\fL x_0\|_\mX\cdot|t'-t|^\d.
\de
The proof is thus completed by Theorem \ref{Large} and Remark \ref{R4}.
\end{proof}

\bx{\rm
Let $a$ be a completely monotonic kernel function, i.e.,
\be
a(t)=\int^\infty_0 e^{-st}\dif \rho(s),\ \ t>0,\label{Com}
\ee
where $s\mapsto\rho(s)$ is nondecreasing, and such that $\int^\infty_1\dif \rho(s)/s<\infty$.
Then the resolvent $\{\fS_t:t\geq 0\}$ associated with $a$ is of analyticity type $(0,\theta)$
for some $\theta\in(0,\pi/2)$ (cf. \cite[p.55, Example 2.2]{Pr}), i.e., {\bf (S1)} holds.
For {\bf (S4)}, besides (\ref{Com}) and (\ref{Li2}), we also assume that for some $C,\a_1>0$
\be
C(1+\l)^{-\a_1}\leq \int^\infty_0 e^{-\l t}\cdot a(t) \dif t<+\infty,\ \ \forall \l>0,\label{Sub}
\ee
which implies by \cite[p.221, Lemma 8.1 (v)]{Pr} that (\ref{Li1}) holds.
In particular,
$$
a_\a(t)=\frac{t^{\a-1}}{\Gamma(\a)},\ \ \a\in(0,1]
$$
is completely monotonic, and satisfies (\ref{Li2}) and (\ref{Sub}),
where $\Gamma$ denotes the usual Gamma function.

Moreover, for the kernel function $a_\a$, if
$$
1<\a<2-\frac{2\phi}{\pi}<2,
$$
where $\phi$ comes from (\ref{Sec}), then
$\fS_t$ is analytic  (cf. \cite[p.55, Example 2.1]{Pr}). Notice that in
\cite{Pr}, $-\fL$ is considered. In this case, (\ref{Li1}) and (\ref{Li2}) clearly hold since
$\hat a_\a(\l)=\l^{-\a}$, $\mathrm{Re}\l>0$.
}
\ex

\section{Semilinear stochastic partial differential equations}

When $a=1$ in \textsc{Eq.}(\ref{Eqe7}), one sees that
\textsc{Eq.}(\ref{Eqe7}) contains a class of semilinear SPDEs.
However, it cannot deal with the equation
like stochastic Navier-Stokes equation. In this section, we shall discuss
 strong solutions of a large class of semilinear SPDEs by using the properties of
analytic semigroups.

\subsection{Mild solutions of SPDEs driven by Brownian motions}
Consider the following semilinear stochastic partial differential equation:
\be
\dif X(t)=[-\fL X(t)+\Phi(t,X(t))]\dif t+\Psi(t,X(t))\dif W(t),\ \ X(0)=x_0.\label{Eq7}
\ee
We study two cases, in application, which correspond to different types of SPDEs.
First of all, we introduce the following assumptions on the coefficients:
\begin{enumerate}[{\bf (M1)}]
\item For some $\a\in(0,1)$
$$
\Phi: \mR_+\times\Omega\times\mX_{\a}\to\mX\in\cM\times\cB(\mX_{\a})/\cB(\mX)
$$
and
$$
\Psi: \mR_+\times\Omega\times\mX_{\a}\to L_2(l^2;\mX_{\frac{\a}{2}})
\in\cM\times\cB(\mX_{\a})/\cB(L_2(l^2;\mX_{\frac{\a}{2}})).
$$
\item For any $R>0$, there exist $C_R>0$ and $\beta\in[0,1)$ with
$$
\a+\beta<1
$$
such that for all $s>0$, $\omega\in\Omega$ and $x,y\in\mX_\a$ with
$\|x\|_{\mX_\a},\|y\|_{\mX_\a}\leq R$,
$$
\|\Phi(s,\omega,x)\|_{\mX}+\|\Psi(s,\omega,x)\|^2_{L_2(l^2;\mX_{\frac{\a}{2}})}
\leq \frac{C_R}{(s\wedge 1)^\beta}
$$
and
\ce
\|\Phi(s,\omega,x)-\Phi(s,\omega,y)\|_{\mX}&\leq& \frac{C_R}{(s\wedge 1)^\beta}\|x-y\|_{\mX_\a},\\
\|\Psi(s,\omega,x)-\Psi(s,\omega,y)\|^2_{L_2(l^2;\mX_{\frac{\a}{2}})}&\leq&
\frac{C_R}{(s\wedge 1)^\beta}\|x-y\|_{\mX_\a}^2.
\de
\item For all $s>0$, $\omega\in\Omega$ and $x\in\mX_\a$, it holds that
\ce
\|\Phi(s,\omega,x)\|_{\mX}
&\leq& \frac{C}{(s\wedge 1)^\beta}(1+\|x\|_{\mX_\a}),\\
\|\Psi(s,\omega,x)\|^2_{L_2(l^2;\mX_{\frac{\a}{2}})}
&\leq& \frac{C}{(s\wedge 1)^\beta}(1+\|x\|_{\mX_\a}^2).
\de
\end{enumerate}
By a mild solution of  equation (\ref{Eq7})
we mean that $X(t)$ solves the following stochastic Volterra integral equation:
\be
X(t)=\fT_t x_0+\int^t_0 \fT_{t-s}\Phi(s,X(s))\dif s
+\int^t_0\fT_{t-s}\Psi(s,X(s))\dif W(s).\label{Eq4}
\ee
\bt\label{Th5}
Under {\bf (M1)} and {\bf (M2)}, for any $x_0\in\mX_\a$ ($\a$ is from {\bf (M1)}),
there exists a unique maximal solution
$(X,\tau)$ for \textsc{Eq.}(\ref{Eq4}) so that
\begin{enumerate}[(i)]
\item $t\mapsto X(t)\in\mX_\a$ is continuous on $[0,\tau)$ almost surely;
\item $\lim_{t\uparrow\tau}\|X(t)\|_{\mX_\a}=+\infty$ on $\{\omega:\tau(\omega)<+\infty\}$;
\item it holds that, $P$-a.s, on $[0,\tau)$
\ce
X(t)=\fT_{t} x_0+\int^{t}_0\fT_{t-s}\Phi(s,X(s))\dif s
+\int^{t}_0\fT_{t-s}\Psi(s,X(s))\dif W(s).
\de
\end{enumerate}
Moreover, if {\bf (M3)} holds, then $\tau=+\infty$, a.s..
\et
\begin{proof}
We first consider the following stochastic Volterra integral equation
\be
Y(t)&=&\fL^\a\fT_t x_0+\int^t_0 \fL^\a\fT_{t-s}\Phi(s,\fL^{-\a}Y(s))\dif s\no\\
&&+\int^t_0\fL^\a\fT_{t-s}\Psi(s,\fL^{-\a}Y(s))\dif W(s).\label{Eqq4}
\ee
Define
\ce
g(t)&:=&\fL^\a\fT_t x_0,\\
A(t,s,\omega,y)&:=&\fL^\a\fT_{t-s}\Phi(s,\omega,\fL^{-\a}y),\\
B(t,s,\omega,y)&:=&\fL^\a\fT_{t-s}\Psi(s,\omega,\fL^{-\a}y).
\de
Let us verify {\bf (H1)$'$}-{\bf (H4)$'$}. Clearly, {\bf (H1)$'$} holds since $x_0\in\mX_\a$.

By (iii) of Proposition \ref{Pr1} and {\bf (M2)},
for all $t>s>0$, $\omega\in\Omega$ and
$x,y\in\mX$ with $\|x\|_\mX,\|y\|_\mX\leq R$ we have
\be
&&\|A(t,s,\omega,x)\|_\mX+\|B(t,s,\omega,x)\|^2_{L_2(l^2;\mX)}\no\\
&&\qquad\preceq\frac{1}{(t-s)^\a}\Big(\|\Phi(s,\omega,\fL^{-\a}x)\|_{\mX}
+\|\Psi(s,\omega,\fL^{-\a}x)\|^2_{L_2(l^2;\mX_{\frac{\a}{2}})}\Big)\no\\
&&\qquad\leq\frac{C_R}{(t-s)^\a (s\wedge 1)^\beta},\label{Ep1}
\ee
and
\ce
\|A(t,s,\omega,x)-A(t,s,\omega,y)\|_\mX
&\preceq&\frac{1}{(t-s)^\a}\|\Phi(s,\omega,\fL^{-\a}x)-\Phi(s,\omega,\fL^{-\a}y)\|_{\mX}\\
&\leq&\frac{C_R}{(t-s)^\a(s\wedge 1)^\beta}\|\fL^{-\a}x-\fL^{-\a}y\|_{\mX_\a}\\
&=&\frac{C_R}{(t-s)^\a (s\wedge 1)^\beta}\|x-y\|_{\mX}
\de
as well as
\ce
&&\|B(t,s,\omega,x)-B(t,s,\omega,y)\|^2_{L_2(l^2;\mX)}\\
&&\qquad\preceq\frac{1}{(t-s)^\a}\|\Psi(s,\omega,\fL^{-\a}x)
-\Psi(s,\omega,\fL^{-\a}y)\|_{L_2(l^2;\mX_{\frac{\a}{2}})}^2\\
&&\qquad\leq\frac{C_R}{(t-s)^\a (s\wedge 1)^\beta}\|x-y\|^2_{\mX}.
\de
Hence, if we  take
$$
\kappa_{1,R}(t,s)=\kappa_{2,R}(t,s):=\frac{C_R}{(t-s)^\a (s\wedge 1)^\beta}\in\sK_{>1},
$$
then {\bf (H2)$'$} and {\bf (H3)$'$} hold.

Let $0<\gamma<1-(\a+\beta)$. By (iv) of Proposition \ref{Pr1} and {\bf (M2)} we  have
\ce
\|A(t',s,\omega,x)-A(t,s,\omega,x)\|_\mX&=&\|(\fT_{t'-t}-1)\fL^\a\fT_{t-s}
\Phi(s,\omega,\fL^{-\a}x)\|_\mX\\
&\preceq&(t'-t)^\g\|\fL^{\a+\g}\fT_{t-s}\Phi(s,\omega,\fL^{-\a}x)\|_{\mX}\\
&\leq&\frac{C_R(t'-t)^{\g}}{(t-s)^{\a+\g} (s\wedge 1)^\beta}
\de
and
\ce
&&\|B(t',s,\omega,x)-B(t,s,\omega,x)\|^2_{L_2(l^2;\mX)}\\
&&\qquad\preceq\|(\fT_{t'-t}-1)\fL^\a\fT_{t-s}\Psi(s,\omega,\fL^{-\a}x)\|_{L_2(l^2;\mX)}^2\\
&&\qquad\preceq(t'-t)^\g\|\fL^{\a+\g/2}\fT_{t-s}\Psi(s,\fL^{-\a}x)\|^2_{L_2(l^2;\mX)}\\
&&\qquad\preceq\frac{(t'-t)^\g}{(t-s)^{\a+\g}}
\|\fL^{\frac{\a}{2}}\Psi(s,\fL^{-\a}x)\|^2_{L_2(l^2;\mX)}\\
&&\qquad\leq\frac{C_R(t'-t)^{\g}}{(t-s)^{\a+\g} (s\wedge 1)^\beta}.
\de
So, if we take
$$
\lambda_R(t',t,s):=\frac{C_R(t'-t)^{\g}}{(t-s)^{\a+\g} (s\wedge 1)^\beta},
$$
then {\bf (H4)$'$} holds.

Hence, by Theorem \ref{Main} there is a unique maximal solution $(Y,\tau)$ for \textsc{Eq.}(\ref{Eqq4})
in the sense of Definition \ref{Def1}. Set
$$
X(t)=\fL^{-\a}Y(t).
$$
It is easy to see that $(X,\tau)$ a unique maximal solution for \textsc{Eq.}(\ref{Eq4}),
which satisfies (i), (ii) and (iii) in the theorem.

Lastly, if {\bf (M3)} is satisfied, then  as estimating (\ref{Ep1}),
for the above $A$ and $B$, {\bf (H2)} holds with some $\kappa_1\in\sK_{>1}$, and
also  {\bf (H4)} holds.
So, by Theorem \ref{Th4} we have $\tau=\infty$ a.s..
\end{proof}

\br
The solution $(X,\tau)$ in Theorem \ref{Th5} is clearly a local solution of \textsc{Eq.}(\ref{Eq4}) in $\mX$.
However, it may be not a maximal solution in $\mX$ because it may happen that
$$
\lim_{t\uparrow\tau(\omega)}\|X(t,\omega)\|_\mX<+\infty \mbox{ on $\{\omega: \tau(\omega)<+\infty\}$}.
$$
\er

Next, we study the large deviation estimate for \textsc{Eq.}(\ref{Eq7}), and
assume that $\Phi$ and $\Psi$ are non-random.
Consider the following small perturbation of \textsc{Eq.}(\ref{Eq7}):
\be
\dif X_\eps(t)=[-\fL X_\eps(t)+\Phi(t,X_\eps(t))]\dif t+\sqrt{\eps}
\Psi(t,X_\eps(t))\dif W(t),\ \ X_\eps(0)=x_0.\label{Eqq7}
\ee

In order to apply Theorem \ref{Large} to this situation, we need the non-explosion
assumptions as {\bf (C2)} and {\bf (C3)}.
For a family of processes $\{h^\eps,\eps\in(0,1)\}$ in $\cA^T_N$
(see (\ref{Op2}) for the definition of $\cA^T_N$),  consider
\ce
X^\eps(t)&=&\fT_t x_0+\int^t_0\fT_{t-s}\Phi(s,X^\eps(s))\dif s
+\int^t_0\fT_{t-s}\Psi(s,X^\eps(s))\dot h^\eps(s)\dif s\no\\
&&+\sqrt{\eps}\int^t_0\fT_{t-s}\Psi(s,X^\eps(s))\dif W(s),
\de
and for $h\in\ell^2_T$ (see (\ref{H2}))
$$
X^h(t)=\fT_t x_0+\int^t_0\fT_{t-s}\Phi(s,X^h(s))\dif s
+\int^t_0\fT_{t-s}\Psi(s,X^h(s))\dot h(s)\dif s.
$$
Below, for $n\in\mN$ we define
$$
\tau^\eps_n:=\inf\{t>0: \|X^\eps(t)\|_{\mX_\a}>n\}.
$$
Our large deviation principle can be stated as follows:
\bt\label{Th6}
Assume {\bf (M1)} and {\bf (M2)}.
Let $x_0\in\mX_{\d}$ for some $1\geq\d>\a$, where $\a$ is from {\bf (M1)}.
We also assume that $\sD(\fL)=\mX_1\subset\mX$ is compact,
and
\be
\lim_{n\rightarrow\infty}\sup_{\eps\in(0,1)}
P\big\{\omega:\tau^\eps_n(\omega)<T\big\}=0\label{Au1}
\ee
and for any $N>0$
\be
\sup_{h\in \mD_N}\sup_{t\in[0,T]}\|X^h(t)\|_{\mX_\a}<+\infty.\label{Au2}
\ee
Then, $\{X_\eps,\eps\in(0,1)\}$ satisfies the large deviation principle in $\mC_T(\mX_\a)$
with the rate function $I(f)$ given by
\be
I(f):=\frac{1}{2}\inf_{\{h\in\ell^2_T:~f=X^h\}}\|h\|^2_{\ell^2_T},\ \ f\in\mC_T(\mX_\a).\label{Rate1}
\ee
\et
\begin{proof}
By Theorem \ref{Large}, it only need to check {\bf(H1)$''$} and {\bf(H2)$''$} for \textsc{Eq.}(\ref{Eqq4}).
Since $x_0\in\mX^\d$ with $\d>\a$,
by (iv) of Proposition \ref{Pr1}, {\bf(H1)$''$} holds with $\d'=\d-\a$ and
$\a'\in(0,\d-\a)$. As the calculations given in (\ref{Ep1}), one finds that
{\bf(H2)$''$} holds with $\a'\in(0,1-\a-\beta)$.
\end{proof}
\br
If {\bf (M3)} is satisfied, one can see that
(\ref{Au1}) and (\ref{Au2}) hold by Remark \ref{R4}.
\er

We now consider another group of assumptions on the coefficients:
\begin{enumerate}[{\bf (M1)$'$}]
\item For some $\a\in(0,1)$
$$
\Phi: \mR_+\times\Omega\times\mX\to\mX_{-\a}\in\cM\times\cB(\mX)/\cB(\mX_{-\a})
$$
and
$$
\Psi: \mR_+\times\Omega\times\mX\to L_2(l^2;\mX_{-\frac{\a}{2}})
\in\cM\times\cB(\mX)/\cB(L_2(l^2;\mX_{-\frac{\a}{2}})).
$$
\item For any $R>0$, there exist $C_R>0$ and $\beta\in(0,1)$ with
$$
\a+\beta<1
$$
such that for all $s>0$, $\omega\in\Omega$ and $x,y\in\mX$ with $\|x\|_\mX,\|y\|_\mX\leq R$,
$$
\|\Phi(s,\omega,x)\|_{\mX_{-\a}}+\|\Psi(s,\omega,x)\|^2_{L_2(l^2;\mX_{-\frac{\a}{2}})}
\leq \frac{C_R}{(s\wedge 1)^\beta}
$$
and
\ce
\|\Phi(s,\omega,x)-\Phi(s,\omega,y)\|_{\mX_{-\a}}&\leq& \frac{C_R}{(s\wedge 1)^\beta}\|x-y\|_{\mX},\\
\|\Psi(s,\omega,x)-\Psi(s,\omega,y)\|^2_{L_2(l^2;\mX_{-\frac{\a}{2}})}&\leq&
\frac{C_R}{(s\wedge 1)^\beta}\|x-y\|_{\mX}^2.
\de
\item For all $s>0$, $\omega\in\Omega$ and $x\in\mX$,  it holds that
\ce
\|\Phi(s,\omega,x)\|_{\mX_{-\a}}
&\leq& \frac{C}{(s\wedge 1)^\beta}(1+\|x\|_\mX),\\
\|\Psi(s,\omega,x)\|^2_{L_2(l^2;\mX_{-\frac{\a}{2}})}
&\leq& \frac{C}{(s\wedge 1)^\beta}(1+\|x\|^2_\mX).
\de
\end{enumerate}

The following two results are parallel to Theorems \ref{Th5} and \ref{Th6}, we omit the details.
\bt\label{Th9}
Under {\bf (M1)$'$} and {\bf (M2)$'$}, for any $x_0\in\mX$,
there exists a unique maximal mild solution $(X,\tau)$ for \textsc{Eq.} (\ref{Eq4})
in the sense of Definition \ref{Def1}. Moreover, if {\bf (M3)$'$} holds,
then $\tau=+\infty$, a.s..
\et

\bt\label{Th10}
Assume that {\bf (M1)$'$}, {\bf (M2)$'$} and {\bf (C1)}-{\bf (C3)} hold.
Let $x_0\in\mX_{\d}$ for some $\d>0$
Then, $\{X_\eps,\eps\in(0,1)\}$ satisfies the large deviation principle in $\mC_T(\mX)$
with the rate function $I(f)$ given by (\ref{Rate}).
\et
\br
Theorem \ref{Th9} is due to Brze\'zniak \cite{b2}. Compared with Theorem \ref{Th5},
the solution in Theorem \ref{Th5} has better regularity, and is in fact a strong solution
under a slightly stronger assumption {\bf (M4)} below.
\er

\subsection{Strong solutions of SPDEs driven by Brownian motions}

In this subsection, following the method used in the deterministic case (cf. \cite{He,Pa}),
we prove the existence of strong solutions for \textsc{Eq.}(\ref{Eq7}). For this aim, in addition to
{\bf (M1)} and {\bf (M2)} with $\beta=0$, we also assume
\begin{enumerate}[{\bf (M4)}]
\item For any $R,T>0$, there exist $\d>0$ and $\a'>1$ such that
for all $s,s'\in[0,T]$, $\omega\in\Omega$ and $x\in\mX_\a$ with $\|x\|_{\mX_\a}\leq R$
\be
\|\Phi(s',\omega,x)-\Phi(s,\omega,x)\|_{\mX}&\leq& C_{T,R}|s'-s|^\d,\label{Po5}\\
\|\Psi(s,\omega,x)\|^2_{L_2(l^2;\mX_{\frac{\a'}{2}})}&\leq& C_{T,R}.\label{Po6}
\ee
\end{enumerate}
Let us first recall the following result (cf. \cite[Theorem 3.2.2]{He}
or \cite[p.114, Theorem 3.5]{Pa}).
\bl\label{Le6}
Let $[0,T]\ni s\mapsto f(s)\in\mX$ be a H\"older continuous function.
Then
$$
t\mapsto\int^t_0\fT_{t-s}f(s)\dif s\in C([0,T];\mX_1).
$$
\el
Using this lemma, we can prove the  existence of strong solutions for \textsc{Eq.}(\ref{Eq7}).
\bt\label{Th55}
Assume that {\bf (M1)}, {\bf (M2)} and {\bf (M4)} hold.
For any $x_0\in\mX_1$, let $(X,\tau)$ be the unique maximal solution of \textsc{Eq.}(\ref{Eq4})
in Theorem \ref{Th5}. Then
\begin{enumerate}[(i)]
\item $t\mapsto X(t)\in\mX_1$ is continuous on $[0,\tau)$ a.s.;
\item it holds that in $\mX$
\ce
X(t)=x_0-\int^{t}_0\fL X(s)\dif s+\int^{t}_0\Phi(s,X(s))\dif s
+\int^{t}_0\Psi(s,X(s))\dif W(s)
\de
for  all $t\in[0,\tau)$, $P$-a.s..
\end{enumerate}
We shall call $(X,\tau)$  the unique maximal strong solution of \textsc{Eq.}(\ref{Eq7}).
\et
\begin{proof}
For $n\in\mN$, set
$$
\tau_n:=\inf\{t>0: \|X(t)\|_{\mX_\a}>n\}
$$
and
$$
G(t,s):=\fT_{t-s}\Psi(s,X(s)).
$$
Then by  (iii) and (iv) of Proposition \ref{Pr1} we have
$$
\|G(t,s)\|_{L_2(l^2;\mX_1)}^2\preceq \frac{1}{(t-s)^{2-\a'}}
\|\Psi(s,X(s))\|^2_{L_2(l^2;\mX_{\a'/2})},
$$
and in view of $\a'>1$
$$
\|G(t',s)-G(t,s)\|_{L_2(l^2;\mX_1)}^2\preceq
\frac{(t'-t)^{(\a'-1)/2}}{(t-s)^{(3-\a')/2}}\|\Psi(s,X(s))\|^2_{L_2(l^2;\mX_{\a'/2})}.
$$
Hence, by Lemma \ref{Le00} and (\ref{Po6}),
$$
t\mapsto\int^t_0\fT_{t-s}\Psi(s,X(s))\dif W(s)\in\mX_1
$$
admits a continuous modification on $[0,\tau_n)$.

Moreover, starting from (\ref{Eqq4}),
as in the proof of  Theorem \ref{Th1},
there exists an $a>0$ such that for $p$ sufficiently large
$$
\mE\left(\sup_{t\not= t'\in[0,T\wedge\tau_n]}
\frac{\|X(t')-X(t)\|^p_{\mX^\a}}
{|t'-t|^{a p}}\right)\leq C_{n,T,p}.
$$
Thus, by {\bf (M2)} and {\bf (M4)} we know that
$$
s\mapsto \Phi(s,X(s))\in\mX \mbox{ is H\"older continuous on $[0,T\wedge\tau_n]$ $P$-a.s.}.
$$
Therefore, by Lemma \ref{Le6} we have
$$
t\mapsto\int^t_0\fT_{t-s}\Phi(s,X(s))\dif s\in C([0,T\wedge\tau_n],\mX_1), \ \ P-a.s..
$$
Noting that $x_0\in\mX_1$ and
\ce
1_{\{t\leq \tau_n\}}\cdot X(t)&=&1_{\{t\leq \tau_n\}}\cdot \fT_{t} x_0+
1_{\{t\leq \tau_n\}}\cdot \int^{t}_0\fT_{t-s}\Phi(s,X(s))\dif s\\
&&+1_{\{t\leq \tau_n\}}\cdot \int^{t}_0\fT_{t-s}\Psi(s,X(s))\dif W(s),
\ \ \forall t\geq 0, \ \ P-a.s,
\de
by $\tau_n\nearrow\tau$, we therefore have
that $t\mapsto X(t)\in\mX_1$ is continuous on $[0,\tau)$ $P$-a.s..

Lastly, by stochastic Fubini's theorem (cf. \cite[Section 6]{On})
we have
\ce
\int^{t}_0\fL X(s)\dif s&=&
\int^{t}_0\fL\fT_{s} x_0\dif s
+\int^{t}_0\int^{s}_0\fL\fT_{s-r}\Phi(r,X(r))\dif r\dif s\\
&&+\int^{t}_0\int^{s}_0\fL\fT_{s-r}\Psi(r,X(r))\dif W(r)\dif s\\
&=&x_0-\fT_{t}x_0
+\int^{t}_0\int^{t}_r\fL\fT_{s-r}\Phi(r,X(r))\dif s\dif r\\
&&+ \int^{t}_0\int^{t}_r\fL\fT_{s-r}\Psi(r,X(r))\dif s\dif W(r)\\
&=&x_0-\fT_{t}x_0
+\int^{t}_0\Big[\Phi(r,X(r))-\fT_{t-r}\Phi(r,X(r))\Big]\dif r\\
&&+\int^{t}_0\Big[\Psi(r,X(r))-\fT_{t-r}\Psi(r,X(r))\Big]\dif W(r)\\
&=&x_0-X(t)+\int^{t}_0\Phi(s,X(s))\dif s+\int^{t}_0\Psi(s,X(s))\dif W(s)
\de
on $\{t\leq \tau_n\}$. The proof is completed by letting $n\to\infty$.
\end{proof}

\subsection{SPDEs driven by fractional Brownian motions}
In this subsection, we study the existence-uniqueness and large deviation
for SPDEs driven by additive fractional Brownian motions.
Let for  $H\in(0,1)$
$$
K_H(t,s):=\left(c_H(t-s)^{H-\frac{1}{2}}+s^{H-\frac{1}{2}}F(t/s)\right)1_{\{s<t\}}, \quad s,t\in[0,1],
$$
where $c_H:=\left(\frac{2H\Gamma(3/2-H)}{\Gamma(H+1/2)\Gamma(2-2H)}\right)^{1/2}$,
$\Gamma$ denotes the usual
Gamma function, and
$$
F(u):=c_H(\frac{1}{2}-H)\int^u_1(r-1)^{H-\frac{3}{2}}(1-r^{H-\frac{1}{2}})\dif r.
$$
The sequence of independent
fractional Brownian motions with Hurst parameter $H\in(0,1)$ may be
defined by  (cf. \cite{De-Us})
$$
W^k_H(t):=\int^t_0K_H(t,s)\dif W^k(s),\ \ k=1,2,\cdots,
$$
which has the covariance function
$$
R_H(t,s)=\mE(W^k_H(t)W^k_H(s))=\frac{1}{2}(s^{2H}+t^{2H}-|t-s|^{2H}).
$$

Consider the following stochastic partial differential equation driven by $\{W^k_H,k\in\mN\}$
\be
\dif X(t)=[-\fL X(t)+\Phi(t,X(t))]\dif t+\Psi(t)\dif W_H(t),\ \ X(0)=x_0\in\mX,\label{Eq6}
\ee
where $\Psi(t)$ is a deterministic function and
will be specified in Theorem \ref{THf}.

As above, we consider the  mild solution:
\be
X(t)=\fT_t x_0+\int^t_0 \fT_{t-s}\Phi(s,X(s))\dif s+\int^t_0\fT_{t-s}\Psi(s)\dif W_H(s).\label{Eq5}
\ee
Here the stochastic integral is defined by the integration by parts formula as
\ce
\int^t_0\fT_{t-s}\Psi(s)\dif W_H(s)&:=&\Psi(t)W_H(t)+\int^t_0 W_H(s)
[\fL\fT_{t-s}\Psi(s)-\fT_{t-s}\dot\Psi(s)] \dif s\\
&=&\int^t_0B(t,s) \dif W(s),
\de
where
$$
B(t,s):=\Psi(t) K_H(t,s)+\int^t_sK_H(u,s)
\big[\fL\fT_{t-u}\Psi(u)-\fT_{t-s}\dot\Psi(s)\big]\dif u.
$$
We also define
$$
A(t,s,x):=\fT_{t-s}\Phi(s,x).
$$
Then we have:
\bt\label{THf}
Assume that $\Phi$ satisfies {\bf (M1)$'$} and {\bf (M2)$'$} and
$\Psi$ satisfies for some $\gamma>0$
\be
\|\Psi(t')-\Psi(t)\|_{L_2(l^2;\mX)}\leq C|t'-t|^{\gamma}, \ \ t,t'\in[0,1]
\label{Po2}
\ee
and for some $\d\in(0,1)$
\be
\sup_{t\in[0,1]}
\Big(\|\Psi(t)\|_{L_2(l^2;\mX_\d)}
+\|\dot\Psi(t)\|_{L_2(l^2;\mX_{\d-1})}\Big)<+\infty, \ \ t\in[0,1].\label{Po3}
\ee
Then for any $x_0\in\mX$,
there exists a unique maximal solution for \textsc{Eq.}(\ref{Eq5})  in the sense of
Definition \ref{Def1}. In particular, if  $\Phi$ also satisfies {\bf (M3)$'$},
then there is no explosion for \textsc{Eq.}(\ref{Eq5}).
\et
\begin{proof}
As the proof of Theorem \ref{Th5}, one can check that
$A$ satisfies {\bf (H2)$'$}-{\bf (H4)$'$}.
In order to finish the proof by Theorem \ref{Main}, we need to verify that
$B$ also satisfies {\bf (H2)} and {\bf (H4)}.
We first check that for some $\gamma'>0$
\be
\int^t_0\|B(t',s)-B(t,s)\|^2_{L_2(l^2;\mX)}\dif s
\preceq |t-t'|^{\gamma'},\ \ 0\leq t<t'\leq 1.\label{Po4}
\ee

Noting that
\be
K_H(t,s)\leq C|t-s|^{H-\frac{1}{2}}+C s^{-|H-\frac{1}{2}|}\label{Po1}
\ee
and
\ce
\int^t_0[K_H(t,s)-K_H(t',s)]^2\dif s&\leq&R_H(t,t)-2R_H(t,t')+R_H(t',t')\\
&=&t^{2H}-(t^{2H}+(t')^{2H}-|t'-t|^{2H})+(t')^{2H}\\
&=&|t'-t|^{2H},
\de
by (\ref{Po2}) we have
$$
\int^t_0\|\Psi(t') K_H(t',s)-\Psi(t) K_H(t,s)\|^2_{L_2(l^2;\mX)}\dif s
\preceq |t'-t|^{2H\wedge\gamma}.
$$

Observe that
\ce
&&\int^t_0\left\|\int^{t'}_sK_H(u,s)\fL\fT_{t'-u}\Psi(u)\dif u-
\int^{t}_sK_H(u,s)\fL\fT_{t-u}\Psi(u)\dif u\right\|_{L_2(l^2;\mX)}^2\dif s\\
&&\quad\preceq
\int^t_0\left|\int^{t'}_tK_H(u,s)\cdot\|\fL\fT_{t'-u}\Psi(u)
\|_{L_2(l^2;\mX)}\dif u\right|^2\dif s\\
&&\quad\quad+\int^t_0\left|\int^{t}_sK_H(u,s)
\cdot\|\fL(\fT_{t'-u}-\fT_{t-u})\Psi(u)\|_{L_2(l^2;\mX)}\dif u\right|^2\dif s\\
&&\quad=:J_1+J_2.
\de
By (\ref{Po1}), (iii)  (iv) of Proposition \ref{Pr1} and (\ref{Po3}), we have
\ce
J_1&\preceq&\int^t_0\left|\int^{t'}_t\Big[(u-s)^{H-\frac{1}{2}}+s^{-|H-\frac{1}{2}|}\Big]
\cdot(t'-u)^{\d-1}\dif u\right|^2\dif s\\
&\preceq&\int^t_0\Big[(t-s)^{H-\frac{1}{2}}+s^{-|H-\frac{1}{2}|}\Big]^2
\cdot\left|\int^{t'}_t(t'-u)^{\d-1}\dif u\right|^2\dif s\\
&\preceq&|t'-t|^{2\d}
\de
and
\ce
J_2&\preceq&\int^t_0\left|\int^{t}_s\Big[(u-s)^{H-\frac{1}{2}}+s^{-|H-\frac{1}{2}|}\Big]
\cdot(t'-t)^{\frac{\d}{2}}\cdot(t-u)^{\frac{\d}{2}-1}\dif u\right|^2\dif s\\
&\preceq&(t'-t)^{\d}\int^t_0\left|\int^{t}_s(u-s)^{H-\frac{1}{2}}
(t-u)^{\frac{\d}{2}-1}\dif u\right|^2\dif s\\
&&+(t'-t)^{\d}\int^t_0s^{-2|H-\frac{1}{2}|}\left|\int^{t}_s
(t-u)^{\frac{\d}{2}-1}\dif u\right|^2\dif s\\
&=:&J_{21}+J_{22}.
\de
It is clear that
$$
J_{22}\preceq (t'-t)^{\d}\int^t_0s^{-2|H-\frac{1}{2}|}
(t-s)^{\d}\dif s\preceq (t'-t)^{\d}.
$$
For $J_{21}$, let us make the following elementary estimation:
\ce
&&\int^t_0\left|\int^{t}_s(u-s)^{H-\frac{1}{2}}
(t-u)^{\frac{\d}{2}-1}\dif u\right|^2\dif s\\
&&\quad\preceq\int^t_0\left|\int^{\frac{t+s}{2}}_s(u-s)^{H-\frac{1}{2}}
(t-u)^{\frac{\d}{2}-1}\dif u\right|^2\dif s\\
&&\quad\quad+\int^t_0\left|\int^{t}_{\frac{t+s}{2}}(u-s)^{H-\frac{1}{2}}
(t-u)^{\frac{\d}{2}-1}\dif u\right|^2\dif s\\
&&\quad\preceq\int^t_0(t-s)^{\d-2}\left|\int^{\frac{t+s}{2}}_s(u-s)^{H-\frac{1}{2}}
\dif u\right|^2\dif s\\
&&\quad\quad+\int^t_0(t-s)^{2H-1}\left|\int^{t}_{\frac{t+s}{2}}
(t-u)^{\frac{\d}{2}-1}\dif u\right|^2\dif s\\
&&\quad\preceq \int^t_0(t-s)^{2H+\d-1}\dif s\leq C.
\de
Hence
$$
J_1+J_2\leq C(t'-t)^{\d}.
$$
Similarly, we have
\ce
&&\int^t_0\left\|\int^{t'}_sK_H(u,s)\fT_{t'-u}\dot\Psi(u)\dif u-
\int^{t}_sK_H(u,s)\fT_{t-u}\dot\Psi(u)\dif u\right\|_{L_2(l^2;\mX)}^2\dif s\\
&&\qquad\qquad\leq C(t'-t)^{\d}.
\de
Summing up the above calculations, we get (\ref{Po4}). Thus,
$B$ satisfies {\bf (H4)$'$}.

Moreover, from the above calculations, one can see that
$$
\|B(t,s)\|_{L_2(l^2;\mX)}^2\leq C(|t-s|^{2H-1}+s^{-|2H-1|}+|t-s|^{2H+\d-1})=:\kappa(t,s).
$$
So, $B$ satisfies  {\bf (H2)$'$} with $\kappa\in\sK_{>1}$.

Lastly, if  $\Phi$ also satisfies {\bf (M3)$'$}, then the non-explosion
follows from Theorem \ref{Th4}.
\end{proof}

Consider the following small perturbation of \textsc{Eq.}(\ref{Eq5}):
$$
X_\eps(t)=\fT_t x_0+\int^t_0 \fT_{t-s}\Phi(s,X_\eps(s))\dif s+\sqrt{\eps}
\int^t_0\fT_{t-s}\Psi(s)\dif W_H(s).
$$
A direct application of Theorem \ref{Large} yields that
\bt
Keep the same assumptions as Theorem \ref{THf}, where $\Phi$ satisfies
{\bf (M1)$'$}-{\bf (M3)$'$}. Then for any $x_0\in\mX_\a$ ($\a>0$), $\{X_\eps,\eps\in(0,1)\}$
satisfies the large deviation principle in $\mC_1(\mX)$ with the rate function $I(f)$ given by
$$
I(f):=\frac{1}{2}\inf_{\{h\in\ell^2_1:~f=X^h\}}\|h\|^2_{\ell^2_1},\ \ f\in\mC_1(\mX),
$$
where $X^h$ solves the following integral equation
$$
X^h(t)=\fT_tx_0+\int^t_0\fT_{t-s}\Phi(s,X^h(s))\dif s+\int^t_0B(t,s)\dot h(s)\dif s.
$$
\et
\br
Let $\Psi_0\in L_2(l^2;\mX_\d)$ for some $\d\in(0,1)$. Then $\Psi(t):=\fT_t\Psi_0$ satisfies (\ref{Po2})
and (\ref{Po3}) by Proposition \ref{Pr1}.
Moreover, under stronger assumptions on $\Psi(t)$, we can also
prove the existence of strong solutions for \textsc{Eq.}(\ref{Eq5}) as Theorem \ref{Th55}.
\er

\section{Application to SPDEs in bounded domains of $\mR^d$}

Let $\cO$ be an open bounded domain of $\mR^d$ with smooth boundary.
For $m\in\mN$, by $C^m(\cO)$ (resp. $C^m_0(\cO)$)
we denote the set of all $m$-times continuously differentiable functions
in $\cO$ (resp. with compact support in $\cO$).
For $u\in C^m(\cO)$ and $p\geq 1$ we define
$$
\|u\|_{m,p}:=\left(\sum_{j=0}^m\int_\cO |D^j u(x)|^p\dif x\right)^{1/p},
$$
where $D^j$ is the usual derivative operator.
The Sobolev spaces $W^{m,p}(\cO)$ and $W^{m,p}_0(\cO)$ are defined respectively as
the completions of $C^m(\cO)$ and $C^m_0(\cO)$ with respect to the norm $\|\cdot\|_{m,p}$.

Let $\sA(x,D)$ be a strongly elliptic differential operator in $\cO$ of
the form (cf. \cite{Fr,Pa}):
$$
\sA(x,D)u:=\sum_{k=0}^{2m}\sum_{\a_1+\cdots\a_d=k}a_{\a_1\cdots\a_d}(x)D^{\a_1}_1\cdots
D^{\a_d}_d u,\ \ m\geq 1,
$$
where $a_{\a_1\cdots\a_d}(x)\in C^\infty(\bar\cO)$,
and $D^{\a_j}_j$ is the $\a_j$-order derivative with respect to the $j$-th variable.
We consider the following stochastic partial differential equation:
\be\label{Sys}
\left\{
\begin{aligned}
\dif u(t,x)=&\Big[\sA(x,D)u(t)+\varphi(t,x,u,Du,\cdots,D^{2m-1}u)\Big]\dif t\\
&+ \psi(t,x,u,Du,\cdots,D^{m-1}u)\dif W(t),\\
\frac{\p^j u(t,x)}{\p\nu^j}=&0,\ \ j=0,1,\cdots,m-1,\ x\in\p\cO,\\
u(0,x)=&u_0(x),
\end{aligned}
\right.
\ee
where $\frac{\p^j}{\p \nu^j}$ denotes the $j$-th outward normal derivative,
$\varphi$ and $\psi$ are two measurable functions with the entries:
\ce
\varphi:\mR_+\times\cO\times\mR\times\mR^d\times\cdots\times\mR^{(2m-1)d}&\to&\mR,\\
\psi:\mR_+\times\cO\times\mR\times\mR^d\times\cdots\times\mR^{(m-1)d}&\to& l^2.
\de

Define for $p>1$ and $\lambda>0$
$$
\fL_p u:=\lambda u-\sA(x,D)u
$$
with
$$
u\in\sD(\fL_p):=W^{2m,p}(\cO)\cap W^{m,p}_0(\cO).
$$
It is well known that for $u\in\sD(\fL_p)$ (cf. \cite[p. 212, Theorem 3.1]{Pa} or \cite{Tr})
\be
\|u\|_{2m,p}\preceq \|\fL_p u\|_{L^p}+\|u\|_{L^p},\label{Po7}
\ee
and $(\fL_p,\sD(\fL_p))$ is a sectorial operator on $\mX^p_0=L^p(\cO)$
with $0\in\rho(\fL_p)$ for $\lambda$ large enough (cf. \cite[p. 213, Theorem 3.5]{Pa}).
Below we shall write for $p>1$ and $\a\geq 0$
$$
\mX_\a^p:=\sD(\fL^\a_p).
$$
We first recall the following well known result (cf. \cite[p.243]{Pa}).
\bl\label{Le5}
For any $p>1$, $j<2m$ and any $0\leq\a'<\frac{j}{2m}<\a\leq 1$ we have
\be
\|\fL^{\a'}_p u\|_{L^p}\preceq\|u\|_{j,p}\preceq \|\fL_p^\a u\|_{L^p}, \ \ u\in\sD(\fL_p^\a).
\label{PP4}
\ee
Moreover,
$$
\mX_\a^p\hookrightarrow W^{k,q} \ \ \mbox{ for $k-\frac{d}{q}<2m\a-\frac{d}{p}$,\ $q\geq p$}
$$
and
\be
\mX_\a^p\hookrightarrow C^\nu(\bar\cO) \ \ \mbox{ for $0\leq\nu<2m\a-\frac{d}{p}$},\label{PL2}
\ee
where $C^\nu(\bar\cO)$ is the usual H\"older space (cf. \cite{Ad, Pa}).
\el

In this section, we fix
\be
p>d\mbox{ and }\frac{2m-1+\frac{d}{p}}{2m}<\a_0<\a<1\label{Con0}
\ee
so that
\be
m(1-\a)^2<(\a-\a_0).\label{Con1}
\ee
Suppose that
\begin{enumerate}[{\bf(F1)}]
\item For any $T,R>0$, there exist $\d>0$ and $C_{R,T}>0$
such that for all $s,t\in[0,T]$, $x\in\cO$ and $U,V\in \mR^{m(2m-1)d+1}$ with $|U|,|V|\leq R$
$$
|\varphi(t,x,U)-\varphi(s,x,V)|\leq C_{R,T}(|t-s|^\d+|U-V|).
$$
Moreover, $\sup_{x\in\cO}|\varphi(0,x,0)|<+\infty$.
\item For each $t\in\mR_+$, $ \psi(t,*)\in C^{m+1}(\bar\cO\times\mR^{m(m+1)d/2+d+1};l^2)$.
Here and below, the asterisk stands  for the rest variables.
\item For each $u\in\mX_{\a_0}^p$,
$$
 \psi(t,\cdot,u,Du,\cdots,D^{m-1}u)\in \mX_{\frac{\a}{2}}^p.
$$
\item
For any $T>0$, there exist  constant $C_T>0$ and $\l_0\in L^p(\cO)$
such that for all $t\in[0,T]$, $x\in\cO$
and $U\in \mR^{m(2m-1)d+1}$
\be
|\varphi(t,x,U)|\leq C_T(\l_0(x)+|U|),\label{Es7}
\ee
and
\be
\label{Es8}
\left\{
\begin{aligned}
& \psi(t,x,u,Du,\cdots, D^{m-1}u)=\sum_{j=0}^{m-1} g_{j}(t)
\cdot D^{j}u+ \psi_0(t,x),\ \ &m\geq 2,\\
&\left.
\begin{aligned}
&\mbox{for some $\d>0$ and each $r\in\mR$},\ \mbox{supp}( \psi(t,\cdot,r))\subset\cO_\d,\\
&\psi(t,*)\in C^2(\bar\cO\times\mR^{d+1};l^2),\ \|\p_r \psi(t,x,r)\|_{l^2}\leq C_T,\\
&\|D_x \psi(t,x,r)\|_{l^2\times\mR^d}+\| \psi(t,x,r)\|_{l^2}\leq C(f_0(t,x)+|r|)
\end{aligned}
\right\},&m=1,
\end{aligned}
\right.
\ee
where $\cO_\d\subset\bar\cO_\d\subset\cO$ is an open subset, and for each $j=0,\cdots,m-1$,
\ce
&&t\mapsto g_{j}(t)\in l^2\times\mR^{jd},\\
&&t\mapsto \psi_0(t,\cdot)\in l^2\times\mX^p_{\a},\\
&&t\mapsto f_0(t,\cdot)\in L^p(\cO)
\de
are bounded measurable functions.
\end{enumerate}

We remark that {\bf (F3)} is related to the boundary conditions, e.g., $\Psi\equiv$constant
 does not satisfy {\bf (F3)}.
 It is easy to see that (\ref{Es8}) implies {\bf (F3)}.

Set
\be
\Phi(t,u)(x)&:=&\varphi(t,x,u,Du,\cdots,D^{2m-1}u)+\lambda u,\label{Phi}\\
\Psi(t,u)(x)&:=& \psi(t,x,u,Du,\cdots, D^{m-1}u).\label{Psi}
\ee
Then the system (\ref{Sys}) can be written as the following abstract form:
\be
\dif u(t)=[-\fL_p u+\Phi(t,u(t))]\dif t+\Psi(t,u(t))\dif W(t),\ \
u(0)=u_0.\label{Eq8}
\ee

Using Theorem \ref{Th55}, we have the following result.
\bt\label{62}
Let $p>d$ and $\a,\a_0$ satisfy (\ref{Con0}) and (\ref{Con1}).
Assume that {\bf (F1)}-{\bf (F3)} hold.
For any $u_0\in\mX_1^p$,
there exists a unique maximal strong solution $(u,\tau)$ for \textsc{Eq.} (\ref{Eq8}) so that
\begin{enumerate}[(i)]
\item $t\mapsto u(t)\in\mX^p_1$ is continuous on $[0,\tau)$  almost surely;
\item $\lim_{t\uparrow\tau}\|u(t)\|_{\mX^p_\a}=+\infty$ on $\{\omega:\tau(\omega)<+\infty\}$;
\item it holds that in $L^p(\cO)$
\ce
u(t)&=&u_0-\int^{t}_0\fL_p u(s)\dif s+\int^{t}_0\Phi(s,u(s))\dif s
+\int^{t}_0\Psi(s,u(s))\dif W(s)\\
&=&u_0+\int^{t}_0\sA(x,D)u(s)\dif s+\int^{t}_0\varphi(s,x,u(s),Du(s),\cdots,D^{2m-1}u(s))\dif s\\
&&+\int^{t}_0 \psi(s,x,u(s),Du(s),\cdots,D^{m-1}u(s))\dif W(s)
\de
for all $t<\tau$, $P$-a.s..
\end{enumerate}
Moreover, if {\bf (F4)} holds, then
$$
\tau=+\infty,\ \  a.s..
$$
\et
\begin{proof}
We only need to verify that {\bf (M1)}-{\bf (M4)} hold
for $\Phi$ and $\Psi$ defined by (\ref{Phi}) and (\ref{Psi}).
In virtue of (\ref{Con0}), by (\ref{PL2}) we have
\be
\|D^j u\|_{C(\bar\cO)}\preceq \|u\|_{\mX_{\a_0}^p},\ \ j=0,1,\cdots, 2m-1.\label{PO7}
\ee
It is easy to see by {\bf (F1)} that $\Phi$ given by (\ref{Phi})
is locally Lipschitz continuous and locally bounded
with respect to $u$ on $\mX_{\a}^p$, and is $\d$-order H\"older continuous with respect to $t$.

Note that by the chain rule, for $j=1,\cdots, m+1$
\be
D^j\Psi(t,u)&=&(\p_{D^{m-1}u} \psi)(t,x,u,Du,\cdots, D^{m-1}u)\cdot D^{m-1+j}u\no\\
&&+\psi_j(t,x,u,Du,\cdots, D^{m-2+k}u),\label{PL1}
\ee
where $\psi_j$ is an $l^2$-valued  continuously differentiable
function of all its variables with the exception
of the $t$-variable.
For any $u,v\in\mX_{\a_0}^p$ with $\|u\|_{\mX_{\a_0}^p}, \|v\|_{\mX_{\a_0}^p}\leq R$,
by {\bf(F2)} and {\bf(F3)} we have
\ce
\|\fL^{\frac{\a}{2}}(\Psi(t,u)-\Psi(t,v))\|^2_{L_2(l^2; \mX^p_0)}
&\stackrel{(\ref{Exa})}{\preceq}&\|\fL^{\frac{\a}{2}}(\Psi(t,u)-\Psi(t,v))\|^2_{L^p(\cO)}\\
&\stackrel{(\ref{PP4})}{\preceq}&\sum_{k=1}^\infty
\sum_{j=0}^m\|D^j(\Psi_k(t,u)-\Psi_k(t,v))\|^2_{L^p(\cO;l^2)}\\
&\preceq&\sum_{k=1}^\infty\sum_{j=0}^m\|D^j(\Psi_k(t,u)-\Psi_k(t,v))\|^2_{C(\bar\cO)}\\
&\stackrel{(\ref{PL1},\ref{PO7})}{\leq}&
C_R\sum_{j=0}^{2m-1}\|D^j(u-v)\|_{C(\bar \cO)}^2\\
&\stackrel{(\ref{PO7})}{\preceq}& C_R\|u-v\|^2_{\mX_{\a_0}^p}\leq C_R\|u-v\|_{\mX_{\a}^p}^2.
\de
Thus, {\bf(M2)} holds.

We next look at {\bf(M4)}. As above, by (\ref{PL1}) and (\ref{PO7}) we have
\be
\|D^{m+1}\Psi(t,u)\|_{L^p(\cO;l^2)}\leq C_R(1+\|D^{2m}u\|_{L^p})
\stackrel{(\ref{Po7})}{\leq} C_R(1+\|u\|_{\mX^p_1})\label{PO8}
\ee
for all $u\in\mX_1^p$ with $\|u\|_{\mX_{\a_0}^p}\leq R$.
By (\ref{Con1}), we may choose
$$
1<\a'<\a''<\frac{m+1}{m}
$$
such that
$$
\theta:=\frac{\a-\a_0}{1-\a_0}>\frac{\a'-\a}{\a''-\a}=:\theta'.
$$
Thus, for all $u\in\mX_1^p$ with $\|u\|_{\mX_{\a_0}^p}\leq R$, we have
\ce
\|\Psi(t,u)\|^2_{L_2(l^2;\mX_{\a''/2}^p)}&\stackrel{(\ref{Exa})}{\preceq}&
\|\fL^{\frac{\a''}{2}}\Psi(t,u)\|^2_{L^p(\cO;l^2)}\\
&\stackrel{(\ref{PP4})}{\preceq}&\sum_{j=0}^{m+1}\|D^j\Psi(t,u)\|^2_{L^p(\cO;l^2)}\\
&\preceq&\sum_{j=0}^{m}\sum_{k=1}^\infty\|D^j\Psi_k(t,u)\|^2_{C(\bar\cO)}\\
&&+\|D^{m+1}\Psi(t,u)\|_{L^p(\cO;l^2)}^2\\
&\stackrel{(\ref{PO8})}{\leq}& C_R(1+\|u\|^2_{\mX_1^p}).
\de
Using Lemma \ref{Le7} with the data $\a_0$, $\theta$ and $\theta'$ as above and
$$
\a_1=1,\ \ \a_2=\frac{\a}{2},\ \ \a_3=\frac{\a''}{2},
$$
we obtain that for all $u\in\mX^p_\a$ with $\|u\|_{\mX_{\a}^p}\leq R$
\ce
\|\Psi(t,u)\|_{L_2(l^2;\mX^p_{\a'/2})}
=\|\Psi(t,u)\|_{L_2(l^2;\mX^p_{\theta'(\a_3-\a_2)+\a_2})}\leq C_R.
\de
Thus, {\bf(M4)} holds.

We now verify {\bf(M3)} under {\bf (F4)}.
First of all, by the linear growth of $\varphi(t,x,*)$ with respect to $*$, we have
\ce
\|\Phi(t,u)\|_{\mX^p_0}
&=&\|\varphi(t,\cdot,u,Du,D^2u,\cdots,D^{2m-1}u)\|_{L^p}\\
&\preceq& 1+\sum_{j=0}^{2m-1}\|D^ju\|_{L^p}\\
&\preceq& 1+\|u\|_{2m-1,p}\\
&\stackrel{(\ref{PP4})}{\preceq}& 1+\|u\|_{\mX^p_\a}.
\de
For $\Psi$, we only consider the case of $m=1$, and have
\ce
\|\Psi(t,u)\|^2_{L_2(l^2;\mX^p_{\frac{\a}{2}})}
&\stackrel{(\ref{Exa})}{\preceq}&\|\fL_p^{\frac{\a}{2}}\Psi(t,u)\|^2_{L^p(\cO;l^2)}\\
&\stackrel{(\ref{PP4})}{\preceq}&\|\Psi(t,u)\|^2_{L^p(\cO;l^2)}+\|D\Psi(t,u)\|^2_{L^p(\cO;l^2)}.
\de
Noting that
$$
D\Psi(t,u)=(D_x \psi)(t,x,u)+(\p_u \psi)(t,x,u)Du
$$
by (\ref{Es8}) we have
$$
\|D\Psi(t,u)\|^2_{L^p(\cO;l^2)}\preceq 1+\|u\|_{L^p}+\|Du\|_{L^p}
\preceq 1+\|u\|_{\mX^p_\a}.
$$
So
$$
\|\Psi(t,u)\|^2_{L_2(l^2;\mX^p_{\frac{\a}{2}})}\preceq 1+\|u\|_{\mX^p_\a}.
$$
Thus, {\bf(M3)} holds.
\end{proof}

Consider the small perturbation of equation (\ref{Eq8}):
\be
\dif u_\eps(t)=[-\fL_p u_\eps(s)+\Phi(t,u_\eps(t))]\dif t+\sqrt{\eps}\Psi(t,u_\eps(t))\dif W(t),\ \
u_\eps(0)=u_0.\label{Eqq8}
\ee
Using Theorem \ref{Th6}, we have
\bt
Let $p>d$ and $\a,\a_0$ satisfy (\ref{Con0}) and (\ref{Con1}).
Assume that {\bf (F1)}  and {\bf (F4)}  hold. Let $u_0\in\mX_1^p$.
Then $\{u_\eps,\eps\in(0,1)\}$ satisfies the large deviation principle in $\mC_T(\mX_\a^p)$
with the rate function $I(f)$ given by (\ref{Rate1}).
\et

\section{Application to SPDEs on complete Riemannian manifolds}

Let $(M,\fg)$ be a $d$-dimensional complete Riemannian manifold without boundary.
The Riemannian volume is denoted by $\dif_\fg x$. Let $\nabla$ denote the gradient or
covariant derivative associated with $\fg$,
$\Delta$ the Laplace Beltrami operator, $T(M)$ the tangent bundle.
Let $L^p(M,\dif_\fg x)$ be the usual real $L^p$-space on $M$ with respect to $\dif_\fg x$.
It is well known that the symmetric heat semigroup $(\fT_t)_{t\geq 0}$ associated with $\Delta$
is strongly continuous and contracted
on $L^p(M,\dif_\fg x)$ for $1\leq p<+\infty$, which is also contracted on
$L^\infty(M,\dif_\fg x)$ (cf. Strichartz \cite[Theorem 3.5]{Str}).
Therefore, for each $1<p<+\infty$, $(\fT_t)_{t\geq 0}$ forms an analytic semigroup on $L^p(M,\dif_\fg x)$
 (cf. Stein \cite[p.67 Theorem 1]{St}).
The Bessel spaces over $M$ are defined by
$$
\mH^p_\a:=(I-\Delta)^{-\a/2}(L^p(M,\dif_\fg x)).
$$

In this section, we make the following geometric assumptions:

{\bf (G)$_n$}: The Ricci curvature $\Ric_\fg$ and
curvature $R_\fg$ tensors together with their covariant derivatives
up to $n$-th order  are bounded.

{\bf (G)$_{inj}$}: The injectivity radius of $(M,\fg)$ is strictly positive.

It was proved by Yoshida  \cite{Yo} that under {\bf (G)$_n$},
an equivalent norm of $\mH^p_n$ is given by the covariant derivatives up to $n$-th order, i.e., there
are two positive constants $C_1$ and $C_2$ such that for any $u\in C_0^\infty(M)$
\be
C_1\sum_{k=0}^n\|\nabla^ku\|_{L^p}\leq \|(I-\Delta)^{n/2}u\|_{L^p}
\leq C_2\sum_{k=0}^n\|\nabla^ku\|_{L^p},\label{re}
\ee
where $\nabla^k$ denotes the $k$-th covariant derivative. As an example, the components
of $\nabla u$ in local coordinates are given by $(\nabla u)_i=\p_i u$, while
the components of $\nabla^2 u$ in local coordinates are given by
$$
(\nabla^2 u)_{ij}=\p_{ij}u-\Gamma^k_{ij}\p_k u,
$$
where $\Gamma^k_{ij}$ are Christoffel symbols. By definition one has that
$$
|\nabla^ku|^2=g^{i_1j_1}\cdots g^{i_kj_k}(\nabla^k)_{i_1\cdots i_k}(\nabla^k)_{j_1\cdots j_k},
$$
where $g_{ij}=\fg(\p_i,\p_j)$ and $(g^{ij})$ denotes the inverse matrix of $(g_{ij})$.

We remark that when $n=1$, (\ref{re}) was first proved by Bakry \cite{Ba} under the assumption
that Ricci curvature is bounded from below.

The following embedding result was proved in \cite{Zh4}. We refer to \cite{Au, He0, He1}
for a detailed study of integer order Sobolev spaces over $M$.
\bt\label{Emd}
Under {\bf (G)$_{n+1}$} and {\bf (G)$_{inj}$}, for $\a\in(0,1)$ and $p>d/\a$ we have
$$
\mH^p_{n+\a}\hookrightarrow C^{n}_b(M),
$$
where $C^{n}_b(M)$ denotes the Banach space of all $n$-times
continuously differentiable functions on $M$ with
$$
\|u\|_{C^n_b}:=\sup_{x\in M}\sum_{k=0}^n|\nabla^k u(x)|<+\infty.
$$
\et
Consider the following SPDE:
\be
\label{Eqm}\left\{
\begin{aligned}
\dif u(t,x)=&\big[\Delta u(t,x)+\varphi(t,x,u(t),\fg(Y(x),\nabla u(t)))\big]\dif t\\
&+\psi(t,x,u(t,x))\dif W(t),\\
u(0,x)=&u_0(x),
\end{aligned}
\right.
\ee
where $Y:M\to T(M)$ is a measurable vector field with
\be
\sup_{x\in M}\fg(Y(x),Y(x))<+\infty\label{YY}
\ee
and
\ce
\varphi:\mR_+\times\Omega\times M\times\mR^2\to \mR&\in&\cM\times\cB(M)\times\cB(\mR^2)/\cB(\mR),\\
\psi:\mR_+\times\Omega\times M\times \mR\to l^2&\in&\cM\times\cB(M)\times\cB(\mR)/\cB(l^2).
\de

In this section, we fix
\be
p>d\mbox{ and }\frac{3}{2}-\sqrt{\frac{3}{4}-\frac{d}{2p}}<\a<1.\label{Con00}
\ee
Assume that
\begin{enumerate}[{\bf (R1)}]
\item For each $T,R>0$, there exist constants $C_{R,T},\d>0$
and $\l^{\varphi,1}_{R,T},\l^{\varphi,2}_{R,T}\in L^p(M,\dif_\fg x)$ such that
for all $s,t\in[0,T], \omega\in\Omega$, $x\in M$
and $|u_1|,|v_1|,|u_2|,|v_2|\leq R$
\ce
&&|\varphi(t,\omega,x,u_1,u_2)-\varphi(t,\omega,x,v_1,v_2)|\\
&&\qquad\leq C_{R,T}\big(\l^{\varphi,1}_{R,T}(x)\cdot|t-s|^\d+|u_1-v_1|+|u_2-v_2|\big)
\de
and
$$
|\varphi(t,\omega,x,u_1,u_2)|\leq \l^{\varphi,2}_{R,T}(x).
$$
\item For each $(t,\omega)\in\mR_+\times\Omega$,
$\psi(t,\omega,\cdot,\cdot)\in C^2(M\times\mR;l^2)$.
For each $T,R>0$, there exist constant $C_{R,T}>0$
and $\l^\psi_{R,T}\in L^p(M,\dif_\fg x)$ such that for all $t\in[0,T], \omega\in\Omega$, $x\in M$
and $|u|\leq R$
\ce
\|\nabla_x\p_u\psi(t,\omega,\cdot,u)\|_{l^2}+\|\p^j_u\psi(t,\omega,\cdot,u)\|_{l^2}\leq C_{R,T}
,\ \ j=1,2
\de
and
\ce
\|\psi(t,\omega,\cdot,u)\|_{l^2}+\|\nabla^j_x\psi(t,\omega,\cdot,u)\|_{l^2}\leq\l^\psi_{R,T}(x),\ \ j=1,2.
\de
\end{enumerate}

\bt\label{Th7}
Let $p>d$ and $\a$ satisfy (\ref{Con00}).
Under {\bf (G)$_{2}$}-{\bf (G)$_{inj}$} and {\bf (R1)}-{\bf (R2)},
for each $u_0\in\mH^p_2$, there exists a unique maximal
strong solution $(u,\tau)$ for \textsc{Eq.}(\ref{Eqm}) so that
\begin{enumerate}[(i)]
\item $t\mapsto u(t)\in\mH^p_2$ is continuous on $[0,\tau)$ almost surely;
\item $\lim_{t\uparrow\tau}\|u(t)\|_{\mH^p_{2\a}}=+\infty$ on $\{\omega:\tau(\omega)<+\infty\}$;
\item it holds that, $P$-a.s., on $[0,\tau)$
\ce
u(t)&=&u_0+\int^{t}_0\big[\Delta u(s)+\varphi(s,\cdot,u(s),\fg(Y(\cdot),\nabla u(s)))\big]\dif s\\
&&+\int^{t}_0\psi(s,\cdot,u(s))\dif W(s)\ \ \mbox{in $L^p(M,\dif_\fg x)$.}
\de
\end{enumerate}
\et
\begin{proof}
Choose $\a_0$ such that
\be
\frac{1}{2}+\frac{d}{2p}<\a_0<3\a-\a^2-1<\a<1.\label{Lp1}
\ee
Let $u,v\in\mH^p_{2\a_0}$ with $\|u\|_{\mH^p_{2\a_0}},\|v\|_{\mH^p_{2\a_0}}\leq R$.
By Theorem  \ref{Emd} we have
\be
\|u\|_{C^1_b}+\|v\|_{C^1_b}\leq C_R.
\label{V5}
\ee
Set
\ce
\Phi(t,\omega,u)&:=&\varphi(t,\omega,\cdot,u,\fg(Y(\cdot),\nabla u)),\\
\Psi(t,\omega,u)&:=&\psi(t,\omega,\cdot,u).
\de
By {\bf (R1)} and (\ref{re}) (\ref{YY}) (\ref{V5}), we have
\ce
\|\Phi(t,\omega,u)-\Phi(s,\omega,v)\|_{L^p}&\preceq& |t-s|^\d+\|u-v\|_{L^p}+
\|\nabla (u-v)\|_{L^p}\\
&\preceq&|t-s|^\d+ \|u-v\|_{\mH^p_1}\\
&\preceq&|t-s|^\d+\|u-v\|_{\mH^p_{2\a}}
\de
and
$$
\|\Phi(t,\omega,u)\|_{L^p}\leq C_R.
$$

Note that
\be
\nabla_x\Psi(t,\omega,u)=(\nabla_x \psi)(t,\omega,\cdot,u)+(\p_u \psi)(t,\omega,x,u)\nabla_x u
\label{V6}
\ee
and
\be
\nabla^2_x\Psi(t,\omega,u)&=&(\nabla^2_x \psi)(t,\omega,\cdot,u)
+2(\nabla_x \p_u\psi)(t,\omega,\cdot,u)\otimes\nabla_x u\no\\
&&+(\p_u \psi)(t,\omega,x,u)\nabla_x u\otimes\nabla_x u+(\p_u \psi)(t,\omega,x,u)\nabla^2_x u.\label{V7}
\ee
By {\bf (R2)} and (\ref{V5}) we have
$$
\|\Psi(t,\omega,u)-\Psi(t,\omega,v)\|_{L^p}\preceq \|u-v\|_{L^p}
$$
and by (\ref{re}) and (\ref{V6})
$$
\|\nabla(\Psi(t,\omega,u)-\Psi(t,\omega,v))\|_{L^p}\preceq \|u-v\|_{\mH^p_1}.
$$
Hence,
$$
\|\Psi(t,\omega,u)-\Psi(t,\omega,v)\|_{\mH^p_\a}\preceq
\|\Psi(t,\omega,u)-\Psi(t,\omega,v)\|_{\mH^p_1}\preceq\|u-v\|_{\mH^p_{2\a_0}}.
$$
Moreover, by  {\bf (R2)} and (\ref{re}) (\ref{V6}) (\ref{V7}),
$$
\|\Psi(t,\omega,u)\|_{\mH^p_{1}}\leq C_{R,T}
$$
and
$$
\|\Psi(t,\omega,u)\|_{\mH^p_{2}}\leq
C_{R,T}(1+\|u\|_{\mH^p_2}).
$$
Using Lemma \ref{Le7} with the data $\a_0$ as above and
$$
\a_3=\a_1=1,\ \ \a_2=\frac{\a}{2},\ \ \frac{\a-\a_0}{1-\a_0}=:\theta
\stackrel{(\ref{Lp1})}{>}\theta'>\frac{1-\a}{2-\a},\ \
$$
we find that for all $u\in\mX^p_\a$ with $\|u\|_{\mX_{\a}^p}\leq R$
\ce
\|\Psi(t,\omega,u)\|_{L_2(l^2;\mH^p_{\a'})}
=\|\Psi(t,\omega,u)\|_{L_2(l^2;\mH^p_{2\theta'(\a_3-\a_2)+2\a_2})}\leq C_R,
\de
where $\a'=2\theta'(\a_3-\a_2)+2\a_2>1$.
Thus, {\bf(M2)} and {\bf(M4)} hold, and the theorem follows from Theorem \ref{Th55}.
\end{proof}

For the non-explosion, we assume that
\begin{enumerate}[{\bf(R3)}]
\item For each $T>0$, there exist $\l_i\in L^p(M,\dif_\fg x), i=0,1,2$
and $k\in\mN$ such that for all
$(t,\omega)\in[0,T]\times\Omega, u,v\in\mR$ such that
\be
u\cdot \varphi(t,\omega,x,u,v)&\leq& C_T|u|\cdot(|u|+|v|+\l_0(x)),\label{V8}\\
|\varphi(t,\omega,x,u,v)|&\leq& C_T(|u|^k+|v|+\l_1(x))\label{V9}
\ee
and
\be
\|\p_u\psi(t,\omega,\cdot,u)\|_{l^2}&\leq&C_T,\\
\|\psi(t,\omega,\cdot,u)\|_{l^2}+\|\nabla_x\psi(t,\omega,\cdot,u)\|_{l^2}
&\leq& C_T(|u|+\l_2(x)).\label{V10}
\ee
\end{enumerate}

The following theorem will follow from the proof of Lemma \ref{Le11} below.
\bt
Keep the same assumptions as in Theorem \ref{Th7}, and also assume
{\bf (R3)}. Let $(u,\tau)$ be the unique maximal strong
solution of \textsc{Eq.}(\ref{Eqm}) in Theorem \ref{Th7}.
Then $\tau=+\infty$ a.s..
\et

Let $\varphi$ and $\psi$ be independent of $\omega$.
Consider now the small perturbation of \textsc{Eq.}(\ref{Eqm}):
\ce\left\{
\begin{aligned}
\dif u_\eps(t,x)=&\big[\Delta u_\eps(t,x)+\varphi(t,x,u_\eps(t),\fg(Y(x),\nabla u_\eps(t)))\big]\dif t\\
&+\sqrt{\eps}\psi(t,x,u_\eps(t,x))\dif W(t),\\
u_\eps(0,x)=&u_0(x)\in\mH^p_2,
\end{aligned}
\right.
\de
as well as the control equation
\be
\label{MP1}\left\{
\begin{aligned}
\dif u^\eps(t,x)=&\big[\Delta u^\eps(t,x)+\varphi(t,x,u^\eps(t),\fg(Y(x),\nabla u^\eps(t)))\big]\dif t\\
&+\psi(t,x,u^\eps(t,x))\dot h^\eps(s)\dif s+\sqrt{\eps}\psi(t,x,u^\eps(t,x))\dif W(t),\\
u^\eps(0,x)=&u_0(x)\in\mH^p_2,
\end{aligned}
\right.
\ee
where $h^\eps\in\cA^{T}_N$ (see (\ref{Op2}) for the definition of $\cA^{T}_N$),
and $T>0$ is fixed below.

Let $(u^\eps,\tau^\eps)$ be the unique maximal strong solution of \textsc{Eq.}(\ref{MP1}).
Define
$$
\tau^\eps_n:=\inf\big\{t: \|u^\eps(t)\|_{\mH^p_{2\a}}>n\big\}.
$$
Then we have:
\bl\label{Le11}
Assume {\bf (R3)}. Then
\ce
\lim_{n\rightarrow\infty}\sup_{\eps\in(0,1)}
P\big\{\omega:\tau^\eps_n(\omega)<T\big\}=0.
\de
\el
\begin{proof}
For the simplicity of notations, we  drop the superscript $\eps$ in $u^\eps$  in the following.
First of all, note that (cf. \cite{Au})
$$
u_0\in\mH^p_2\subset \cap_{q>1}L^q(M,\dif_\fg x).
$$
For $q,r\geq 2$, by the usual It\^o formula (cf. \cite[Theorem A.2]{bp}), we have
$$
\|u(t)\|_{L^q}^{rq}=\|u_0\|_{L^q}^{rq}+J_1(t)+J_2(t)+J_3(t)+J_4(t)+J_5(t)+J_6(t)
$$
on $[0,\tau^\eps_n]$, where
\ce
J_1(t)&:=&rq\int^t_0\|u(s)\|_{L^q}^{(r-1)q}\<|u(s)|^{q-2}u(s),\Delta u(s)\>_{L^2}\dif s,\\
J_2(t)&:=&rq\int^t_0\|u(s)\|_{L^q}^{(r-1)q}\<|u(s)|^{q-2}u(s),\varphi(s,\cdot,u(s),
\fg(Y(\cdot),\nabla u(s)))\>_{L^2}\dif s,\\
J_3(t)&:=&rq\sum_k\int^t_0\|u(s)\|_{L^q}^{(r-1)q}\<|u(s)|^{q-2}u(s),
\psi_k(s,\cdot,u(s))\>_{L^2}\dif W^k_s,\\
J_4(t)&:=&\frac{rq(q-1)}{2}\sum_k\int^t_0\|u(s)\|_{L^q}^{(r-1)q}\<|u(s)|^{q-2},
|\psi_k(s,\cdot,u(s))|^2\>_{L^2}\dif s,\\
J_5(t)&:=&\frac{q^2 r(r-1)}{2}\sum_k\int^t_0\|u(s)\|_{L^q}^{(r-2)q}\<|u(s)|^{q-2}u(s),
\psi_k(s,\cdot,u(s))\>_{L^2}^2\dif s,\\
J_6(t)&:=&rq\int^t_0\|u(s)\|_{L^q}^{(r-1)q}\<|u(s)|^{q-2}u(s),
\psi(s,\cdot,u(s))\dot h^\eps(s)\>_{L^2}\dif s.
\de
For $J_1(t)$ we have
\ce
J_1(t)=-rq(q-1)\int^t_0\|u(s)\|_{L^q}^{(r-1)q}
\int_M|u(s)|^{q-2}|\nabla u(s)|^2\dif_\fg x\dif s.
\de
For $J_2(t)$, by (\ref{V8}) (\ref{YY}) and Young's inequality we have
\ce
J_2(t)&\leq& rq\int^t_0\|u(s)\|_{L^q}^{(r-1)q}
\int_M|u(s)|^{q-1}\big(|u(s)|+|\fg(Y,\nabla u(s))|+\l_0\big)\dif_\fg x\dif s\\
&\leq&-\frac{J_1(t)}{2}+C\int^t_0(\|u(s)\|^{rq}_{L^q}+1)\dif s.
\de
Similarly, by (\ref{V10}) we have
$$
J_4(t)+J_5(t)\leq C\int^t_0(\|u(s)\|^{rq}_{L^q}+1)\dif s,
$$
and by Young's inequality
\be
J_6(t)&\preceq& \int^t_0\|u(s)\|_{L^q}^{(r-1)q}(\|u(s)\|_{L^q}^{q}+1)\cdot \|\dot h^\eps(s)\|_{l^2}\dif s\no\\
&\preceq& N\left(\int^t_0\|u(s)\|_{L^q}^{2(r-1)q}(\|u(s)\|_{L^q}^{q}+1)^2\dif s\right)^{1/2}\no\\
&\preceq& N\sup_{s\in[0,t]}\|u(s)\|_{L^q}^{rq/2}
\cdot\left(\int^t_0(\|u(s)\|_{L^q}^{rq}+1)\dif s\right)^{1/2}\no\\
&\leq& \frac{1}{2}\sup_{s\in[0,t]}\|u(s)\|_{L^q}^{rq}
+C_N\cdot\int^t_0(\|u(s)\|_{L^q}^{rq}+1)\dif s.\label{MP2}
\ee
Combining the above calculations, we obtain
\ce
\sup_{s\in[0,t]}\|u(s\wedge\tau^\eps_n)\|^{rq}_{L^q}\leq 2\|u_0\|_{L^q}^{rq}+
2\sup_{s\in[0,t]}J_3(s\wedge\tau^\eps_n)+C_N\int^{t\wedge\tau^\eps_n}_0(\|u(s)\|^{rq}_{L^q}+1)\dif s.
\de

Set
$$
f_1(t):=\mE\left(\sup_{s\in[0,t]}\|u(s\wedge\tau^\eps_n)\|^{rq}_{L^q}\right).
$$
By BDG's inequality and as (\ref{MP2}) we have
\ce
\mE\left(\sup_{s\in[0,t]}|J_3(s\wedge\tau^\eps_n)|\right)&\preceq&
\mE\left(\int^{t\wedge\tau^\eps_n}_0\|u(s)\|_{L^q}^{2(r-1)q}(\|u(s)\|^{q}_{L^q}+1)^2
\dif s\right)^{1/2}\\
&\leq&\frac{1}{2}f_1(t)+C\mE\left(\int^{t\wedge\tau^\eps_n}_0(\|u(s)\|^{rq}_{L^q}+1)\dif s\right).
\de
Therefore,
\ce
f_1(t)&\leq& 4\|u_0\|_{L^q}^{rq}+C_N\mE\int^{t\wedge\tau^\eps_n}_0(\|u(s)\|^{rq}_{L^q}+1)\dif s\\
&\leq& 4\|u_0\|_{L^q}^{rq}+C_N\int^t_0(f_1(s)+1)\dif s,
\de
which yields by Gronwall's inequality that
\be
\mE\left(\sup_{t\in[0,T]}\|u(t\wedge\tau^\eps_n)\|_{L^q}^{rq}\right)\leq C_{T,N}.\label{V11}
\ee
Here and below, the constant $C_{T,N}$ is independent of $n$ and $\eps$.

Set
$$
\xi^\eps_n(t):=t\wedge\tau^\eps_n
$$
and for $q\geq 2$
$$
f_2(t):=\mE\left(\sup_{t'\leq\xi^\eps_n(t)}\|u(t')\|^q_{\mH^p_{2\a}}\right).
$$
Note that
\ce
u(t)&=&\fT_tu_0+\int^t_0\fT_{t-s}\varphi(s,\cdot,u(s),\fg(Y(\cdot),\nabla u(s)))\dif s\\
&&+\int^t_0\fT_{t-s}\psi(s,\cdot,u(s))\dot h^\eps(s)\dif s+\sqrt{\eps}
\int^t_0\fT_{t-s}\psi(s,\cdot,u(s))\dif W(s)\\
&=:&\fT_tu_0+\cJ_1(t)+\cJ_2(t)+\cJ_3(t).
\de
By (iii) of Proposition \ref{Pr1} and H\"older's inequality
we have, for $q>\frac{1}{1-\a}$
\ce
&&\mE\left(\sup_{t'\in[0,\xi^\eps_n(t)]}\|\cJ_1(t)\|^q_{\mH^p_{2\a}}\right)\\
&&\quad\preceq\mE\left(\sup_{t'\in[0,\xi^\eps_n(t)]}
\int^{t'}_0\frac{1}{(t'-s)^\a}\|\varphi(s,\cdot,u(s),\fg(Y(\cdot),\nabla u(s)))\|_{L^p}\dif s\right)^q\\
&&\quad\preceq\mE\left(\int^{\xi^\eps_n(t)}_0\|\varphi(s,\cdot,u(s),
\fg(Y(\cdot),\nabla u(s)))\|^q_{L^p}\dif s\right)\\
&&\quad\stackrel{(\ref{V9}, \ref{YY})}{\preceq}
\mE\left(\int^{\xi^\eps_n(t)}_0(1+\|u(s)\|^{kq}_{L^{pk}}+
\|\nabla u(s)\|^q_{L^p})\dif s\right)\\
&&\quad\stackrel{(\ref{re})}{\preceq} \int^t_0(\|u(s)\|^q_{\mH^p_1}+1)\dif s
\stackrel{(\ref{V11})}{\preceq} \int^t_0(f_2(s)+1)\dif s.
\de

On the other hand, by (\ref{re}), (\ref{Exa}) and {\bf (R3)} we have, for $u\in\mH^p_1$
\be
\|\psi(s,\cdot,u)\|^2_{L_2(l^2;\mH^p_{\a})}&\preceq&
\|\psi(s,\cdot,u)\|^2_{L^p(M;l^2)}+
\|\nabla \psi(s,\cdot,u)\|^2_{L^p(M;l^2)}\no\\
&\preceq&\|u\|^2_{L^p}+ \|\nabla u\|^2_{L^p}+1\preceq \|u\|^2_{\mH^p_1}+1.\label{YY0}
\ee
Thus, as above,
by (iii) of Proposition \ref{Pr1} and H\"older's inequality we have, for $q>\frac{2}{1-\a}$
\ce
&&\mE\left(\sup_{t'\in[0,\xi^\eps_n(t)]}\|\cJ_2(t)\|^q_{\mH^p_{2\a}}\right)\\
&&\qquad\leq\mE\left(\sup_{t'\in[0,\xi^\eps_n(t)]}\int^t_0\|\fT_{t-s}\psi(s,\cdot,u(s))\dot h^\eps(s)\|_{\mH^p_{2\a}}
\dif s\right)^q\\
&&\qquad\leq N\mE\left(\sup_{t'\in[0,\xi^\eps_n(t)]}\int^t_0\|\fT_{t-s}\psi(s,\cdot,u(s))\|^2_{L_2(l^2;\mH^p_{2\a})}
\dif s\right)^{q/2}\\
&&\qquad \leq C_N\mE\left(\sup_{t'\in[0,\xi^\eps_n(t)]}\int^t_0\frac{1}{(t-s)^\a}
\|\psi(s,\cdot,u(s))\|^2_{L_2(l^2;\mH^p_{\a})}\dif s\right)^{q/2}\\
&&\qquad\stackrel{(\ref{YY0})}{\leq}C_{T,N}\mE\left(\int^{\xi^\eps_n(t)}_0
(\|u(s)\|^q_{\mH^p_1}+1)\dif s\right)\\
&&\qquad\leq C_{T,N} \int^t_0(f_2(s)+1)\dif s.
\de

Set
$$
G(t,s):=\sqrt{\eps}\fT_{t-s}\psi(s,\cdot,u(s)).
$$
Then  by (iii) and (iv) of Proposition \ref{Pr1} we have
$$
\|G(t,s)\|^2_{\mH^p_{2\a}}\leq\frac{C}{(t-s)^\a}\|\psi(s,\cdot,u(s))\|^2_{L_2(l^2;\mH^p_\a)}
$$
and for $\g\in(0,(1-\a)/2)$
$$
\|G(t',s)-G(t,s)\|^2_{\mH^p_{2\a}}\leq\frac{|t'-t|^\g}{(t-s)^{\a+2\g}}
\|\psi(s,\cdot,u(s))\|^2_{L_2(l^2;\mH^p_\a)}.
$$
Therefore, using Lemma \ref{Le00}, for $q$ large enough, we have
\ce
\mE\left(\sup_{t'\in[0,\xi^\eps_n(t)]}\|\cJ_3(t)\|^q_{\mH^p_{2\a}}\right)
&=&\mE\left(\sup_{t'\in[0,T\wedge\xi_n^\eps(t)]}\left\|\int^{t'}_0G(t',s)\dif W(s)
\right\|^q_{\mH^p_{2\a}}\right)\\
&\leq& C_T\mE\left(\int^{T\wedge\xi_n^\eps(t)}_0
\|\psi(s,\cdot,u(s))\|_{L_2(l^2;\mH^p_\a)}^q\dif s\right)\\
&\stackrel{(\ref{YY0})}{\leq}& C_T \int^t_0(f_2(s)+1)\dif s.
\de
Combining the above estimates, we get
$$
f_2(t)\leq C\|u_0\|^q_{\mH^p_2}+C_{T,N} \int^t_0(f_2(s)+1)\dif s.
$$
By Gronwall's inequality again, we find
\ce
\mE\left(\sup_{t\leq T\wedge\tau^\eps_n}\|u(t)\|^q_{\mH^p_{2\a}}\right)=f_2(T)\leq C_{T,N},
\de
which in turn implies that
$$
\lim_{n\to\infty}\sup_{\eps\in(0,1)}P\{\tau^\eps_n<T\}=0.
$$
The proof is complete.
\end{proof}
Moreover, under {\bf (R3)}, similar to the above lemma, we can check that (\ref{Au2}) holds.
Thus, using Theorem \ref{Th6} we obtain
\bt
Let $(M,\fg)$ be a compact Riemannian manifold, and $p>d$,
$\a$ satisfy (\ref{Con00}). Let $u_0\in\mH_2^p$. Under {\bf (R1)}-{\bf (R3)},
$\{u_\eps,\eps\in(0,1)\}$ satisfies the large deviation principle in $\mC_T(\mH_{2\a}^p)$
with the rate function $I(f)$ given by
$$
I(f):=\frac{1}{2}\inf_{\{h\in\ell^2_T:~f=u^h\}}\|h\|^2_{\ell^2_1},\ \ f\in\mC_T(\mH^p_{2\a}),
$$
where  $u^h$ solves the following equation:
\ce
u^h(t)&=&u_0+\int^{t}_0\big[\Delta u^h(s)+\varphi(s,\cdot,u^h(s),
\fg(Y(\cdot),\nabla u^h(s)))\big]\dif s\\
&&+\int^{t}_0\psi(s,\cdot,u^h(s))\dot h(s)\dif s.
\de
\et

\section{Application to stochastic Navier-Stokes equations}

\subsection{Unique maximal strong solution for SNSEs}

Let $\cO$ be a bounded smooth domain in $\mR^d$($d\geq 2$), or the whole space
$\mR^d$, or $d$-dimensional torus $\mT^d$. Let
$$
\W^{m,p}(\cO):=(W^{m,p}(\cO))^d,\ \ \W^{m,p}_0(\cO):=(W^{m,p}_0(\cO))^d
$$
and
$$
\C^\infty_{0,\sigma}(\cO):=\{\u\in (C^\infty_0(\cO))^d: \div(\u)=0\}.
$$
Notice that $\W^{m,p}(\mR^d)=\W^{m,p}_0(\mR^d)$ and  $\W^{m,p}(\mT^d)=\W^{m,p}_0(\mT^d)$.

Let $\L^p_\sigma(\cO)$ be the closure of $\C^\infty_{0,\sigma}(\cO)$ with respect to
the norm in $\L^p(\cO):=(L^p(\cO))^d$.
Let $\sP_2$ be the orthonormal projection from $\L^2(\cO)$ to $\L^2_\sigma(\cO)$. It is well known that
$\sP_2$ can be extended to a bounded linear operator from $\L^p(\cO)$ to $\L^p_\sigma(\cO)$
 (cf. \cite{Fu-Mo}) so that  for every $\u\in \L^p(\cO)$
$$
\u=\sP_p\u+\nabla\pi,\ \ \pi\in (L^p_{loc}(\cO))^d.
$$

The stokes operator is defined by
\be
A_p\u:=-\sP_p\Delta\u,\ \ \sD(A_p):=\mH^p_2\cap\L^p_\sigma(\cO),\label{Es9}
\ee
where
$$
\mH^p_2:=\W^{2,p}(\cO)\cap\W^{1,p}_0(\cO)=\sD(I-\Delta_p)
$$
and $\Delta_p$ is the Laplace operator on $\L^p(\cO)$.

It is well known that $(A_p,\sD(A_p))$ is a sectorial operator on $\L^p_\sigma(\cO)$  (cf. \cite{Gi1}).
It should be noticed that when $\cO=\mR^d$ or $\mT^d$, since the projection $\sP_p$ can
commute with $\nabla$  (cf. \cite[p.84]{Li}), we have
$$
A_p\u=-\Delta\sP_p\u=-\Delta\u,\ \ \u\in\sD(A_p).
$$
That is, the stokes operator is just the restriction of $-\Delta_p$ on
$\W^{2,p}(\cO)\cap \L^p_\sigma(\cO)$, where $\cO=\mR^d$ or $\mT^d$.

Below, we write
$$
\fL_p:=I+A_p
$$
and
$$
\bH_\a^p:=\sD(\fL^{\a/2}_p).
$$
Giga \cite{Gi2} proved that for $\a\in[0,1]$
\be
\bH_\a^p=[\L^p_\sigma(\cO),\sD(A_p)]_\a=\mH^p_\a\cap\L^p_\sigma(\cO),\label{Do}
\ee
where $\mH^p_\a=[\L^p(\cO),\mH^p_2]_\a$ and $[\cdot,\cdot]_\a$ stands for the complex interpolation
space between two Banach spaces. In particular, the following embedding results hold (see (\ref{PP4}) and (\ref{PL2})):
for $p>1$ and $0\leq\a'<\frac{1}{2}<\a\leq 1$
\be
\|\u\|_{\bH^p_{2\a'}}\preceq\|\u\|_{1,p}\preceq \|\u\|_{\bH^p_{2\a}}, \ \ \u\in\bH^p_\a,
\label{PP04}
\ee
and for $q\geq p$, $k-\frac{d}{q}<2\a-\frac{d}{p}$
\be
\bH^p_{2\a}\hookrightarrow \W^{k,q}(\cO),\label{Em}
\ee
and for $\a>\frac{d}{p}$
\be
\bH_\a^p\hookrightarrow C_b(\cO).\label{PL02}
\ee

In what follows, we fix
\be
p>d,\ \ \
\frac{1}{2}<\a<1,\label{Con2}
\ee
and consider the following stochastic Navier-Stokes equation with Dirichlet boundary
(only for bounded smooth domain):
\be\label{NS}
\left\{
\begin{aligned}
&\dif\u(t)=\big[\Delta\u(t)+(\u(t)\cdot\nabla)\u(t)+\nabla \pi(t)\big]\dif t\\
&\qquad\qquad+F(t,\u(t))\dif t+\Psi(t,\u(t))\dif W(t)\\
&\u(t,\cdot)|_{\p\cO}=0,\ \ \ \ \ \ \div\u(t)=0,\\
&\u(0,x)=\u_0(x),
\end{aligned}
\right.
\ee
where $\u$ and $\pi$ are unknown functions, and
$$
F: \mR_+\times\bH^p_{2\a}\to\bH^p_0\mbox{ and }\Psi:\mR_+\times\bH^p_{2\a}\to\bH^p_{\a}
$$
are two measurable functions.

We assume that
\begin{enumerate}[{\bf(N1)}]
\item For each $T,R>0$, there exist  $\d>0$ and $C_{T,R,\d}>0$
such that for all $t,s\in[0,T]$ and $\u,\v\in\bH^p_{2\a}$ with
$\|\u\|_{\bH^p_{2\a}},\|\v\|_{\bH^p_{2\a}}\leq R$
$$
\|F(t,\u)-F(s,\v)\|_{\bH^p_0}\leq C_{T,R,\d}\Big(|t-s|^\d+\|\u-\v\|_{\bH^p_{2\a}}\Big).
$$

\item For each $T,R>0$, there exist  $\a'>1$ and $C_{T,R}>0$
such that for all $t\in[0,T]$ and $\u,\v\in\bH^p_{2\a}$ with
$\|\u\|_{\bH^p_{2\a}},\|\v\|_{\bH^p_{2\a}}\leq R$
$$
\|\Psi(t,\u)-\Psi(t,\v)\|_{L_2(l^2;\bH^p_{\a})}\leq C_{T,R}\|\u-\v\|_{\bH^p_{2\a}}
$$
and
\be
\|\Psi(t,\u)\|_{L_2(l^2;\bH^p_{\a'})}\leq C_{T,R}.\label{V0}
\ee
\end{enumerate}

Set
\be
\Phi(t,\u):=\u+\sP_p[(\u \cdot\nabla)\u] +F(t,\u).\label{Phi1}
\ee
Then \textsc{Eq.}(\ref{NS}) can be written as the following abstract form:
\be
\dif \u(t)=[-\fL_p\u(t)+\Phi(t,\u)]\dif t+\Psi(t,\u)\dif W(s),\ \ \u(0)=\u_0.\label{NS1}
\ee

\bt\label{THNS} Let $p>d$ and $\frac{1}{2}<\a<1$.
Under {\bf (N1)} and {\bf (N2)}, for any $\u_0\in\bH^p_2$,
there exists a unique maximal strong solution $(\u,\tau)$ for \textsc{Eq.}(\ref{NS1}) so that
\begin{enumerate}[(i)]
\item $t\mapsto\u(t)\in\bH^p_2$ is continuous on $[0,\tau)$ a.s.;
\item $\lim_{t\uparrow\tau}\|\u(t)\|_{\bH^p_{2\a}}=\infty$ on $\{\tau<+\infty\}$;
\item it holds that in $\L^p_\sigma(\cO)=\bH^p_0$
\ce
\u(t)&=&\u_0+\int^t_0[-\fL_p\u(s)+\Phi(s,\u(s))]\dif s+\int^t_0\Psi(s,\u(s))\dif W(s)\\
&=&\u_0+\int^t_0[A_p\u(s)+\sP_p((\u(s)\cdot\nabla)\u(s))]\dif s\\
&&+\int^t_0F(s,\u(s))\dif s+\int^t_0\Psi(s,\u(s))\dif W(s),
\de
for all  $t\in[0,\tau)$, $P$-a.s..
\end{enumerate}
\et
\begin{proof}
In view of (\ref{Con2}), (\ref{PP04}) and (\ref{PL02}), for any $\u,\v\in\bH^p_{2\a}$ we have

\ce
\|\sP_p[(\u \cdot\nabla)\u-(\v \cdot\nabla)\v]\|_{\L^p_\sigma}&\preceq&
\|(\u \cdot\nabla)\u-(\v \cdot\nabla)\v\|_{\L^p}\\
&\preceq&\|\u-\v\|_{\L^\infty} \cdot\|\nabla\u\|_{\L^p}\\
&&+\|\v\|_{\L^\infty}\cdot\|\nabla(\u-\v)\|_{\L^p}\\
&\preceq&\|\u-\v\|_{\bH^p_{2\a}} \cdot\|\u\|_{\bH^p_{2\a}}\\
&&+\|\v\|_{\bH^p_{2\a}}\cdot\|\u-\v\|_{\bH^p_{2\a}},
\de
Thus, by {\bf (N1)} and {\bf (N2)}, it is easy to see
that {\bf (M2)} and {\bf (M4)} hold for the above $\Phi$ and $\Psi.$
The result now follows by  Theorem \ref{Th55}.
\end{proof}

We now give two concrete functionals so that {\bf (N1)} and {\bf (N2)} are satisfied.
Let $\f:\mR_+\times\cO\times\mR^d\to\mR^d$ be a measurable function, and satisfy that:
for any $T,R>0$, there exist constants $\d, C_{T,R}>0$ and $\lambda^{\f}_{R,T}\in L^p(\cO)$
such that for all $t,s\in[0,T], x\in\cO$ and $\u,\v\in\mR^d$ with $|\u|,|\v|\leq R$
$$
|\f(t,x,\u)-\f(s,x,\v)|\leq C_{T,R}\big(\lambda^{\f}_{R,T}(x)\cdot|t-s|^\d+|\u-\v|\big).
$$
Let $\b:\mR_+\times\cO\times\mR^d\to l^2\times\mR^d$ be a measurable function, and satisfy that:
\ce
\left\{
\begin{aligned}
\left.
\begin{aligned}
&\b(t,x,\u)=c(t)\u+\b_2(t,x),\ \\
&\exists\a'>1\ \ s.t.\ \ \sup_{t\in[0,T]}(|c(t)|+\|\b_2(t,\cdot)\|_{\bH^p_{\a'}})\leq C_T,
\end{aligned}
\right\},&\mbox{ $\cO$  bounded;}\\
\left.
\begin{aligned}
&\b(t,\cdot,\cdot)\in C^2(\cO\times\mR^d;l^2),\mbox{ and for $|\u|\leq R$}, \\
&\|\nabla_x\p_{\u}\b(t,x,\u)\|_{l^2}+\|\p^j_{\u}\b(t,x,\u)\|_{l^2}\leq C_{R,T}, j=1,2,\\
&\sup_{t\in[0,T]}\|\nabla^j_x\b(t,\cdot,0)\|_{\mR^d\times l^2}\leq \lambda^{\b}_{R,T}(x),\ j=0,1,2
\end{aligned}
\right\},&\mbox{ $\cO=\mR^d$ or $\mT^d$},
\end{aligned}
\right.
\de
where $\lambda^{\b}_{R,T}\in L^p(\cO)$.

We define
\be
F(t,\u):=\sP_p(\f(t,\cdot,\u))\label{FE}
\ee
and
\be
\Psi(t,\u):=\sP_p(\b(t,\cdot,\u)).\label{FE1}
\ee
One can see that {\bf (N1)} and {\bf (N2)} hold.
Indeed, for $\u,\v\in\bH^p_{2\a}$ with
$\|\u\|_{\bH^p_{2\a}},\|\v\|_{\bH^p_{2\a}}\leq R$, we have
\ce
\|F(t,\u)-F(s,\v)\|_{\L^p_\sigma}
&\preceq&\|\f(t,\cdot,\u)-\f(s,\cdot,\v)\|_{\L^p}\\
&\stackrel{(\ref{PL02})}{\leq}& C_{T,R}(|t-s|^\d+\|\u-\v\|_{\L^p})\\
&\leq& C_{T,R}(|t-s|^\d+\|\u-\v\|_{\bH^p_{2\a}}).
\de
Thus, {\bf (N1)} holds. For {\bf (N2)}, let us look at
the case of $\cO=\mR^d$ or $\mT^d$. Since $\Delta_p$ can commute with $\sP_p$,
we have, for $\u,\v\in\bH^p_{2\a}$ with
$\|\u\|_{\bH^p_{2\a}},\|\v\|_{\bH^p_{2\a}}\leq R$
\ce
\|\Psi(t,\u)-\Psi(t,\v)\|_{L_2(l^2;\bH^p_\a)}^2&\stackrel{(\ref{Exa})}{\preceq}&
\|\fL^{\frac{\a}{2}}\sP_p[\b(t,\u)-\b(t,\v)]\|_{\L^p_\sigma(\cO;l^2)}^2\\
&\preceq&\|(I-\Delta_p)^{\frac{\a}{2}}[\b(t,\u)-\b(t,\v)]\|_{\L^p(\cO;l^2)}^2\\
&\stackrel{(\ref{PP04})}{\preceq}&\sum_k\|\b_k(t,\u)-\b_k(t,\v)\|_{1,p}^2\\
&\preceq&C_R\|\u-\v\|^2_{\bH^p_{2\a}}.
\de
Using Lemma \ref{Le7}, as the calculations given in Theorem \ref{Th7}, one can verify that
(\ref{V0}) holds under (\ref{Con00}). Thus, {\bf (N2)} holds.

\subsection{Non-explosion and large deviation for 2D SNSEs}

In this subsection, we study the non-explosion and large
deviation for SNSE in the case of two dimension. For this aim, in addition to {\bf (N1)} and {\bf (N2)},
we also suppose that
\begin{enumerate}[{\bf (N3)}]
\item For any $T>0$, there exists  $C_T>0$ such that for all $t\in[0,T]$ and $\u\in\bH^p_2$
\ce
\|F(t,\u)\|_{\bH^2_0}&\leq& C_T(\|\u\|_{\bH^2_1}+1),\\
\|F(t,\u)\|_{\bH^p_0}&\leq& C_T(\|\u\|_{\bH^p_{2\a}}+1)
\de
and for $i=0,1$
\ce
\|\Psi(s,\u)\|_{L_2(l^2;\bH^2_i)}&\leq& C_T(1+\|\u\|_{\bH^2_i}), \\
\|\Psi(s,\u)\|_{L_2(l^2;\bH^p_\a)}&\leq& C_T(1+\|\u\|_{\bH^p_{2\a}}),
\de
where $p$ and $\a$ satisfy (\ref{Con2}).
\end{enumerate}
We remark that $F$ and $\Psi$ defined by (\ref{FE}) and (\ref{FE1})
satisfy {\bf (N3)} when $\f$  satisfies
$$
|\f(t,x,\u)|\leq C_T(|\u|+\l_0(x))
$$
and $\b$ satisfies ($\cO=\mR^2$ or $\mT^2$)
\ce
\|\p_{\u}\b(t,x,\u)\|_{l^2}&\leq& C_T,\\
\|\b(t,x,\u)\|_{l^2}+\|\nabla_x\b(t,x,\u)\|_{l^2}&\leq& C_T(|\u|+\l_1(x)),
\de
where $\l_0,\l_1\in L^p(\cO)$.

We have the following result, the proof will be given in Lemma \ref{Le9} below.
\bt Let $p>d$ and $\frac{1}{2}<\a<1$.
Assume that {\bf (N1)}-{\bf (N3)} hold. Let $(\u,\tau)$ be the
unique maximal solution of \textsc{Eq.}(\ref{NS2}) in Theorem \ref{THNS}.
Then $\tau=+\infty$ a.s..
\et

We now consider the small perturbation for 2D stochastic Navier-Stokes equation:
\ce
\dif \u_\eps(t)=\big[-\fL_p\u_\eps(t)+\Phi(t,\u_\eps(t))\big]\dif t
+\sqrt{\eps}\Psi(t,\u_\eps(t))\dif W(t),\ \ \u_\eps(0)=\u_0
\de
as well as the control equation:
\be
\dif \u^\eps(t)&=&\big[-\fL_p\u^\eps(t)+\Phi(t,\u^\eps(t))
+\Psi(t,\u^\eps(t))\dot h^\eps(t)\big]\dif t\no\\
&&+\sqrt{\eps}\Psi(t,\u^\eps(t))\dif W(t),\ \ \u^\eps(0)=\u_0,\label{NS2}
\ee
where $h^\eps\in\cA^{T}_N$ (see (\ref{Op2}) for the definition of $\cA^{T}_N$),
and $T>0$ is fixed below.

Let $(\u^\eps,\tau^\eps)$ be the unique maximal strong solution of \textsc{Eq.} (\ref{NS2})
with the properties:
$$
\lim_{t\uparrow\tau^\eps}\|\u^\eps(t)\|_{\bH^p_{2\a}}=+\infty
\mbox{ on $\{\tau^\eps<\infty\}$},
$$
and  $t\mapsto\u^\eps(t)\in\bH^p_2$ is continuous on $[0,\tau^\eps)$.

Before proving the non-explosion result (Lemma \ref{Le9}),
we first prepare a series of lemmas.
\bl
There exists a  constant $C_T>0$ such that for any $t\in[0,T]$ and $\u\in\bH^2_2$
\be
\<\u,-\fL_2\u+\Phi(s,\u)\>_{\bH^2_0}&\leq& -\frac{1}{2}\|\u\|_{\bH^2_1}^2
+C_T(\|\u\|^2_{\bH^2_0}+1),\label{LP1}\\
\<\fL_2\u,-\fL_2\u+\Phi(s,\u)\>_{\bH^2_0}&\leq& C\|\u\|^2_{\bH^2_0}\|\u\|^4_{\bH^2_1}
+C_T\big(1+\|\u\|^2_{\bH^2_1}\big) \label{LP2}
\ee
and
\be
\|\Phi(t,\u)\|_{\bH^p_0}\leq C_T\big(1+\|\u\|_{\bH^2_1}\big)\cdot\big(1+\|\u\|_{\bH^p_{2\a}}\big).
\label{LL4}
\ee
\el
\begin{proof}
Let $\u\in\bH^2_2$. Noting that
$$
\<\u,\sP_2((\u\cdot\nabla)\u)\>_{\bH^2_0}=
\<\u,(\u\cdot\nabla)\u\>_{\L^2}=\frac{1}{2}\int_\cO\u(x)\cdot\nabla|\u(x)|^2\dif x=0,
$$
by {\bf (N3)} and Young's inequality we have
\ce
\<\u,-\fL_2\u+\Phi(s,\u)\>_{\bH^2_0}
&=&-\|\u\|^2_{\bH^2_1}+\<\u,\u+F(t,\u))\>_{\bH^2_0}\\
&\leq&-\frac{1}{2}\|\u\|^2_{\bH^2_1}+C_T(\|\u\|^2_{\bH^2_0}+1).
\de
Thus, (\ref{LP1}) is proved.

For (\ref{LP2}), noting that by Gagliado-Nirenberge's inequality
 (cf. \cite[p.24 Theoerem 9.3]{Fr}) and (\ref{Do})
$$
\|\u\|^2_{\L^\infty}\preceq \|\u\|_{\mH^2_2}\cdot\|\u\|_{\mH^2_0}
\preceq \|\u\|_{\bH^2_2}\cdot\|\u\|_{\bH^2_0},
$$
by Young's inequality we have
\ce
\<\fL_2\u,\sP_2((\u\cdot\nabla)\u)\>_{\bH^2_0}
&\leq&\frac{1}{4}\|\u\|_{\bH^2_2}^2+\|\sP_2((\u\cdot\nabla)\u)\|^2_{\bH^2_0}\\
&\leq&\frac{1}{4}\|\u\|_{\bH^2_2}^2+C\|(\u\cdot\nabla)\u)\|^2_{\L^2}\\
&\leq&\frac{1}{4}\|\u\|_{\bH^2_2}^2+C\|\u\|_{\L^\infty}^2\cdot\|\nabla\u\|^2_{\L^2}\\
&\leq&\frac{1}{4}\|\u\|_{\bH^2_2}^2+C\|\u\|_{\bH^2_0}\cdot\|\u\|_{\bH^2_2}\cdot\|\u\|^2_{\bH^2_1}\\
&\leq&\frac{1}{2}\|\u\|_{\bH^2_2}^2+C\|\u\|^2_{\bH^2_0}\cdot\|\u\|_{\bH^2_1}^4
\de
and by {\bf (N3)}
$$
\<\fL_2\u,F(s,\u)\>_{\bH^2_0}\leq\frac{1}{2}\|\u\|_{\bH^2_2}^2
+C_T(1+\|\u\|^2_{\bH^2_1}).
$$
Thus, (\ref{LP2}) holds.

Let
$$
p<q<\frac{d}{1+\frac{d}{p}-2\a},\ \ q^*=\frac{qp}{q-p}.
$$
By H\"older's inequality we have
\ce
\|\sP_p(\u\cdot\nabla)\u\|_{\bH^p_0}&\preceq&\|\u\cdot\nabla\u\|_{\L^p}\\
&\preceq&\|\u\|_{\L^{q^*}}\cdot\|\nabla\u\|_{\L^q}\\
&\stackrel{(\ref{Em})}{\preceq}&\|\u\|_{\bH^2_1}\cdot\|\u\|_{\bH^p_{2\a}}.
\de
The estimate (\ref{LL4}) now follows by {\bf (N3)}.
\end{proof}

Below, set for $n\in\mN$
\ce
\tau^\eps_n:=\inf\Big\{t\geq 0: \|\u^\eps(t)\|_{\bH_{2\a}^p}>n\Big\}.
\de
\bl\label{LL1}
There exists a  constant $C_T>0$ such that for all $\eps\in(0,1)$ and $n\in\mN$
\ce
\mE\left(\sup_{s\in[0,T\wedge\tau^\eps_n]}\|\u^\eps(s)\|^2_{\bH^2_0}\right)+
\mE\left(\int^{T\wedge\tau^\eps_n}_0
\|\u^\eps(s)\|^2_{\bH^2_1}\dif s\right)\leq C_T.
\de
\el
\begin{proof}
By Ito's formula we have
\ce
\|\u^\eps(t)\|_{\bH^2_0}^2&=&\|\u_0\|_{\bH^2_0}^2
+2\int^t_0\<\u^\eps(s),-\fL_2\u^\eps(s)+\Phi(s,\u^\eps(s))\>_{\bH^2_0}\dif s\\
&&+2\int^t_0\<\u^\eps(s),\Psi(s,\u^\eps(s))\dot h^\eps(s)\>_{\bH^2_0}\dif s\\
&&+2\sqrt{\eps}\sum_k\int^t_0\<\u^\eps(s),\Psi_k(s,\u^\eps(s))\>_{\bH^2_0}\dif W^k(s)\\
&&+\eps\sum_k\int^t_0\|\Psi_k(s,\u^\eps(s))\|^2_{\bH^2_0}\dif s\\
&=:&\|\u_0\|_{\bH^2_1}^2+J_1(t)+J_2(t)+J_3(t)+J_4(t).
\de

Set
$$
f(t):=\mE\left(\sup_{s\in[0,t\wedge\tau^\eps_n]}\|\u^\eps(s)\|^2_{\bH^2_0}\right).
$$
First of all, noting that by (\ref{LP1})
$$
J_1(t)\leq -\int^t_0\|\u^\eps(s)\|_{\bH^2_1}^2+C_T\int^t_0(\|\u^\eps(s)\|^2_{\bH^2_0}+1)\dif s,
$$
we have
\ce
\mE\left(\sup_{s\in[0,t\wedge\tau^\eps_n]}J_1(s)\right)+
\mE\left(\int^{t\wedge\tau^\eps_n}_0
\|\u^\eps(s)\|^2_{\bH^2_1}\dif s\right)\leq C_T\int^t_0(f(s)+1)\dif s.
\de
By {\bf (N3)} and Young's inequality we have
\ce
\mE\left(\sup_{s\in[0,t\wedge\tau^\eps_n]}J_2(s)\right)&\leq&2\mE\left(\int^{t\wedge\tau^\eps_n}_0
\|\u^\eps(s)\|_{\bH^2_0}\cdot\|\Psi(s,\u^\eps(s))\|_{L_2(l^2;\bH^2_0)}
\cdot\|\dot h^\eps(s)\|_{l^2}\dif s\right)\\
&\leq&2N\mE\left(\int^{t\wedge\tau^\eps_n}_0\|\u^\eps(s)\|^2_{\bH^2_0}\cdot
\|\Psi(s,\u^\eps(s))\|^2_{L_2(l^2;\bH^2_0)}\dif s\right)^{1/2}\\
&\leq&\frac{1}{4}f(t)+C_N\mE\left(\int^{t\wedge\tau^\eps_n}_0(1+\|\u^\eps(s)\|^2_{\bH^2_0})
\dif s\right)\\
&\leq&\frac{1}{4}f(t)+C_N\int^t_0(1+f(s))\dif s.
\de
Similarly, we also have
$$
\mE\left(\sup_{s\in[0,t\wedge\tau^\eps_n]}J_3(s)\right)\leq\frac{1}{4}f(t)+C\int^t_0(1+f(s))\dif s
$$
and
$$
\mE\left(\sup_{s\in[0,t\wedge\tau^\eps_n]}J_4(s)\right)\leq C\int^t_0(1+f(s))\dif s.
$$

Combining the above calculations we get
$$
f(t)+2\mE\int^{t\wedge\tau^\eps_n}_0\|\u^\eps(s)\|^2_{\bH^2_1}\dif s
\leq 2\|\u_0\|_{\bH^2_0}^2+C_N+C_N\int^t_0(1+f(s))\dif s.
$$
The desired estimate follows by Gronwall's inequality.
\end{proof}

Set for $n\in\mN$
\ce
\eta^\eps_n(t)&:=&\int^{t\wedge\tau^\eps_n}_0\|\u^\eps(s)\|^2_{\bH^2_1}\cdot
\|\u^\eps(s)\|^2_{\bH^2_0}\dif s+t\\
&=&\int^{t}_0\|\u^\eps(s)\|^2_{\bH^2_1}\cdot
\|\u^\eps(s)\|^2_{\bH^2_0}\cdot 1_{[0,\tau^\eps_n]}(s)\dif s+t
\de
and
$$
\theta^{\eps}_n(t):=\inf\left\{s\geq 0: \eta^\eps_n(s)\geq t\right\}.
$$
Clearly, $t\mapsto\eta^\eps_n(t)$ is a continuous and strictly increasing function, and
the inverse function of $t\mapsto\theta^\eps_n(t)$ is just given by $\eta^\eps_n$.
Moreover, since $\eta^\eps_n(t)>t$, we have
$$
\theta^{\eps}_n(t)<t.
$$

\bl\label{LL2}
For any $K>0$, there exists a  constant $C_{K,N}>0$ such that for all $\eps\in(0,1)$ and $n\in\mN$
$$
\mE\left(\sup_{s\in[0,\theta^\eps_n(K)\wedge\tau^{\eps}_n]}
\|\u^\eps(s))\|^2_{\bH^2_1}\right)\leq C_{K,N}.
$$
\el
\begin{proof}
Consider the following evolution triple
$$
\bH^2_2\subset\bH^2_1\subset\bH^2_0.
$$
By Ito's formula  (cf. \cite{Ro}), we have
\ce
\|\u^\eps(t)\|_{\bH^2_1}^2&=&\|\u_0\|_{\bH^2_1}^2
+2\int^t_0\<\fL_2\u^\eps(s),-\fL_2\u^\eps(s)+\Phi(s,\u^\eps(s))\>_{\bH^2_0}\dif s\\
&&+2\int^t_0\<\fL_2\u^\eps(s),\Psi(s,\u^\eps(s))\dot h^\eps(s)\>_{\bH^2_0}\dif s\\
&&+2\sqrt{\eps}\sum_k\int^t_0\<\u^\eps(s),\Psi_k(s,\u^\eps(s))\>_{\bH^2_1}\dif W^k(s)\\
&&+\eps\sum_k\int^t_0\|\Psi_k(s,\u^\eps(s))\|^2_{\bH^2_1}\dif s\\
&=:&\|\u_0\|_{\bH^2_1}^2+J_1(t)+J_2(t)+J_3(t)+J_4(t).
\de

Set
\ce
f(t)&:=&\mE\left(\sup_{s\in[0,t]}
\|\u^\eps(\theta^\eps_n(s)\wedge\tau^{\eps}_n)\|^2_{\bH^2_1}\right)\\
&=&\mE\left(\sup_{s\in[0,\theta^\eps_n(t)\wedge\tau^{\eps}_n]}
\|\u^\eps(s)\|^2_{\bH^2_1}\right).
\de
For $J_1(t)$, by (\ref{LP2}) we have, for $t\in[0,K]$
\ce
J_1(\theta^\eps_n(t)\wedge\tau^{\eps}_n)
&\leq&\int^{\theta^\eps_n(t)\wedge\tau^{\eps}_n}_0
\Big[C\|\u^\eps(s)\|^2_{\bH^2_0}\cdot\|\u^\eps(s)\|_{\bH^2_1}^4
+C_K(1+\|\u^\eps(s)\|_{\bH^2_1}^2)\Big]\dif s\\
&\leq& C\int^{\theta^\eps_n(t)}_0\|\u^\eps(s\wedge\tau^{\eps}_n)\|_{\bH^2_1}^2\dif
\eta_n^\eps(s)+C_K\\
&=&C\int^{t}_0\|\u^\eps(\theta^\eps_n(s)\wedge\tau^{\eps}_n)\|_{\bH^2_1}^2\dif s+C_K,
\de
where the last step is due to the substitution of variable formula .
So,
$$
\mE\left(\sup_{s\in[0,t]}J_1(\theta^\eps_n(s)\wedge\tau^{\eps}_n)\right)
\leq C\int^{t}_0f(s)\dif s+C_K.
$$
Using the same trick as used in Lemma \ref{LL1} and by {\bf (N3)}, we also have
$$
\mE\left(\sup_{s\in[0,t]}J_i(\theta^\eps_n(s)\wedge\tau^{\eps}_n)\right)
\leq \frac{1}{2}f(t)+C_{N,K}\int^{t}_0(f(s)+1)\dif s,\ \ \ i=2,3,4.
$$
Thus, we get
$$
f(t)\leq 2\|\u_0\|_{\bH^2_1}^2+C_{N,K}\int^{t}_0(f(s)+1)\dif s,
$$
which yields the desired estimate by Gronwall's inequality.
\end{proof}

Set for $M>0$
$$
\zeta^\eps_n(M):=\inf\Big\{t\geq 0:\|\u^\eps(t\wedge\tau^\eps_n)\|_{\bH^2_1}\geq M\Big\}.
$$
\bl \label{Le99}
For any $M>0$ and $q\geq 2$, there exists a  constant $C_{T,M,N}>0$ such that
for all $\eps\in(0,1)$ and $n\in\mN$
$$
\mE\left[\sup_{t\in[0,T\wedge\tau^\eps_n\wedge\zeta^\eps_n(M)]}
\|\u^\eps(t)\|_{\bH_{2\a}^p}^q\right]\leq C_{T,M,N}.
$$
\el
\begin{proof}
Set for $t\in[0,T]$
$$
\xi_n^\eps(t):=t\wedge\tau^\eps_n\wedge\zeta^\eps_n(M)
$$
and for $q\geq 2$
$$
f(t):=\mE\left[\sup_{t'\in[0,\xi^\eps_n(t)]}\|\u^\eps(t)\|_{\bH_{2\a}^p}^q\right].
$$
Note that
\ce
\u^\eps(t)&=&\fT_t\u_0+\int^t_0\fT_{t-s}\Phi(s,\u^\eps(s))\dif s+
\int^t_0\fT_{t-s}\Psi(s,\u^\eps(s))\dot h^\eps(s)\dif s\\
&&+\sqrt{\eps}\int^t_0\fT_{t-s}\Psi(s,\u^\eps(s))\dif W(s).
\de
By (iii) of Proposition \ref{Pr1}, H\"older's inequality
and Lemma \ref{LL4}, we have, for $q>\frac{1}{1-\a}$
\ce
&&\mE\left[\sup_{t'\in[0,\xi^\eps_n(t)]}\left
\|\int^{t'}_0\fT_{t'-s}\Phi(s,\u^\eps(s))\dif s\right\|^q_{\bH^p_{2\a}}\right]\\
&&\qquad\preceq\mE\left[\sup_{t'\in[0,\xi^\eps_n(t)]}
\left(\int^{t'}_0\frac{1}{(t'-s)^\a}\|\Phi(s,\u^\eps(s))\|_{\bH_0^p}\dif s\right)^q\right]\\
&&\qquad\preceq\mE\left[\int^{\xi^\eps_n(t)}_0
\|\Phi(s,\u^\eps(s))\|_{\bH_0^p}^q\dif s\right]  \\
&&\qquad\stackrel{(\ref{LL4})}{\preceq}\mE\left[\int^{\xi^\eps_n(t)}_0
\Big[\big(1+\|\u^\eps(s)\|_{\bH^2_1}^q\big)\cdot
\big(1+\|\u^\eps(s)\|^q_{\bH^p_{2\a}}\big)\Big]\dif s\right]\\
&&\qquad\leq C_{M}\int^t_0(f(s)+1)\dif s.
\de
On the other hand, set
$$
G(t,s):=\fT_{t-s}\Psi(s,\u^\eps(s)).
$$
Then by (iii) and (iv) of Proposition \ref{Pr1}, we have
$$
\|G(t,s)\|^2_{\bH^p_{2\a}}\leq\frac{C}{(t-s)^\a}\|\Psi(s,\u^\eps(s))\|^2_{L_2(l^2;\bH^p_\a)}
$$
and for $\g\in(0,(1-\a)/2)$
$$
\|G(t',s)-G(t,s)\|^2_{\bH^p_{2\a}}\leq\frac{|t'-t|^\g}{(t-s)^{\a+2\g}}
\|\Psi(s,\u^\eps(s))\|^2_{L_2(l^2;\bH^p_\a)}.
$$
Therefore, using Lemma \ref{Le00} for $q$ large enough, we get
\ce
&&\mE\left(\sup_{t'\in[0,T\wedge\xi_n^\eps(t)]}\left\|\int^{t'}_0G(t',s)\dif W(s)
\right\|^q_{\bH^p_{2\a}}\right)\\
&&\qquad\leq C_T\mE\left(\int^{T\wedge\xi_n^\eps(t)}_0
\|\Psi(s,\u^\eps(s))\|_{L_2(l^2;\bH^p_\a)}^q\dif s\right)\\
&&\qquad\stackrel{{\bf (N3)}}{\leq}C_T \int^t_0(f(s)+1)\dif s.
\de
Similarly, we have
\ce
&&\mE\left(\sup_{t'\in[0,T\wedge\xi_n^\eps(t)]}\left\|
\int^t_0\fT_{t-s}\Psi(s,\u^\eps(s))\dot h^\eps(s)\dif s\right\|^q_{\bH^p_{2\a}}\right)\\
&&\qquad\leq C_{T,N} \int^t_0(f(s)+1)\dif s.
\de
Combining the above calculations, we obtain
$$
f(t)\leq C_{T,M,N}\int^t_0f(s)\dif s+C_{T,M,N},
$$
which yields the desired estimate by Gronwall's inequality.
\end{proof}

\bl \label{Le9}
It holds that
\be
\lim_{n\rightarrow\infty}\sup_{\eps\in(0,1)}
P\Big\{\omega:\tau^\eps_n(\omega)\leq T\Big\}=0.\label{Auu1}
\ee
\el
\begin{proof}
First of all, for any $M, K>0$ we have
\ce
P\{\zeta^\eps_n(M)<T\}&\leq&P\{\zeta^\eps_n(M)<T;\theta^\eps_n(K)\geq T\}
+P\{\theta^\eps_n(K)< T\}\\
&=&P\left\{\sup_{t\in[0,T)}\|\u^\eps(t\wedge
\tau^\eps_n)\|_{\bH^2_1}>M;\theta^\eps_n(K)\geq T\right\}\\
&&+P\left\{\sup_{s\in[0,T)}\eta^\eps_n(s)>K\right\}\\
&\leq&P\left\{\sup_{t\in[0,\theta^\eps_n(K)\wedge
\tau^\eps_n]}\|\u^\eps(t)\|_{\bH^2_1}>M\right\}
+P\big\{\eta^\eps_n(T)>K\big\}\\
&\leq&\mE\left(\sup_{t\in[0,\theta^\eps_n(K)\wedge
\tau^\eps_n)}\|\u^\eps(t)\|_{\bH^2_1}^2\right)/M^2
+\mE\left(\eta^\eps_n(T)\right)/K.
\de
Hence, by Lemmas \ref{LL1} and \ref{LL2} we have
$$
\lim_{M\to\infty}\sup_{n,\eps}P\{\zeta^\eps_n(M)<T\}=0.
$$

Secondly, we also have
\be
P\{\tau^\eps_n<T\}\leq P\{\tau^\eps_n<T;\zeta^\eps_n(M)\geq T\}
+P\{\zeta^\eps_n(M)<T\}.\label{LP3}
\ee
For the first term, by Lemma \ref{Le99} we have
\ce
P\{\tau^\eps_n<T;\zeta^\eps_n(M)\geq T\}&=& P\left\{\sup_{t\in[0,T)}\|\u^\eps(t)
\|_{\bH^p_{2\a}}>n;\zeta^\eps_n(M)\geq T\right\}\\
&\leq& P\left\{\sup_{t\in[0,T\wedge\tau^\eps_n]}\|\u^\eps(t)
\|_{\bH^p_{2\a}}\geq n;\zeta^\eps_n(M)\geq T\right\}\\
&\leq& P\left\{\sup_{s\in[0,T\wedge\zeta^\eps_n(M)\wedge\tau^\eps_n]}\|\u^\eps(t)
\|_{\bH^p_{2\a}}\geq n\right\}\\
&\leq& \mE\left(\sup_{s\in[0,T\wedge\zeta^\eps_n(M)\wedge\tau^\eps_n]}
\|\u^\eps(t)\|^q_{\bH^p_{2\a}}\right)/n^q\\
&\leq&\frac{C_{T,M,N}}{n^q},
\de
where $C_{T,M,N}$ is independent of $\eps$ and $n$. The desired limit  now follows by taking limits
for (\ref{LP3}), first $n\to\infty$, then $M\to\infty$.
\end{proof}

Thus, using Theorem \ref{Th6} we get:
\bt\label{Th66}
Let  $\cO=\mT^2$ or a bounded smooth domain in $\mR^2$. Under {\bf (N1)}-{\bf (N3)},
for $\u_0\in\bH^p_2$, $\{\u_\eps,\eps\in(0,1)\}$ satisfies the large deviation principle
in $\mC_T(\bH_{2\a}^p)$ with the rate function $I(f)$ given by
$$
I(f):=\frac{1}{2}\inf_{\{h\in\ell^2_T:~f=\u^h\}}\|h\|^2_{\ell^2_T},\ \ f\in\mC_T(\bH_{2\a}^p),
$$
where $\u^h$ solves the following equation:
\ce
\u^h(t)&=&\u_0+\int^t_0\Delta\u^h(s)\dif s+\int^t_0\sP_p((\u^h(s)\cdot\nabla)\u^h(s))\dif s\\
&&+\int^t_0F(s,\u^h(s))\dif s+\int^t_0\Psi(s,\u^h(s))\dot h(s)\dif s.
\de
\et

\vspace{5mm}

{\bf Acknowledgements:}

The author would like to thank Professor Benjamin Goldys for
providing him an excellent environment to work in the University of New South Wales.
His work is supported by ARC Discovery grant DP0663153 of Australia and
NSF of China (No. 10871215).

\end{document}